\documentclass[11pt]{article}
\usepackage[]{amsmath,amssymb,amsfonts,latexsym,amsthm,enumerate,fullpage}
\usepackage[latin9]{inputenc}
\usepackage{amsmath}
\usepackage{amssymb}
\usepackage{breqn}
\usepackage[pdfstartview=FitH]{hyperref}
\usepackage{amsthm}
\usepackage{hyperref}
\usepackage{url}
\usepackage{enumitem}
\usepackage{authblk}
\usepackage{bbm}
\usepackage{verbatim}
\usepackage{cleveref}
\usepackage{tikz}
\usetikzlibrary{cd}
\usetikzlibrary{nfold}
\usetikzlibrary{arrows.meta, graphs, graphs.standard, quotes}
\usetikzlibrary{spath3}

\usepackage{todonotes}

\newcommand{\Aut}{\operatorname{Aut}}

\newcommand{\supp}{\operatorname{supp}}
\newcommand{\Cone}{\operatorname{Cone}}

\newcommand{\Rad}{\operatorname{Rad}}

\newcommand{\dist}{\operatorname{dist}}

\newcommand{\flag}{\operatorname{Flag}}
\newcommand{\orifl}{\operatorname{Orifl}}
\newcommand{\st}{\operatorname{st}}
\newcommand{\lk}{\operatorname{lk}}
\newcommand{\mat}{\operatorname{Mat}}

\newcommand{\cha}{\operatorname{char}}
\newcommand{\N}{\mathbb{N}}
\newcommand{\R}{\mathbb{R}}
\newcommand{\Z}{\mathbb{Z}}
\newcommand{\EE}{\mathcal{E}}

\renewcommand{\epsilon}{\varepsilon}
\newcommand*{\uak}{\mathcal{U}_A(K)}

\newcommand*{\CC}{\mathcal{CC}}

\newcommand*{\chev}{\operatorname{Chev}}

\newcommand{\Ao}{\mathring A}

\newcommand*{\cobound}{\operatorname{cb}}
\newcommand*{\cosys}{\operatorname{cs}}
\newcommand*{\SL}{\operatorname{SL}}
\newcommand*{\St}{\operatorname{St}}
\newcommand*{\height}{\operatorname{ht}}
\newcommand*{\abar}{\overline{\alpha}}

\newcommand{\ord}{\operatorname{ord}}

\newcommand{\diagnode}[1]{\fill #1 circle (.1);}

\DeclareMathOperator*{\freeprod}{\raisebox{-2pt}{\scalebox{2}{$\ast$}}}

\makeatletter

\begin{document}
\newtheorem{theorem}{Theorem}[section]
\newtheorem{lemma}[theorem]{Lemma}

\newtheorem{proposition}[theorem]{Proposition}
\newtheorem{corollary}[theorem]{Corollary}
\newtheorem{conjecture}[theorem]{Conjecture}

\theoremstyle{definition}
\newtheorem{definition}[theorem]{Definition}
\newtheorem{claim}[theorem]{Claim}
\newtheorem{example}[theorem]{Example}
\newtheorem{remark}[theorem]{Remark}
\newtheorem{observation}[theorem]{Observation}
\newtheorem{fact}[theorem]{Fact}

\newcommand{\subscript}[2]{$#1 _ #2$}

\definecolor{Mycolor1}{HTML}{082A54}
\definecolor{Mycolor2}{HTML}{E02B35}
\definecolor{Mycolor3}{HTML}{F0C571}
\definecolor{Mycolor4}{HTML}{59A89C}
\definecolor{Mycolor5}{HTML}{A559AA}

\renewcommand{\labelenumi}{(\roman{enumi})}

\author{Izhar Oppenheim  \footnote{The author is partially supported by ISF grant no.  242/24.}}
\affil{Department of Mathematics, Ben-Gurion University of the Negev, Be'er Sheva 84105, Israel, izharo@bgu.ac.il}

\author{Inga Valentiner-Branth \footnote{The author is supported by the FWO and the F.R.S.--FNRS under the Excellence of Science (EOS) program (project ID~40007542).}
}
\affil{Department of Mathematics, Computer science and Statistics, Ghent University, Krijgslaan 281 - S9, 9000 Ghent, Belgium, Inga.ValentinerBranth@ugent.be}

\title{Non-Abelian expansion of congruence KMS complexes}


\maketitle

\begin{abstract}
Coboundary expansion with non-Abelian coefficients is a strong version of high-dimensional expansion for simplicial complexes. One motivation for studying this notion is that it was recently shown to have deep connections to problems in theoretical computer science. However, very few examples of families of simplicial complexes with this type of expansion are known. Namely, prior to our work, the only known examples were quotients of symplectic buildings and a slight variation of the Kaufman-Oppenheim coset complexes construction associated with $\operatorname{SL}_{n} (\mathbb{F}_p [t])$.

In this paper, we show that the Grave de Peralta and Valentiner-Branth constructions of KMS complexes have coboundary expansion with non-Abelian coefficients when it is performed with respect to congruence subgroups of Chevalley groups of classical type, i.e., of type $A_n, B_n, C_n$ and $D_n$. This gives four new sources of examples to this expansion phenomenon, thus significantly enriching our list of constructions. 
\end{abstract}

\section{Introduction}
Given a simplicial complex $X$, the cohomology of $X$ is usually defined with coefficients in an Abelian group.  However,  when one considers the low dimensional cohomology,  the definition can be generalized to computing the $0$-th and first cohomologies of $X$ with respect to any coefficient group $\Lambda$ where $\Lambda$ need not be Abelian.  Given a group $\Lambda$,  $1$-coboundary expansion of $X$ is a quantization of the vanishing of $H^1 (X,  \Lambda)$ (see exact definition in \Cref{subsubsec: cosystolic/coboundary expansion subsec} below).  This notion of coboundary expansion with non-Abelian coefficients has received much attention lately due to its deep connections with theoretical computer science \cite{BM,  DDcovers, BLM,  DDL, BMV}  where the main focus is choosing the coefficient group $\Lambda$ to be a finite symmetric group.  

As in the case of expander graphs,  the focus is not on a single simplicial complex, but rather on a family of complexes.  Given a group $\Lambda$,  a family of simplicial complexes $\lbrace X^{(s)} \rbrace$ is said to be a family of $1$-coboundary expanders with $\Lambda$ coefficients if it has uniformly bounded degree,  the size of $X^{(s)}$ grows with $s$ to infinity and all the complexes $X^{(s)}$ have a uniform bound on their $\Lambda$ $1$-coboundary expansion (see exact definition in \Cref{subsubsec: cosystolic/coboundary expansion subsec} below).   

Prior to our work,  for a fixed finite symmetric group $\Lambda$,  only two constructions of families of $1$-coboundary expanders with $\Lambda$ coefficients were known: A construction coming from quotients of symplectic buildings (see \cite{BLM, DDL}) and a construction coming from a coset complexes over the group $\SL_{n+1} (\mathbb{F}_p [t])$ (see \cite{KO24cbeofcoco}).  In this work, we expand the list of known constructions and show how the constructions of certain types of KMS complexes that were introduced in \cite{hdxfromkms} give rise to $1$-coboundary expanders with $\Lambda$ coefficients.   

In \cite{hdxfromkms},  it was shown how to utilize KMS groups to construct high-dimensional spectral expanders.  As a specific instantiation of their general method,  the authors of \cite{hdxfromkms} showed how to construct a family of special KMS complexes, which we will call here \textit{congruence KMS complexes} (this is our terminology and it did not appear in \cite{hdxfromkms}).  The construction of congruence KMS complexes in \cite{hdxfromkms} can be roughly described as follows:  Given a spherical root system $\mathring \Phi$,  a finite field $K$ and an irreducible polynomial $f \in K [t]$,  the construction in \cite{hdxfromkms} outputs a finite simplicial complex $X=X (\mathring \Phi,  K,  f)$ that is a coset complex (see exact details in \Cref{subsubsec: Congruence KMS} and a first example in \Cref{A_n intro subsec} below). This construction is best behaved when $\cha(K) \neq 2$, hence we will assume is for the rest of this paper, although the assumption is only necessary for type $B_n$ and $C_n$ and even the characteristic two case could be handled with a bit more technical work, see \cite[Remark 5.10]{hdxfromkms}. 

With this notation,  our main result can be stated as follows:
\begin{theorem}[See Corollary \ref{cor: main result CBE} for a more general formulation]
Let $\Lambda$ be a group,  $n \geq 3$ an integer and $\mathring \Phi$ a spherical root system of type $A_n,  B_n,  C_n$ or $D_n$.  Also,  let $K$ be a finite field of characteristic $p \neq 2$ such that $\vert K \vert$ is sufficiently large with respect to $n$ and that $\Lambda$ does not have an element of order $p$ (e.g., if $\Lambda$ is finite and $p > \vert \Lambda \vert$).  Fix $f_s \in K[t]$ to be a sequence of irreducible polynomials such that $\lim_s \deg (f_s) = \infty$ and denote $X^{(s)} = X (\mathring \Phi,  K,  f_s)$.  Then the family $\lbrace X^{(s)} \rbrace_s$ is a family of $1$-coboundary expanders with $\Lambda$ coefficients.
\end{theorem}

This Theorem is deduced from two separate results that are  independent from each other:
\begin{itemize}
\item A vanishing result for the first non-Abelian cohomology for congruence KMS complexes.
\item A $1$-cosystolic expansion with non-Abelian coefficients for congruence KMS complexes.
\end{itemize}

\subsection{Vanishing of the first cohomology}

For congruence KMS complexes we prove the following vanishing result:
\begin{theorem}[See Corollary \ref{cor: homology of KMS complex is trivial}]
Let $\Lambda$ be a finite group,  $n \geq 3$ an integer and $\mathring \Phi$ a spherical root system.  Also,  let $K$ be a finite field of characteristic $p \neq 2$ such that $\vert K \vert \geq 5$ and such that $\Lambda$ does not have an element of order $p$ (e.g., if $\Lambda$ is finite and $p > \vert \Lambda \vert$).
 For any irreducible polynomial $f \in K[t]$ of degree at least $2$,  it holds that $H^1 (X (\mathring \Phi,  K,  f),  \Lambda) =0$.
\end{theorem}

It is worthwhile to compare our proof of the vanishing of cohomology with recent results in \cite{BLM, DDL} and in \cite{KO24cbeofcoco}.  For a fixed finite group $\Lambda$ all the works (\cite{BLM, DDL, KO24cbeofcoco} and ours) use the same basic methodology to construct a sequence of finite simplicial complexes with vanishing first $\Lambda$-cohomology.  We all start with  a universal ``algebraic'' simplicial complex that has trivial 1-cohomology with respect to $\Lambda$,  and then passes to a family of finite complexes by taking quotients via congruence subgroups, which also have trivial 1-cohomology with $\Lambda$-coefficients.

Applying this methodology yields different difficulties that are dependent on the specific construction:  
\begin{itemize}
\item In \cite{BLM, DDL}, the initial complex is a contractible symplectic $\widetilde{C}_n$-building, which has trivial 1-cohomology for any group. The difficulty lies in proving that the congruence subgroups exhibit trivial 1-cohomology with respect to $\Lambda$,  requiring deep results from group theory, such as the congruence subgroup property and the strong approximation theorem.
\item In \cite{KO24cbeofcoco},  the quotients are taken over a coset complex over $\SL_{n+1} (\mathbb{F}_p [t])$.  For $n \geq 3$,  there is an explicit description of the congruence subgroups, which follows from the fact that $\SL_{n+1} (\mathbb{F}_p [t])$ has the same presentation as  a Steinberg group.  This explicit description gives rise to a straightforward argument for proving the vanishing of 1-cohomology with $\Lambda$-coefficients given that $\Lambda$ has no element of order $p$ (e.g., when $\Lambda$ is finite and $p> \vert \Lambda \vert$).  Thus,  most of the work in \cite{KO24cbeofcoco} focuses on proving that the coset complex over $\SL_{n+1} (\mathbb{F}_p [t])$ has trivial 1-cohomology with respect to $\Lambda$.
\item At least with respect to proving vanishing of cohomology,  the congruence KMS complexes enjoy the best of both worlds: First,  the initial complex from which we take quotients is simply connected and thus has trivial 1-cohomology for any group. Second,  showing the vanishing of 1-cohomology for congruence subgroups with $\Lambda$-coefficients can be done via the presentation theory of the Chevalley groups (the argument shares some similarity with that of \cite{KO24cbeofcoco}).  The takeaway is that the proof of vanishing of the first cohomology with $\Lambda$ coefficients for congruence KMS complexes is perhaps the simplest in all known examples.   
\end{itemize}

\subsection{Non-Abelian cosystolic expansion}

Cosystolic expansion is a weaker variant of coboundary expansion.  There is a known local to global machinery for proving cosystolic expansion developed by \cite{KKL,  EK} for Abelian coefficients and \cite{KM,  DD-cosys} for non-Abelian coefficients.  Namely,  in order to show cosystolic expansion for a finite simplicial complex,  it is enough to show local spectral expansion and coboundary expansion of the links.  Local spectral expansion for congruence KMS complexes was shown in \cite{hdxfromkms} so we are left to prove non-Abelian coboundary expansion of the links.  

In \cite{CSEofKMS},  we already showed that links of congruence KMS complexes have Abelian coboundary expansion in all relevant dimensions.  Here, we generalize our method from \cite{CSEofKMS} to show that the links of congruence KMS complexes have $1$-coboundary expansion with respect to any group $\Lambda$:
\begin{theorem}[See Theorem \ref{thm: proper links are CBE} for a more general formulation]
Let $n \geq 3$ be an integer and $\mathring \Phi$ a spherical root system of type $A_n,  B_n,  C_n$ or $D_n$.  Also,  let $K$ be a finite field of characteristic $>2$ such that $\vert K \vert$ is sufficiently large with respect to $n$.  There is a constant $\varepsilon ' >0$ such that any group $\Lambda$,  any irreducible polynomial $f \in K [t]$ of degree at least $2$ and any vertex $v$ in $X (\mathring \Phi,  K,  f_s)$,  the $\Lambda$ $1$-coboundary expansion of $\lk (v)$ is $\geq \varepsilon '$.
\end{theorem}

Combining it with the results of \cite{DD-cosys} yields that congruence KMS complexes are $1$-cosystolic expanders:
\begin{theorem}[\Cref{cor: cosystolic exp of classcial KMS}]
Let $n \geq 3$ be an integer and $\mathring \Phi$ a spherical root system of type $A_n,  B_n,  C_n$ or $D_n$.  Also,  let $K$ be a finite field such that $\vert K \vert$ is sufficiently large with respect to $n$ and $\cha(K) \neq 2$.  There is a constant $\varepsilon >0$ such that for any group $\Lambda$ and any irreducible polynomial $f \in K [t]$ of degree at least $2$,   the $\Lambda$ $1$-cosystolic expansion of the complex $X (\mathring \Phi,  K,  f_s)$ is $\geq \varepsilon$.
\end{theorem}

\begin{remark}
Our results regarding $\Lambda$ $1$-cosystolic expansion are stated above for congruence KMS complexes,  but it is worth noting that they hold for a larger class of KMS complexes,  namely of $n$-classical KMS complexes (see \Cref{def: GCM and DyDi}).  As noted above,  the more general results are stated in \Cref{thm: proper links are CBE} and \Cref{cor: cosystolic exp of classcial KMS} below.
\end{remark}

\subsection{$A_n$ congruence KMS complexes}
\label{A_n intro subsec}

Our result stated above regards congruence KMS complexes of type $A_n,  B_n, C_n$ and $D_n$.  In order to give a concrete example,  here we will explicitly describe the congruence KMS complexes of type $A_n$ and state our result for these complexes. 

All congruence KMS complexes that we consider are coset complexes. More details on coset complexes can be found for example in \cite{KO2023high},  and we will only give a brief overview. 

\begin{definition}
Let $G$ be a group and $\mathcal{H}=(H_0,\dots , H_n)$ be a family of subgroups of $G$. Then the coset complex $\CC(G;(H_0, \dots, H_n))$ is defined to be the simplicial complex with 
\begin{itemize}
\item vertex set $\bigsqcup_{i=0}^n G/H_i$
\item a set of vertices $\left\{ g_1H_{i_1},\dots,g_kH_{i_k} \right\}$ forms a $(k-1)$-simplex in $\CC(G; \mathcal{H})$ if and only if $\bigcap_{j=1}^k g_j H_{i_j} \neq \emptyset$.
\end{itemize}
\end{definition}

Note that the set of vertices is partitioned into $n+1$ subsets $G/H_i$ for $i = 0,\dots,n$ such that for any simplex $\sigma \in \CC(G;\mathcal{H})$ we have $|\sigma \cap G/H_i| \leq 1$ for all $i$. 





With this notation,  we will define a coset complex of a subgroup of $\SL_{n+1} (\mathbb{F}_q [t])$ where $n \geq 2$ and $q$ is an odd prime power such that $q \geq 5$. For $1 \leq i,j \leq n+1,  i \neq j$ and $f \in  \mathbb{F}_q [t]$,  we denote $e_{i,j} (f) \in \SL_{n+1} (\mathbb{F}_q [t])$ to be the matrix with $1$'s along the main diagonal,  $f$ in the $(i,j)$-th entry and $0$'s in all other entries.  Let
$$B = \left\{e_{i,i+1} (1); i = 1,\dots n\right\} \cup \{ e_{n+1,1} (t)\}$$ and let $\mathcal{U} (\mathbb{F}_q) < \SL_{n+1} (\mathbb{F}_q [t])$ be the subgroup generated by $B$, i.e.
$$\mathcal{U} (\mathbb{F}_q) = \langle e_{1,2} (1),  e_{2,3} (1),\dots ,e_{n,n+1} (1),  e_{n+1,1} (t) \rangle.$$
We define subgroups $H_0,\dots, H_n < \mathcal{U} (\mathbb{F}_q)$ as 
\begin{align*}
H_{0} &=  \langle e_{1,2} (1),  e_{2,3} (1),\dots,e_{n,n+1} (1) \rangle = \langle B \setminus \{e_{n+1,1}(t)\} \rangle, \\
H_{1} &=  \langle  e_{2,3} (1),\dots,e_{n,n+1} (1),  e_{n+1,1} (t) \rangle = \langle B \setminus \{e_{1,2}(1)\} \rangle, \\
&\vdots \\
H_{i} &=  \langle  B \setminus \{e_{i,i+1}(1)\} \rangle, \\
& \vdots \\
 H_{n} &= \langle  e_{1,2} (1),  e_{2,3} (1),\dots,e_{n-1,n} (1),  e_{n+1,1} (t) \rangle ,
\end{align*}
and denote $\mathcal{H} = (H_0,\dots , H_n)$.

Let $f \in \mathbb{F}_q [t]$ be an irreducible polynomial of degree $s >1$. Denote $\phi_{f} : \SL_{n+1} (\mathbb{F}_q [t]) \rightarrow \SL_{n+1}(\mathbb{F}_q[t]/(f)) \cong \SL_{n+1} (\mathbb{F}_{q^s})$ defined as follows: For every matrix $A \in  \SL_{n+1} (\mathbb{F}_q [t]),  A= (A_{i,j})_{i,j=1}^n$ where $A_{i,j} \in \mathbb{F}_q [t]$ define 
$$\phi_f (A) =  (A_{i,j} + ( f))_{i,j=1}^n.$$
Restricting $\phi_f$ to $\mathcal{U} (\mathbb{F}_q)$, we get $G_f = \phi_f (\mathcal{U} (\mathbb{F}_q))$ which is a finite quotient of $\mathcal{U} (\mathbb{F}_q)$. It is proven in \cite{hdxfromkms}, that $\phi_f$ is injective on $H_0,\dots,H_n$.  The congurence KMS complex of type $A_n$ over $K$ with respect to $f$ is defined to be the coset complex $X_f = \CC (G_f , \mathcal{H}_f)$,  where $\mathcal{H}_f=(\phi_f(H_0),\dots, \phi_f(H_n))$.   

Let $f_s$ be a sequence of irreducible polynomials of degrees $ \geq 2$ such that $\deg (f_s) \rightarrow \infty$.  Define $X^{(s)} =  \CC (G_{f_s} , \mathcal{H}_{f_s})$.   In \cite{hdxfromkms}  (see also \Cref{subsec: preliminiaries KMS} below) it was proven that for every large enough $q$ with respect to $n$,  the sequence $\lbrace X^{(s)} \rbrace_s$ is a family of uniformly bounded degree pure $n$-dimensional simplicial complexes that are $\frac{1}{\sqrt{q}-(n-1)}$-local spectral expanders. 

All the links of all the vertices in all the $X^{(s)}$'s are isomorphic to the coset complex $\mathcal{CC} (H_0 ; (H_0 \cap H_1, ...,H_0 \cap H_n ))$ (where $H_0,...,H_n$ are as above).  The coset complex $\mathcal{CC} (H_0 ; (H_0 \cap H_1, ...,H_0 \cap H_n ))$ is isomorphic to the opposite complex inside the spherical building of type $A_{n}$ over the field $\mathbb{F}_q$. 

With respect to the family $\lbrace X^{(s)} \rbrace_s$ defined above,  our main result is as follows:
\begin{theorem}
Let $\Lambda$ be a finite group,  $n \geq 3$ be an integer,  $q = p^m$ where $p>2$ is prime.  Also let $\lbrace X^{(s)} \rbrace_s$ be the family of congruence KMS complexes constructed above as coset complexes over quotients of $\mathcal{U} (\mathbb{F}_q) < \SL_{n+1} (\mathbb{F}_q [t])$.   If $q > 2^{n-1}$ and $p > \vert \Lambda \vert$,  then the family $\lbrace X^{(s)} \rbrace_s$,  is a family of $1$-coboundary expanders with $\Lambda$ coefficients.
\end{theorem}

\subsection{Organization}
The paper is organized as follows. In \Cref{sec: preliminaries}, we first recall the relevant notions of high-dimensional expansion in \Cref{subsec: preliminiaries - hdx} and then introduce Chevalley groups, KMS groups and the associated complexes in \Cref{subsec: preliminiaries KMS}. In \Cref{sec: vanishing of first cohomology}, we prove vanishing of the first cohomology of the congruence KMS complexes with respect to many coefficient groups. \Cref{sec: n-a cone functions} is devoted to non-Abelian cone functions, first introducing the concept in \Cref{subsec: def cone and relation CBE}, then showing how to construct cone functions for joins and by adding certain vertices. In \Cref{sec: CSE with n-a coefficients}, we show the cosystolic expansion of classical KMS complexes in the non-Abelian setting. \Cref{app: sec buildings} provides further details on buildings, their opposite complexes and the relation to KMS complexes. \Cref{app: An subcomplexes sec}, \Cref{app: Cn subcomplexes section} and \Cref{app: sec dn case} provide further details on the proof of non-Abelian coboundary expansion of the opposite complexes of type $A_n, C_n$ and $D_n$, respectively. 

\subsection*{Acknowledgements}
The authors are grateful to Pierre-Emmanuel Caprace, Tom De Medts and Hendrik van Maldegem for their advice on terminology. 

The first author is partially supported by ISF grant no. 242/24.
The second author is supported by the FWO and the F.R.S.--FNRS under the Excellence of Science (EOS) program (project ID~40007542).

\section{Preliminaries} \label{sec: preliminaries}
In the first part of this section, we recall the definition and some properties of high-dimensional expanders, following the work of \cite{KO2023high} and \cite{DM} among others. In the second part, we recall the construction of the KMS complexes from \cite{hdxfromkms} and \cite{CSEofKMS}, and recall some of their properties.
\subsection{High-dimensional expanders} \label{subsec: preliminiaries - hdx}

\subsubsection{Weighted simplicial complexes} \label{subsubsec: weighted SCs}

Let $X$ be a finite $n$-dimensional simplicial complex. A simplicial complex $X$ is called {\em pure $n$-dimensional} if every face in $X$ is contained in some face of dimension $n$.  The set of all $k$-faces of $X$ is denoted $X(k)$, and we will be using the convention in which $X(-1) = \{\emptyset\}$.  For every $0 \leq k \leq n$,  the \textit{$k$-skeleton of $X$} is the simplicial complex $\bigcup_{i=-1}^k X(i)$.  We will say that $X$ is \textit{connected} if its $1$-skeleton is a connected graph.

Given a pure $n$-dimensional simplicial complex $X$,  the \textit{weight function} $ w : \bigcup_{k=-1}^n X(k) \rightarrow \mathbb{R}_+$ is defined to be
$$\forall \tau \in X(k), w(\tau) = \frac{\vert \lbrace \sigma \in X(n) : \tau \subseteq \sigma \rbrace \vert}{{n+1 \choose k+1} \vert X(n) \vert}.$$

We note that $w$ is normalized such that for each $-1 \leq k \leq n$, the function $w$ can be thought of as a probability function on $X(k)$.

Given a finite pure $n$-dimensional complex $X$ and $\tau \in X$, \textit{the link of $\tau$} is the subcomplex $\lk_X(\tau)$ defined as 
$$\lk_X(\tau) = \lbrace \sigma \in X : \tau \cap \sigma = \emptyset,  \tau \cup \sigma \in X \rbrace.$$
We note that if $\tau \in X(k)$,  then $\lk_X(\tau)$ is a pure $(n-k-1)$-dimensional complex. We call $\lk_X(\tau)$ a \textit{proper link}, if $1 \leq \dim(\lk_X(\tau))\leq n-1$.

Below,  we will use the following abuse of notation and write $v \in X(0)$ instead of $\lbrace v \rbrace \in X(0)$.  

Given a family of pure $n$-dimensional simplicial complexes $\lbrace X^{(s)} \rbrace_s$,  we say that this family has \textit{bounded degree} if 
$$\sup_s \sup_{v \in X^{(s)} (0)} \vert \lbrace \sigma \in X^{(s)} (n) : v \in \sigma  \rbrace \vert < \infty.$$

\subsubsection{Local spectral expansion}
\label{subsubsec: Local spectral expansion subsec}

Let $X$ be a finite pure $n$-dimensional simplicial complex with $n \geq 1$ and $w$ the weight function on $X$ defined above.  We define the stochastic matrix of the random walk on (the $1$-skeleton of) $X$ to be the matrix indexed by $X(0) \times X(0)$ and defined as 
$$M ( v,  u) = \begin{cases}
 \frac{w( \lbrace u,v \rbrace) }{\sum_{\lbrace u' ,v \rbrace \in X(1)} w( \lbrace u' ,v \rbrace)} & \lbrace v,u \rbrace \in X(1) \\
 0 & \lbrace v,u \rbrace \notin X(1)
\end{cases}.$$

For $\lambda <1$,  we say that $X$ is a \textit{(one-sided) $\lambda$-spectral expander} if $X$ is connected and the second largest eigenvalue of $M$ is $\leq \lambda$.  We say that $X$ is a  \textit{(one-sided) $\lambda$-local spectral expander} if for every $-1 \leq k \leq n-2$ and every $\tau \in X(k)$,  the simplicial complex $\lk_X(\tau)$ is a (one-sided) $\lambda$-spectral expander.  We note that every (one-sided) $\lambda$-local spectral expander is a (one-sided) $\lambda$-spectral expander, since the link of the empty set is $X$ itself.

\subsubsection{Cohomological notations, cosystolic and coboundary expansion}
\label{subsubsec: cosystolic/coboundary expansion subsec}

Here we recall the basic definitions regarding $1$-coboundary and $1$-cosystolic expansion of a simplicial complex $X$ with coefficients in a general group $\Lambda$.  These definitions appeared in several works (e.g., \cite{DM}) and we claim no originality here.

Throughout,  we denote $X$ to be a finite pure $n$-dimensional simplicial complex.   For every $-1 \leq k \leq n$, $X(k)$ denotes the set of $k$-dimensional simplices of $X$ and $X_\ord (k)$ to be the set of ordered $k$-simplices. 

Given a group $\Lambda$,  we define the space of $k$-cochains for $k=-1,0,1$ with values in $\Lambda$ by

\begin{align*}
C^{-1} (X,  \Lambda) &= \lbrace \phi : \lbrace \emptyset \rbrace  \rightarrow \Lambda \rbrace,\\
C^{0} (X,  \Lambda) &= \lbrace \phi : X(0) \rightarrow \Lambda \rbrace,\\
C^1 (X,  \Lambda) &= \lbrace \phi :   X_{\ord} (1) \rightarrow \Lambda : \forall (u,v) \in X_{\ord} (1),  \phi ((u,v)) = ( \phi ((v,u)))^{-1} \rbrace.
\end{align*}
We note that we do not assume that $\Lambda$ is Abelian or finite.

For $k=-1,0$,  we define the coboundary maps $d_k : C^{k} (X,  \Lambda) \rightarrow C^{k+1} (X,  \Lambda)$ as follows: First,  for every $\phi \in C^{-1} (X, \Lambda)$ and every $v \in X(0)$,  we define $ d_{-1} \phi (v) = \phi ( \emptyset)$.  Second,  for every $\phi \in C^{0} (X,  \Lambda)$ and every $(v_0, v_1) \in X_\ord (1)$,  we define $d_{0} \phi ((v_0,v_1)) = \phi (v_0) \phi (v_1)^{-1}$.   Last,  we also define the map $d_1 : C^1 (X,  \Lambda) \rightarrow \lbrace \phi : X_\ord (2) \rightarrow \Lambda \rbrace$ as follows: for every $\phi \in C^{1} (X,  \Lambda)$ and every $(v_0, v_1,v_2) \in X_\ord (2)$,  we define
$d_{1} \phi ((v_0,v_1,v_2)) = \phi ((v_0,v_1))  \phi ((v_1,v_2)) \phi ((v_2,v_0))$.

We note that $d_k d_{k-1} \equiv e_\Lambda$ for $k=0,1$ and define cocycles and coboundaries by
\begin{align*}
Z^k (X, \Lambda) &= \lbrace \phi \in C^{0} (X,  \Lambda) : d_k \phi \equiv e_\Lambda \rbrace,\\
B^k (X, \Lambda)  &= d_{k-1} (C^{k-1} (X, \Lambda)).
\end{align*}
We note that $B^k (X, \Lambda)  \subseteq Z^k (X, \Lambda)$ (as sets). 

We further define the an action of $C^{k-1} (X, \Lambda)$ on $Z^{k} (X, \Lambda)$:
\begin{itemize}
\item An element $\psi \in C^{-1} (X, \Lambda)$ acts on $Z^{0} (X, \Lambda)$ by
$$\psi . \phi (v) = \psi (\emptyset)^{-1} \phi (v),  \forall \phi \in Z^{0} (X, \Lambda),  \forall v \in X (0).$$
\item An element $\psi \in C^{0} (X, \Lambda)$ acts on $Z^{1} (X, \Lambda)$ by
$$\psi . \phi ((v_0,v_1)) = \psi (v_0) \phi ((v_0,v_1)) \psi (v_1)^{-1},  \forall \phi \in Z^{1} (X, \Lambda),  \forall (v_0,v_1) \in X_{\ord} (1).$$
\end{itemize}

We denote by $[ \phi ]$ the orbit of $\phi \in Z^{k} (X, \Lambda)$ under the action of $C^{k-1} (X, \Lambda)$ and define the reduced $k$-th cohomology to be
$$H^k (X, \Lambda) = \lbrace [ \phi ] : \phi \in Z^{k} (X, \Lambda) \rbrace.$$
We will say that the $k$-th cohomology is trivial if $H^k (X, \Lambda) = \lbrace [ \phi \equiv e_\Lambda ]  \rbrace$ and note that the cohomology is trivial if and only if $B^k (X, \Lambda) = Z^k (X, \Lambda)$.

When $\Lambda$ is an Abelian group,  $Z^k (X, \Lambda),  B^k (X, \Lambda)$ are Abelian groups and for $k=0,1$,  the $k$-th reduced cohomology is defined by $H^k (X,  \Lambda) = Z^k (X, \Lambda) / B^k (X, \Lambda)$.  We leave it to the reader to verify that the definition in the Abelian setting coincides (as sets) with the definition in the general setting given above.

In order to define expansion, we need to define a norm.  Let $ w : \bigcup_{k=-1}^n X(k) \rightarrow \mathbb{R}_+$ be the weight function defined above.
For $\phi \in C^0 (X, \Lambda)$,  we define $\supp (\phi) \subseteq X(0)$ as
$$\supp (\phi) = \lbrace v \in X(0) : \phi (v) \neq e_\Lambda \rbrace.$$
Also,  for $\phi \in C^{1} (X, \Lambda)$,  we define
$$\supp (\phi) = \lbrace \lbrace u,v \rbrace \in X(1) : \phi ((u,v)) \neq e_\Lambda \rbrace,$$
and
$$\supp (d_1 \phi) =  \lbrace \lbrace u,v, w \rbrace \in X(2): (d_1\phi) ((u,v,w)) \neq e_\Lambda \rbrace.$$
We note that for every permutation $\pi$ on $\lbrace 0,1,2\rbrace$ it holds that $d_1 \phi ((v_0,v_1,v_2)) = e_\Lambda$ if and only if $d_1 \phi ((v_{\pi (0)},v_{\pi (1)},v_{\pi (2)})) = e_\Lambda$ and thus $\supp (d_1 \phi)$ if well-defined. 

With this notation,  we define the following norms: for $\phi \in  C^{k} (X, \Lambda)$,  $k=0,1$,  we define
$$\Vert \phi \Vert = \sum_{\tau \in \supp (\phi)} w (\tau).$$
Also,  for $\phi \in  C^{1} (X, \Lambda)$, we define
$$\Vert d_1 \phi \Vert = \sum_{\tau \in \supp (d_1 \phi)} w (\tau).$$

Last,  for two maps $\phi, \psi \in  C^{k} (X, \Lambda),  k=0,1$,  we define $\dist (\phi, \psi)$ as follows: For $k=0$,  we define
$$\dist (\phi, \psi) = w \left( \lbrace v \in X(0) : \phi (v) (\psi (v))^{-1} \neq e_\Lambda \rbrace \right).$$
For $k=1$,  we define 
$$\dist (\phi, \psi) = w \left( \lbrace \lbrace u,v \rbrace \in X(1) : \phi ((u,v)) (\psi ((u,v)))^{-1} \neq e_\Lambda \rbrace \right).$$

Using the notations above,  we will define coboundary and cosystolic expansion constants (generalizing the Cheeger constant) as follows: For $k =0,1$,  define
$$h^{k}_{\cobound} (X, \Lambda) = \min_{\phi \in C^{k} (X, \Lambda) \setminus B^{k} (X, \Lambda)} \dfrac{\Vert d_k \phi \Vert}{\min_{\psi \in B^k (X,\Lambda)} \dist (\phi,  \psi) },$$
$$h^{k}_{\cosys} (X, \Lambda) = \min_{\phi \in C^{k} (X, \Lambda) \setminus Z^{k} (X, \Lambda)} \dfrac{\Vert d_k \phi \Vert}{\min_{\psi \in Z^k (X,\Lambda)} \dist (\phi,  \psi)}.$$

Observe the following (for $k =0,1$):
\begin{itemize}
\item For any $X$,  $h^{k}_{\cosys} (X, \Lambda) > 0$.
\item If $H^{k}  (X, \Lambda)$ is non-trivial, then $h^{k}_{\cobound} (X, \Lambda) =0$.
\item If $H^{k}  (X, \Lambda)$ is trivial,  then $h^{k}_{\cobound} (X, \Lambda) = h^{k}_{\cosys} (X, \Lambda)$.
\end{itemize}

For a constant $\beta >0$,  we say that
\begin{itemize}
\item The complex $X$ is a \textit{$(\Lambda,  \beta)$ $1$-coboundary expander},  if $h^{0}_{\cobound} (X, \Lambda) \geq \beta$ and $h^{1}_{\cobound} (X, \Lambda) \geq \beta$.
\item The complex $X$ is a \textit{$(\Lambda,  \beta)$ $1$-cosystolic expander},  if $h^{0}_{\cosys} (X, \Lambda) \geq \beta$,  $h^{1}_{\cosys} (X, \Lambda) \geq \beta$ and for every $\phi \in Z^1 (X,  \Lambda) \setminus B^1 (X,  \Lambda)$ it holds that $\Vert \phi \Vert \geq \beta$.
\end{itemize}

We say that a family of pure $n$-dimensional simplicial complexes $\lbrace X^{(s)} \rbrace_s$,  is a family of  \textit{$1$-coboundary expanders with $\Lambda$ coefficients} ( \textit{$1$-cosystolic expanders with $\Lambda$ coefficients}), if the the following holds:
\begin{itemize}
\item The family $\lbrace X^{(s)} \rbrace_s$ is of bounded degree,  i.e., 
$$\sup_s \sup_{v \in X^{(s)} (0)} \vert \lbrace \sigma \in X(n) : v \in \sigma \rbrace \vert < \infty.$$
\item The size of $X^{(s)}$ grows with $s$, i.e., 
$$\lim_s \vert X^{(s)} (0) \vert = \infty.$$
\item There exists $\beta >0$ such that for every $s$,  $X^{(s)}$ is a $(\Lambda,  \beta)$ $1$-coboundary expander ($(\Lambda,  \beta)$ $1$-cosystolic expander).
\end{itemize}

\begin{remark}
A lower bound on $h^{0}_{\cobound}$ is usually easy to verify: If the $1$-skeleton of $X$ is a connected finite graph, then for any non-trivial group $\Lambda$, $h^{0}_{\cobound} (X, \Lambda)$ is exactly the Cheeger constant of the weighted $1$-skeleton (where each edge $\lbrace u,v \rbrace$ is weighted by $w(\lbrace u,v \rbrace$).  Moreover,  the Cheeger constant of the weighted $1$-skeleton can be bounded by the spectral gap of the weighted random walk on the $1$-skeleton, which in concrete examples can be bounded using the trickling down Theorem from \cite{OppLocI}.
\end{remark}

In \cite{DD-cosys},  the following result was proven (generalizing previous results in \cite{KKL,  EK} that applied to the Abelian setting, and improving parameters over the result of \cite{KM}, that was the first cosystolic expansion result in the non-Abelian setting):
\begin{theorem}\cite[Theorem 8]{DD-cosys}
\label{DD for cosys thm}
Let $0 \leq \lambda <1,  \beta >0$ be constants and $\Lambda$ be a group.  For a finite pure $n$-dimensional simplicial complex $X$ with $n \geq 3$,  if for every vertex $v$,  $h^1_{\cobound} (\lk_X(v), \Lambda) \geq \beta$ and $X$ is a $\lambda$-one-sided local spectral expander,  then
$$h^1_{\cosys} (X, \Lambda) \geq \frac{(1-\lambda) \beta}{24} - e \lambda,$$
where $e \approx 2.71$ is the Euler constant.
\end{theorem}

\subsection{Kac--Moody--Steinberg groups and complexes}
\label{subsec: preliminiaries KMS}

We recall the notions of generalized Cartan matrices and root systems which are the building blocks of Chevalley, Steinberg and Kac-Moody-Steinberg groups. We then describe the construction of KMS complexes and explain how we get KMS complexes over Chevalley groups, called congruence KMS complexes, which will be the objects of our main result, \Cref{cor: main result CBE}.

\subsubsection{Root systems and Chevalley groups} \label{subsubsec: root systems and chevalley groups}
We start by recalling some facts about generalized Cartan matrices, the associated root systems and Chevalley groups. These are the building blocks of KMS groups. 

Most of the structure of the KMS groups is encoded in its generalized Cartan matrix, which also gives rise to a Dynkin diagram and a root system.

\begin{definition} \label{def: GCM and DyDi}
A generalized Cartan matrix (GCM) is a matrix $A = (A_{ij})_{i,j \in I} \in \mat_n(\Z)$ such that $A_{ii}=2$ for all $i \in I$, $A_{ij} \leq 0$ for all $i \neq j \in I$ and $A_{ij}= 0 \iff A_{ji}=0$. We call $\lvert I \rvert$ the rank of $A$.

Every GCM $A=(A_{ij})_{i,j \in I}$ gives rise to a Dynkin diagram in the following way. As vertex set, we take the index set $I$ and two vertices $i,j$ are connected by $|A_{ij}|$ edges if $A_{ij}A_{ji} \leq 4$ and $|A_{ij}| \geq |A_{ji}|$, and these edges are equipped with an arrow pointing towards $i$ if $|A_{ij}|>1$. If $A_{ij}A_{ji}>4$, the vertices $i$ and $j$ are connected by a bold-faced labelled edge with the ordered pair of integers $|A_{ij}|,|A_{ji}|$.

A GCM $A$ is \emph{irreducible} if there exists no non-trivial partition $I = I_1 \cup I_2$ such that $A_{i_1 i_2} = 0$ for all $i_1 \in I_1, i_2 \in I_2$.
\end{definition}

\begin{definition}
Let $A$ be a $(d+1)\times (d+1)$ generalized Cartan matrix with index set $I=\{0,\dots,d\}$, and let $J \subseteq I$.
\begin{enumerate}[label=(\roman*)]

    \item  The subset  $J$ is called \emph{spherical} if $A_J = (A_{ij})_{i,j \in J}$ is of spherical type, meaning that the associated Coxeter group (see e.g. \cite[Proposition 4.22]{marquis2018introduction}) is finite (see e.g. \cite[Chapter 6.4.1]{bourbaki2008lie}).
Given $n  \geq 2$, $A$ is \emph{$n$-spherical} if every subset $J \subseteq I$ of size $n$ is spherical. 

\item We denote by $Q_A$ the set of spherical subsets of $I$ associated to the generalized Cartan matrix $A$.

\item A generalized Cartan matrix $A$ is \emph{purely $n$-spherical} if every spherical subset $J \subseteq I$ is contained in a spherical subset of size $n$.  (In particular, no set of size $n+1$ is spherical for a purely $n$-spherical generalized Cartan matrix.)

\item Analogously, we call a GCM \emph{$n$-classical} if each irreducible factor of $A_J = (A_{ij})_{i,j \in J}$ is of classical type $A_k,B_k,C_k,D_k$ (see e.g. \cite[Chapter 6.4.1]{bourbaki2008lie}) for each $J \subseteq I$ with $\lvert J \rvert \leq n$. 

\item To a generalized Cartan matrix $A$ we can associate the following sets:
\begin{itemize}
    \item a set of simple roots $\Pi = \{\alpha_i \mid i \in I\}$,
    \item a set of real roots $\Phi \subseteq \bigoplus_{i \in I} \Z \alpha_i$,
    \item two sets, one of positive and one of negative real roots: $\Phi^+= \bigoplus_{i \in I}\N \alpha_i \cap \Phi, \Phi^- = - \Phi^+.$
\end{itemize}   
Note that $\Phi = \Phi^+ \sqcup \Phi^-$. More details can be found e.g. in \cite[Chapter 3.5]{marquis2018introduction}.
\end{enumerate}
\end{definition}

\begin{figure}

\begin{align*}
	\text{Classical types } & \begin{cases}
	A_n &
	\begin{tikzpicture}[line width=1pt, scale=.8, baseline=-.6ex]
		\draw (0,0) -- (2,0);
		\draw[dotted] (2,0) -- (4,0);
		\draw (4,0) -- (6,0);
		\diagnode{(0,0)}
		\diagnode{(1,0)}
		\diagnode{(2,0)}
		\diagnode{(4,0)}
		\diagnode{(5,0)}
		\diagnode{(6,0)}
	\end{tikzpicture} \\[1ex]
	B_n &
	\begin{tikzpicture}[line width=1pt, scale=.8, baseline=-.6ex]
		\draw (0,0) -- (2,0);
		\draw[dotted] (2,0) -- (4,0);
		\draw (4,0) -- (5,0);
		\draw[double distance=2.3pt] (5,0) -- (6,0);
		\draw[line width=.7pt] (5.3,.3) -- (5.6,0) -- (5.3,-.3);
		\diagnode{(0,0)}
		\diagnode{(1,0)}
		\diagnode{(2,0)}
		\diagnode{(4,0)}
		\diagnode{(5,0)}
		\diagnode{(6,0)}
	\end{tikzpicture} \\[1ex]
	C_n &
	\begin{tikzpicture}[line width=1pt, scale=.8, baseline=-.6ex]
		\draw (0,0) -- (2,0);
		\draw[dotted] (2,0) -- (4,0);
		\draw (4,0) -- (5,0);
		\draw[double distance=2.3pt] (5,0) -- (6,0);
		\draw[line width=.7pt] (5.7,.3) -- (5.4,0) -- (5.7,-.3);
		\diagnode{(0,0)}
		\diagnode{(1,0)}
		\diagnode{(2,0)}
		\diagnode{(4,0)}
		\diagnode{(5,0)}
		\diagnode{(6,0)}
	\end{tikzpicture} \\[1ex]
	D_n &
	\begin{tikzpicture}[line width=1pt, scale=.8, baseline=-.6ex]
		\draw (0,0) -- (2,0);
		\draw[dotted] (2,0) -- (4,0);
		\draw (4,0) -- (5,0);
		\draw (5,0) -- (6,.4);
		\draw (5,0) -- (6,-.4);
		\diagnode{(0,0)}
		\diagnode{(1,0)}
		\diagnode{(2,0)}
		\diagnode{(4,0)}
		\diagnode{(5,0)}
		\diagnode{(6,.4)}
		\diagnode{(6,-.4)}
	\end{tikzpicture}
	\end{cases} \\[2ex]
	\text{Exceptional types } & \begin{cases}
	E_6, E_7, E_8 &
	\begin{tikzpicture}[line width=1pt, scale=.8, baseline=-.6ex]
		\draw (0,0) -- (4,0);
		\draw[dotted] (4,0) -- (6,0);
		\draw (2,0) -- (2,-1);
		\diagnode{(0,0)}
		\diagnode{(1,0)}
		\diagnode{(2,0)} \diagnode{(2,-1)}
		\diagnode{(3,0)}
		\diagnode{(4,0)}
		\diagnode{(5,0)}
		\diagnode{(6,0)}
	\end{tikzpicture} \\[5ex]
	F_4 & 	\begin{tikzpicture}[line width=1pt, scale=.8, baseline=-.6ex]
		\draw (0,0) -- (1,0);
		\draw[double distance=2.3pt] (1,0) -- (2,0);
		\draw (2,0) -- (3,0);
		\draw[line width=.7pt] (1.3,.3) -- (1.6,0) -- (1.3,-.3);
		\diagnode{(0,0)}
		\diagnode{(1,0)}
		\diagnode{(2,0)}
		\diagnode{(3,0)} \end{tikzpicture} \\[1ex]
	G_2 &	\begin{tikzpicture}[line width=1pt, scale=.8, baseline=-.6ex]
		\draw[double distance=3pt] (0,0) -- (1,0);
		\draw (0,0) -- (1,0);
		\draw[line width=.7pt] (0.3,.3) -- (0.6,0) -- (0.3,-.3);
		\fill (0,0) circle (.11);
		\fill (1,0) circle (.11); \end{tikzpicture} \\[1ex] \end{cases}
\end{align*}	
\caption{The Dynkin diagrams of irreducible spherical root systems\label{table:dynkin}}	
\end{figure}
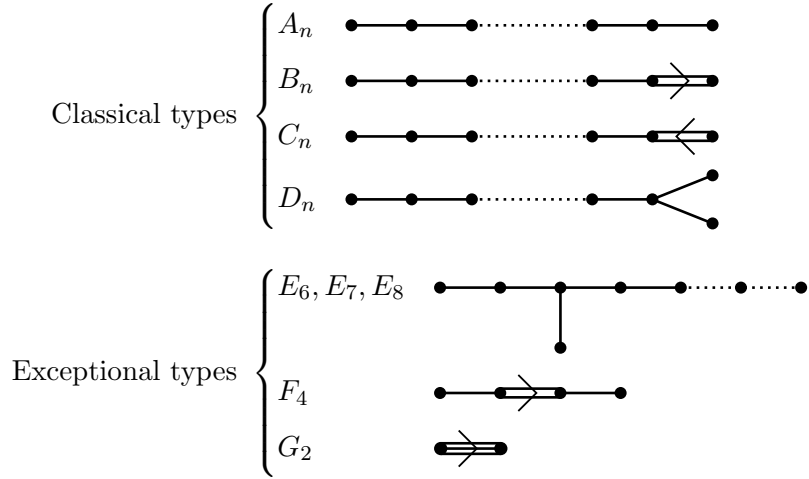

\begin{remark}
\begin{enumerate}[label=(\roman*)]
\item Irreducible spherical  root systems have been classified. All possible diagrams are given in \Cref{table:dynkin}. 
\item All rank $n+1$, $n$-spherical diagrams have been classified and are precisely the affine diagrams and the compact hyperbolic diagrams, the latter have rank at most 5, where the diagram of rank 5 is exactly the one that we have to exclude to get $n$-classicality. See e.g. \cite[Chapter 6.9]{HumphreyReflGrp} for the classification result. 

The non-spherical, rank $n+1$, $n$-classical diagrams, with $n \geq 3$ are precisely the non-spherical rank $n+1$ diagrams which are $n$-spherical but not one of the following
\begin{itemize}
\item affine diagrams $\tilde{E}_6,\tilde{E}_7,\tilde{E}_8,\tilde{F}_4,\tilde{G}_2$,
\item compact hyperbolic diagram \begin{tikzpicture}[line width=1pt, scale=.8, baseline=-.6ex]
		\draw (0,0) -- (1,0);
		\draw[double distance=2.3pt] (1,0) -- (2,0);
		\draw (2,0) -- (3,0);
		\draw[line width=.7pt] (1.3,.3) -- (1.6,0) -- (1.3,-.3);
		\draw (0,0) -- (1.5,1);
		\draw (3,0) -- (1.5,1);
		\diagnode{(0,0)}
		\diagnode{(1,0)}
		\diagnode{(2,0)}
		\diagnode{(3,0)}
		\diagnode{(1.5,1)} 
		\end{tikzpicture}.
\end{itemize}

\end{enumerate}

\end{remark}

\begin{definition} \label{def: chevalley group}
Corresponding to any irreducible, spherical GCM $\Ao$ with root system $\mathring \Phi$ of rank at least 2, and any finite field $K$, there is an associated universal (or simply connected) Chevalley group, denoted $\chev_{\Ao}(K)$. Abstractly, it is generated by symbols $x_\alpha(s)$ for $\alpha \in \mathring \Phi$ and $s \in K$, subject to the relations
$$
\begin{aligned}
x_\alpha(s) x_\alpha(u) & =x_\alpha(s+u) \\
{\left[x_\alpha(s), x_\beta(u)\right] } & =\prod_{i, j>0} x_{i \alpha+j \beta}\left(C_{i j}^{\alpha, \beta} s^i u^j\right) \quad(\text {for } \alpha+\beta \neq 0) \\
h_\alpha(s) h_\alpha(u) & =h_\alpha(s u) \quad(\text { for } s,u \neq 0), \\
\text { where } \quad h_\alpha(s) & =n_\alpha(s) n_\alpha(-1) \\
\text { and } \quad n_\alpha(s) & =x_\alpha(s) x_{-\alpha}\left(-s^{-1}\right) x_\alpha(s) .
\end{aligned}
$$
Note that the $C_{ij}^{\alpha, \beta}$ are integers called structure constants that can be found in \cite{carter1989simple}.
\end{definition}

\begin{remark}
\begin{enumerate}
	\item Let $K$ be a field and let $K[t]$ denote the polynomial ring in one variable over $K$. The simply-connected Chevalley group $\chev_{\Ao}(K[t])$ is, similar to the case of a Chevalley group over a field, generated by elements $x_\alpha(s)$ for $\alpha \in \mathring \Phi, s \in K[t]$ that satisfy the relations above, where for the third relation we have to add the extra assumption that $u,s$ are invertible in $K[t]$. This can be found e.g. in \cite{rehmann1975prasentationen}.
	\item More generally, one can define Chevalley groups by introducing the \emph{Chevalley--Demazure group scheme}, which is a functor from the category of commutative unital rings to the category of groups. Over fields and their polynomial ring in one variable these definitions yield the same groups.
\end{enumerate}	
\end{remark}
Chevalley groups are closely related to Steinberg groups, which are defined as follows.
\begin{definition}
	Given an irreducible, spherical Cartan matrix $\Ao$ with root system $\mathring \Phi$ of rank at least 2, and a ring $R$, we define the Steinberg group of type $\Ao$ over $R$ to be the group given by the following presentation. 
	\begin{align*}
	\St(R):=\St(\Ao,R) := \langle x_\alpha(u); \alpha \in \mathring \Phi, u \in R \mid \mathcal{R} \rangle
	\end{align*}
	where the relations $\mathcal{R}$ are
	\begin{align*}
	x_\alpha(s) x_\alpha(u) & =x_\alpha(s+u) \quad \text{ for } \alpha \in \mathring \Phi, u,s \in R, \\
	{\left[x_\alpha(s), x_\beta(u)\right] } & =\prod_{i, j>0} x_{i \alpha+j \beta}\left(C_{i j}^{\alpha, \beta} s^i u^j\right) \quad \text {for } \alpha, \beta \in \mathring \Phi, \alpha+\beta \neq 0, s,u \in R^\times .
	\end{align*}
	\end{definition}
\begin{remark} \label{rmk:St vs Chev}
	Note that all the relations of the Steinberg group also hold in the corresponding Chevalley group, while the Chevalley group has potentially more relations.
	
	Over finite fields the Chevalley group and Steinberg group coincide, see e.g. \cite[Corollary 6.1.1]{WeibelKBook}.
\end{remark}

\subsubsection{Kac--Moody--Steinberg groups}
We can now put the ingredients together to define what a Kac--Moody--Steinberg group is and present some of its properties. For more details, see \cite[Chapter 3]{hdxfromkms} and references therein. 
\begin{definition}
Let $K$ be a field. Let $A=(A_{ij})_{i,j\in I}$ be a GCM over the index set $I$ which is 2-spherical. Let $Q_A=\{J \subseteq I \mid A_J:=(A_{ij})_{i,j \in J} \text{ is spherical}\}$. Let $\Phi$ be the (real) root system associated to $A$ with simple roots $\Pi = \{\alpha_i \mid i \in I\}$. For $J\in Q_A$ set 
$$U_J:= \langle x_{\alpha_i}(s) \mid s \in K, i \in J \rangle \leq \chev_{A_J}(K) $$
the group generated by all simple roots in the Chevalley group of type $A_J$ over $K$. Note that if $L \subset J$ then we have a natural inclusion $U_L \hookrightarrow U_J$ by sending $x_{\alpha_i}(s)\in U_L$ to the same generator in $U_J$.  
The KMS group of type $A$ over $K$ is defined as the free product of the $U_J, J\in Q_A$ modulo the natural inclusions:
$$\uak = \freeprod_{J \in Q_A} U_J/(U_L \hookrightarrow U_J, L\subseteq J).$$ 
\end{definition}
Note that $U_J \hookrightarrow \uak$ and we will denote the image of $U_J$ in $\uak$ again by $U_J$. These subgroups are called the local groups of $\uak$. For $i,j\in I$ we will write $U_i:= U_{\{i\}}$ and $U_{i,j}=U_{\{i,j\}}$.

\begin{remark} \label{rmk: presentation of KMS}
For two roots $\alpha, \beta \in \Phi$ we write $]\alpha, \beta[_\N = \left\{n_1 \alpha + n_2 \beta \in \Phi \mid n_1,n_2 \in \N^* \right\}$ and $[\alpha,\beta]_\N = \left]\alpha, \beta\right[_\N \cup \left\{\alpha,\beta\right\}$.

We have the following abstract presentation for the KMS groups:
$$\mathcal{U}_A(K)= \Bigl\langle u_\beta(t) \text{ for } t \in K,  i,j \in I , \beta \in [\alpha_i,\alpha_j]_{\mathbb{N}} \ | \ \mathcal{R} \Bigr\rangle$$
where the set of relations $\mathcal{R}$ is defined as:

\indent for all $i,j \in I , \  \{\alpha,\beta \} \subseteq [\alpha_i, \alpha_j]_{\mathbb{N}},   \ t,u\in K$: 
\begin{align*}
u_\alpha(t) u_\alpha(s) &= u_\alpha(t+s) \\
[u_\alpha(t),u_\beta(u)] &= \prod_{\gamma = k\alpha + l\beta \in ]\alpha,\beta[_{\mathbb{N}}} u_\gamma \left(C_{k,l}^{\alpha \beta} t^k u^l\right).
\end{align*}
The constants $C_{k,l}^{\alpha,\beta}$ are the same structure constants as in the presentation of the Chevalley group.

Since $A$ is 2-spherical, the subgroups $U_\beta$, for $ \beta \in \left[ \alpha_i, \alpha_j\right] $, are contained in the group generated by the root groups $U_{\alpha_i}, U_{\alpha_j}$ (see \cite[Proposition 7]{Abr}).
\end{remark}

\subsubsection{KMS complexes}
We now use the KMS groups defined above in the coset complex construction to get pure simplicial complexes called KMS complexes, which were first introduced in \cite{hdxfromkms}. First,  we recall the definition of coset complexes and state one key fact. More details on coset complexes can be found for example in \cite{KO2023high}.

\begin{definition} \label{def: coset complex}
Let $G$ be a group and $\mathcal{H}=(H_0,\dots , H_n)$ be a family of subgroups of $G$. Then the coset complex $\CC(G;(H_0, \dots, H_n))$ is defined to be the simplicial complex with 
\begin{itemize}
\item vertex set $\bigsqcup_{i=0}^n G/H_i$
\item a set of vertices $\left\{ g_1H_{i_1},\dots,g_kH_{i_k} \right\}$ forms a $(k-1)$-simplex in $\CC(G; \mathcal{H})$ if and only if $\bigcap_{j=1}^k g_j H_{i_j} \neq \emptyset$.
\end{itemize}
\end{definition}

Note that the set of vertices is partitioned into $n+1$ subsets $G/H_i$ for $i = 0,\dots,n$ such that for any simplex $\sigma \in \CC(G;\mathcal{H})$ we have $|\sigma \cap G/H_i| \leq 1$ for all $i$. In that case, we say that $\CC(G;\mathcal{H})$ is a $(n+1)$-partite simplicial complex. 

\begin{definition} \label{def: type of a face}
    The \emph{type} of a simplex $\left\{ g_1H_{i_1},\dots,g_kH_{i_k} \right\}$ in $\CC(G;\mathcal{H})$ is the set of indices $\{i_1,\dots,i_k\} \subseteq \{0,\dots,n\}$. Given a type $\emptyset \neq T \subseteq \{0,\dots,n\}$ we write 
    $$H_T:= \bigcap_{i \in T} H_i \qquad \text{ and set } \qquad H_{\emptyset}:= \langle H_0,\dots,H_n \rangle \leq G.$$
\end{definition}

The following well-known fact characterizes links of faces in a coset complex. 

\begin{fact} \label{fact: links in CC} 
Let $\sigma$ be a face in $\CC (G , \mathcal{H})$ of type $T \neq \emptyset$. Then the link of $\sigma$ is isomorphic to the coset complex $\CC \left(H_T ,\left(H_{T \cup\{i\}}: i \notin T\right)\right)$.
\end{fact}

KMS complexes are coset complexes over certain finite quotients of KMS groups. The precise definition is as follows:

\begin{definition} \label{def: KMS complex}
Let $A$ be a GCM over the index set $I$ which is $(\vert I \vert-1)$-spherical but non-spherical. Let $K$ be a finite field of size $q \geq 4$. Let $\phi: \uak \to G$ be a finite quotient of $\uak$ such that
\begin{enumerate}
\item $\phi|_{U_J}$ is injective for all $J \in Q_A$ (i.e. $\phi$ is injective on the local groups), 
\item $\phi(U_J \cap U_L) = \phi(U_J) \cap \phi(U_L)$ for all $J,L \in Q_A$.
\end{enumerate}
Then we set 
$$X= \CC(G; (\phi(U_{I \setminus \{j \}}))_{j \in I})$$ 
and call $X$ a KMS complex of type $A$ over $K$. 
\end{definition}

It is know that KMS complexes satisfying some mild assumptions are local-spectral expander:

\begin{theorem}{\cite[Theorem 4.3]{hdxfromkms}}
\label{thm: KMS spectral thm}
Let $A$ be a GCM over the index set $I$ which is non-spherical but $(\lvert I \rvert -1)$-spherical. Let $K$ be a finite field of size $q \geq 4$. Let $\phi_i: \uak \to G_i, i\in I$ be family of finite quotients of $\uak$ such that $X_i = \CC(G_i, (\phi_i(U_{I \setminus \{j \}}))_{j \in I}), i \in \N$ is a family of KMS complexes and $\lvert G_i \rvert \overset{i \to \infty}{\longrightarrow} \infty$.
If $q$ is such that $\sqrt{\frac{3}{q}} \leq \frac{1}{\lvert I \rvert -1}$ then the family $(X_i)_{i\in \N}$ is a family of bounded degree $\lambda$-local spectral expanders for some $\lambda = \lambda (q) <1$.  

Moreover,  if $q \gg (\vert I \vert -1)^2$,  then the family $(X_i)_{i\in \N}$ is a family of bounded degree $\frac{2}{\sqrt{q}}$-local spectral expanders. 
\end{theorem}

It was proven in \cite{CSEofKMS} that the $(n-1)$-skeleton of an $n$-classical KMS complex (over a large enough field) is a cosystolic expander with Abelian coefficients.

\subsubsection{Congruence KMS complexes} \label{subsubsec: Congruence KMS}
In this part, we explore the relation between affine KMS groups and the corresponding Chevalley groups, following the results of \cite[Chapter 6]{hdxfromkms}.  We will repeat an explicit construction of a family of simplicial complexes given in \cite[Chapter 6]{hdxfromkms} which we will call here \textit{congruence KMS complexes} (this terminology did not appear in \cite{hdxfromkms}).

Let $\Ao = (A_{ij})_{i,j \in \mathring{I}}, \mathring{I} = \{1,\dots n \}$ be an irreducible spherical GCM with spherical root system $\mathring \Phi$ and set of simple roots $\mathring \Pi = \{\alpha_1,\dots,\alpha_n\}$. This root system has a unique highest root $\gamma \in \mathring \Phi$, which means that $\gamma$ is such that for any $\alpha_i \in \mathring \Pi$ we have $\gamma + \alpha_i \notin \mathring \Phi$. 

Note that, whenever we remove one root from the set $S:= \{-\gamma, \alpha_1,\dots,\alpha_n\}$, we get a set of simple roots (although they might not generate the full root system $\mathring \Phi$ but a subsystem).

We set $\alpha_0 = -\gamma$ and $I = \{0,\dots, n\}$. The idea now is to consider the set $\Pi = \{\abar_0,\dots,\abar_n\}$ as set of simple roots of a root system $\Phi$ such that for all $i,j \in I$ the roots $\abar_i,\abar_j$ generate the same rank two subsystem as $\alpha_i, \alpha_j$ in $\mathring{\Phi}$. The resulting root systems $\Phi$ with simple roots $\Pi = \{\abar_0,\dots,\abar_n \}$ and generalized Cartan matrix $A$ are described by the Dynkin diagrams in \cite[\S 4.8]{KacInfDimLieAlg} (or \cite[Table 5.1 Aff 1]{marquis2018introduction}) and are of affine (untwisted) type. Note that $A$ is of rank $n+1$, $n$-spherical and purely $n$-spherical. If $\Ao$ is of classical type, then $A$ is $n$-classical.

For a root $\beta = \sum_{i \in \mathring{I}} \lambda_i \alpha_i \in \mathring{\Phi}$, $\lambda_i \in \Z$, we write $\overline{\beta} = \sum_{i \in \mathring{I}} \lambda_i \abar_i \in \Phi$. We set $\overline{\delta} = \abar_0 + \overline{\gamma} \in \bigoplus_{i \in I} \Z \abar_i$ which is not in $\Phi$ but it is what is called an imaginary root of the root system associated to $A$. For our purposes, it suffices to think of it as an element in $\bigoplus_{i \in I} \Z \abar_i$. From \cite[Section 7.4]{KacInfDimLieAlg} it follows that $\Phi = \{ a_{\alpha,m}:= \abar + m \overline{\delta} \mid m \in \Z, \alpha \in \mathring \Phi \}$. Note that $\abar_0 = a_{-\gamma,1}$. A root $a_{\alpha,m}$ is in $\Phi^+$ if and only if $m\geq 1$ or $m=0$ and $\alpha \in \mathring \Phi^+$ (see e.g. \cite[proof of Theorem 7.90]{marquis2018introduction}).

Before we continue with the construction, we fix a few notations: $K[t]$ denotes the polynomial ring in one variable $t$ over a field $K$. For a polynomial $f \in K[t]$, let $(f) = \{g \cdot f \mid g \in K[t]\}$ denote the ideal generated by $f$ in $K[t]$ (e.g. $(t)$ is the set of polynomials without constant term) and let $Kt^n = \{\lambda t^n \mid \lambda \in K\}$ be the set of all scalar multiples of $t^n$.

Let $K$ be a finite field of order $|K| \geq 5$ with $\cha(K) \neq 2$ and assume that $\operatorname{rank}(\mathring \Phi) \geq 3$ (thus we have that $\uak$ is isomorphic to the maximal positive unipotent subgroup of the corresponding Kac--Moody group, see \cite[Proposition 3.8]{hdxfromkms}). The aforementioned result, together with \cite[Theorem 7.90]{marquis2018introduction} imply that the following is a well-defined injective homomorphism
\begin{alignat*}{1}
    \phi: \uak & \to \chev_{\Ao}(K[t])  \\
    u_{a_{\alpha,m}}(\lambda) & \mapsto x_\alpha(\lambda t^m).
\end{alignat*}

Let $f \in K[t]$ be an irreducible polynomial of degree $\ell \geq 2$. We denote by $ \pi_f: K[t] \to K[t]/(f)$ the projection to the quotient of the polynomial ring by the ideal generated by $f$. 

Using the functoriality of the Chevalley group functor, we can go from the Chevalley group over $K[t]$ to the one over $K[t]/(f)$ by applying $\pi_f$ ``entry-wise''. We get the following well-defined map
\begin{alignat*}{2}
    \varphi_f: \uak & \overset{\phi}{\longrightarrow} \chev_{\Ao}(K[t]) && \overset{\pi_f}{\longrightarrow} \chev_{\Ao}(K[t]/(f))\\
    u_{\alpha_i}(\lambda) & \mapsto \begin{cases} x_{\alpha_i}(\lambda)  & i \neq 0 \\ x_{-\gamma}(\lambda t) & i = 0 \end{cases} && \mapsto \begin{cases} x_{\alpha_i}(\lambda + (f))  & i \neq 0 \\ x_{-\gamma}(\lambda t + (f)) & i = 0 \end{cases}.
\end{alignat*}

In \cite[Chapter 6]{hdxfromkms} it was shown that the maps $\varphi_f$ satisfy the assumptions of \Cref{def: KMS complex}, hence give rise to KMS complexes.  In particular, choosing an infinite family of irreducible polynomials with increasing degree gives rise to a family of bounded degree high-dimensional spectral expanders, see \cite[Theorem 5.13]{hdxfromkms}.  

This leads to the following definition of congruence KMS complexes:  
\begin{definition}[Congruence KMS complexes]
\label{def:congKMS}
For every spherical root system $\mathring \Phi$ of type $\Ao$ and rank $n$, every finite field $K$ of order $\geq 5$ and characteristic different from 2, and every irreducible polynomial   $f \in K[t]$ of degree $\geq 2$,  we define the congruence KMS complex $X (\mathring \Phi,  K,  f)$ as 
$$X (\mathring \Phi,  K,  f) = \mathcal{CC}(\chev_{\Ao}(K[t]/(f));(\pi_f(H_i))_{i \in \{0,\dots , n\}})$$
where $\pi_f: \chev_{\Ao}(K[t]) \to \chev_{\Ao}(K[t]/(f))$ is the projection induced by the map $K[t] \to K[t]/(f)$, and for $i=1, \dots, n$ we set $H_i = \langle x_{-\gamma}(\lambda t), x_{\alpha_j}(\lambda) \mid j \in \{1,\dots,n\} \setminus \{i\}, \lambda \in K  \rangle \leq \chev_{\Ao}(K[t])$ and $H_0 = \langle  x_{\alpha_j}(\lambda) \mid j \in \{1,\dots,n\}, \lambda \in K  \rangle \leq \chev_{\Ao}(K[t])$. 

Given a sequence of irreducible polynomials  $f_s \in K[t]$ of degree $\ell \geq 2$ such that $\deg (f_s)$ tends to infinity,  we define the family of congruence KMS complexes associated to the sequence $\lbrace f_s \rbrace$ to be the family $\lbrace  X^{(s)} = X (\mathring \Phi,  K,  f_s) \rbrace_s$.

Last,  we define congruence KMS complexes of type $X$ over a field $K$ to be any family of congruence KMS complexes $\lbrace  X^{(s)} = X (\mathring \Phi,  K,  f_s) \rbrace_s$ associated to some sequence $\lbrace f_s \rbrace$ as above where $\mathring \Phi$ is of type $X$ (e.g.,  $\mathring \Phi$ is of type $X =A_n$).   
\end{definition}

\section{Vanishing of first cohomology}
\label{sec: vanishing of first cohomology}
In this chapter, we show that for the congruence KMS complexes the first cohomology with many coefficient groups vanishes. Together with the results from \Cref{sec: CSE with n-a coefficients} this shows that the congruence KMS complexes give rise to families of bounded degree 1-coboundary expanders.  

We want to apply \cite[Theorem 3.1]{KO24cbeofcoco} to our setting. We recall the statement: 
\begin{theorem}[\cite{KO24cbeofcoco}] \label{thm: homology of quotients}
Let $X$ be a connected pure 2-dimensional simplicial complex and $N$ be a group acting on $X$. Assume $p: X \to N\backslash X$ is rigid, i.e. it is simplicial and for all $\sigma \in X$ we have $\lvert \sigma \rvert = \lvert p(\sigma)\rvert$. For every group $\Lambda$, if $H^1(X,\Lambda) = 0$ and $H^1(N,\Lambda) = 0$ then $H^1(N \backslash X, \Lambda)=0$ (the cohomology of $N$ is taken with the trivial action of $N$ on $\Lambda$). 
\end{theorem}

In the remainder of the section, we will first investigate the underlying infinite complex in our setting, and then the kernel of the quotient map.

\subsection{The underlying infinite complex} \label{subsec: underlying infinite complex}
The complex that will play the role of $X$ in \Cref{thm: homology of quotients} is the following.
\begin{definition}
Given a $n$-spherical, purely $n$-spherical GCM $A=(A_{ij})_{i,j \in I}$ of rank $n+1$ and a field $K$, we define the \emph{fundamental KMS complex} to be 
$$Y_A(K) = \CC(\uak, (U_{I\setminus \{i\}})_{i \in I}).$$
\end{definition}

\begin{observation}
The complex $Y_A(K)$ is a connected, pure, $n$-dimensional simplicial complex. 
\end{observation}

\begin{proposition} \label{prop: simple connectedness}
Let $A$ be a $n$-spherical, purely $n$-spherical GCM of rank $n+1$, then $Y_A(K)$ is simply connected, in particular $H^1(Y_A(K),\Lambda)=0$ for every group $\Lambda$.
\end{proposition}
\begin{proof}
This follows from \cite[Theorem 2.4]{AbelsHolz} since $\uak$ can be described as the free product of the $U_{I \setminus \{i\}}$ amalgamated along their intersections.
\end{proof}

\begin{lemma} \label{lemma: quotients have correct form}
Given $\varphi_f: \uak \to \chev_{\Ao}(K[t]/(f))$ as in \Cref{subsubsec: Congruence KMS}. Let $N_f = \ker \varphi_f$. Then $N_f \backslash Y_A(K) \cong \CC(\chev_{\Ao}(K[t]/(f)), (\varphi_f(U_{I\setminus \{i\}})_{i \in I})$ and $p_f: Y_A(k) \to N_f \backslash Y_A(K)$ is rigid.  
\end{lemma}

\begin{proof}
This follows from \cite[Proposition 2.5]{KO24cbeofcoco} and \cite[Observation 2.2]{KO24cbeofcoco}.
\end{proof}

\subsection{The kernel of the quotient map}
\label{subsec: The kernel of the quotient map}

Let $A$ by a generalized Cartan matrix of affine type and rank $n\geq 4$, and let $K$ be a finite field with $\lvert K \rvert \geq 5$. For simplicity, we assume that $\operatorname{char}(K) \neq 2$, see \cite[Remark 5.10]{hdxfromkms}. 
Let $\Ao$ denote the spherical GCM related to $A$ as in the construction in \Cref{subsubsec: Congruence KMS}, with root system $\mathring \Phi$. Recall from \Cref{subsubsec: Congruence KMS} that we have an injective homomorphism

\begin{align*}
    \phi: \uak & \to \chev_{\Ao}(K[t]) \\
    u_{\alpha_i}(\lambda) & \mapsto \begin{cases} x_{\alpha_i}(\lambda ) & i \neq 0 \\
    x_{-\gamma}(\lambda t) & i = 0 \end{cases}
\end{align*}
where $\gamma$ is the longest root in $\mathring \Phi$. 

Let $\mathcal{I}$ denote the image of $\phi$ in $\chev(K[t])$ (since the type of the Chevalley group is fixed here, we will omit if from the notation). As shown in \cite[Proposition 5.11]{hdxfromkms}, we have
$$\mathcal{I}= \langle x_{\alpha}(g) \mid \alpha \in \mathring \Phi^+ \text{ and } g \in K[t], \text{ or } \alpha \in \mathring \Phi^- \text{ and } g\in (t) \rangle \leq \operatorname{Chev}(K[t]).$$
Furthermore, let $f\in K[t]$ be irreduible with $\deg(f) \geq 2$.

The quotient map $K[t] \to K[t]/(f)$ induces the map $\pi_f: \chev(K[t]) \to \chev(K[t]/(f))$. We will denote its kernel by $\chev(K[t], (f)):=\ker(\pi_f)$. 

We introduce some further notation, where for a group $G$ and a subset $S$ we denote by $\langle \langle S \rangle \rangle_G$ the normal subgroup of $G$ generated by $S$. 
\begin{align*}
\mathcal{I}(f) &= \langle \langle x_{\alpha}(g) \mid \alpha \in \mathring \Phi^+ \text{ and } g \in (f), \text{ or }  \alpha \in \mathring \Phi^- \text{ and } g\in (tf) \rangle \rangle_{\mathcal{I}}, \\
\operatorname{St} \left(\mathcal{I} \right) &= \langle x_{\alpha}(g) \mid \alpha \in \mathring \Phi^+  \text{ and } g \in K[t], \text{ or }  \alpha \in \mathring \Phi^- \text{ and } g\in (t) \rangle \leq \operatorname{St} \left( K[t] \right), \\
\operatorname{St} \left( \mathcal{I},(f) \right) &= \langle \langle x_{\alpha}(g) \mid \alpha \in \mathring \Phi^+ \text{ and } g \in (f), \text{ or } \alpha \in \mathring \Phi^-  \text{ and } g\in (tf) \rangle \rangle_{\operatorname{St} \left( \mathcal{I} \right)}.
\end{align*}

The following result is of a similar flavor as \cite[Lemma 3.3]{ApteCR} and the proof is inspired by the proof of the aforementioned result.
\begin{theorem}\label{thm: kernel gen by order p elements}
In the above set up, we have
$$\chev(K[t],(f))\cap \mathcal{I} = \mathcal{I}(f).$$
\end{theorem}
\begin{proof}
We have the following commutative diagram.
\begin{center}
\begin{tikzcd}
0 \arrow[r, "\alpha"] \arrow[d]
& 0 \arrow[d ] \\
\St(\mathcal{I},(f)) \arrow[r, "\beta"] \arrow[d]
& \chev(K[t],(f)) \cap \mathcal{I} \arrow[d ] \\
\St(\mathcal{I}) \arrow[r, "\gamma"] \arrow[d]
& \mathcal{I} \arrow[d ] \\
\St(K[t]/(f)) \arrow[r, "\delta" ]
&  \chev(K[t]/(f))
\end{tikzcd}
\end{center}

Here, $\alpha$ is the identity, $\beta, \gamma$ and $\delta$ are given by mapping a generator $x_\alpha(g)$ of the Steinberg group to the corresponding generator (with the same label) in the Chevalley group.

Note that $\mathcal{I}(f) = \beta(\St(\mathcal{I},(f))) \subseteq \chev(K[t],(f))\cap \mathcal{I}$. 

The vertical maps in the diagram are first the trivial map, then the inclusion map followed by the quotient maps $\St(K[t]) \to \St(K[t]/(f))$ and $\chev(K[t]) \to \chev(K[t]/(f))$ respectively, both restricted to the desired domain. 

Our goal is to show that $\beta$ is surjective. To this end, we want to apply the ``Four lemma'' (see e.g. \cite[Lemma I.3.2]{MacLane63}) which says that if the columns are exact, $\delta$ is injective and $\alpha$ and $\gamma$ are surjective then $\beta$ is surjective.  

Note that $\alpha$ is trivially surjective, $\delta$ is an isomorphism, since $K[t]/(f)$ is a finite field (see \Cref{rmk:St vs Chev}) and $\gamma$ is surjective since all generators of $\mathcal{I}$ lie in the image of $\St(\mathcal{I})$.

It remains to show the exactness of the columns. For the right column, this follows from the definitions.

For the left column we need to show that $M:= \ker(\operatorname{St} \left( \mathcal{I} \right) \to \operatorname{St} \left( K[t]/(f) \right))$ is equal to  $N:=\operatorname{St} \left( \mathcal{I},(f) \right)$.

Clearly we have that $N \subseteq M$. For the other inclusion, we proceed in two steps. 

Step 1:
We show that $\operatorname{St} \left( \mathcal{I} \right)/N \cong \operatorname{St} \left( \mathcal{I} \right)/M$. 

Thus, let $\psi: \operatorname{St} \left( \mathcal{I} \right) \to \operatorname{St} \left( \mathcal{I} \right)/N$ denote the projection. 
For $\alpha \in \mathring \Phi^+, g \in K[t], h \in (f)$ we have
$$
\psi(x_{\alpha}(g+h)) = \psi(x_{a}(g)) \psi(x_{\alpha}(h)) = \psi(x_{\alpha}(g))
$$
since $x_{\alpha}(h) \in N$.
Similarly, for $\alpha \in \mathring \Phi^-, g \in (t), h \in (tf)$ we have
$$
\psi(x_{\alpha}(g+h)) = \psi(x_{a}(g)) \psi(x_{\alpha}(h)) = \psi(x_{\alpha}(g)).
$$
Thus the expressions $\psi(x_{\alpha}(g+(f)))$ for $\alpha\in \mathring \Phi^+, g \in K[t]$ and $\psi(x_{\alpha}(g+(tf)))$ for $\alpha \in \mathring \Phi^-, g \in (t)$ are well defined. In particular, they form a set of generators of $\operatorname{St} \left( \mathcal{I} \right)/N$.

Let $\height(\alpha)$ denote the height of the root $\alpha$ with respect to a fixed set of simple roots (i.e. if $\mathring \Pi = \left\{ \alpha_{1},\dots,\alpha_{n} \right\}$ and $\alpha = \sum_{i=1}^n \lambda_{i} \alpha_{i}$ then $\height(\alpha) = \sum_{i=1}^n \lambda_{i}$). Note that if $\alpha\in \mathring \Phi^+$ then $\height(\alpha) >0$ and if $\alpha \in \mathring \Phi^-$ then $\height(\alpha)<0$ (and the height is always an integer).

We define the following map: 

\begin{align*}
\xi: \operatorname{St} \left( K[t]/(f) \right) &\to \operatorname{St} \left( \mathcal{I} \right)/N \\
x_{\alpha}(g+(f)) &\mapsto \begin{cases}
\psi(x_{\alpha}(t^{-\height(\alpha)}g + (f))) & \text{ if } \alpha \in \mathring \Phi^+ \\ 
\psi(x_{\alpha}(t^{-\height(\alpha)}g + (tf))) & \text{ if } \alpha \in \mathring \Phi^- 
\end{cases}.
\end{align*}

Note that since $K[t]/(f)$ is a field, it is well defined to consider $t^{-1}$ modulo $(f)$. 
We need to check that $\xi$ is a well-defined homomorphism.
Let $g \in (f)$. Then $t^{-\height(\alpha)}g \in (f)$ if $\alpha$ is positive and $t^{-\height(\alpha)}g \in (tf)$ if $\alpha$ is negative. 
Hence different representatives of the same class $g+(f)$ get mapped to the same element. 
Next, we check that the relations are respected by the map.

Clearly, $\xi(x_{\alpha}(s+u)) = \xi(x_{\alpha}(s))\xi(x_{\alpha}(u))$.

Note that $\height(i\alpha + j\beta) = i\height(\alpha)+j\height(b)$ for all $\alpha,\beta \in \mathring\Phi, i,j \in \mathbb{Z}$ and hence
$$
[x_{\alpha}(t^{-\height(\alpha)}s), x_{\beta}(t^{-\height(\beta)}u)]=\prod_{i\alpha + j \beta \in \mathring \Phi, i,j>0} x_{i\alpha+j\beta}(C_{i,j}^{{\alpha,\beta}}t^{-\height(i\alpha+j\beta)}s^i u^j).
$$
This implies that the commutator relation is also preserved under $\xi$.

Hence, $\xi$ is a well-defined homomorphism. 
The next step is to show that it is surjective. 
This follows from the fact that 
$$
K[t]/(f) \to K[t]/(f): \lambda \mapsto (t+(f))^n \lambda
$$
is bijective for all $n \in \mathbb{Z}$ and that the map
$$
K[t]/(f) \to (t)/(tf) : \lambda \mapsto t \lambda
$$
is well defined and bijective. 
Hence all the generators of $\operatorname{St} \left( \mathcal{I} \right)/N$ lie in the image of $\xi$.

On the other hand, since $N \subseteq M$, the surjective map $\operatorname{St} \left( \mathcal{I} \right) \to \operatorname{St} \left( K[t]/(f) \right)$ induces a surjective map $\operatorname{St} \left( \mathcal{I} \right)/N \to \operatorname{St} \left( K[t]/(f) \right)$.
Since $\operatorname{St} \left( K[t]/(f) \right)$ is finite and $\operatorname{St} \left( K[t]/(f) \right) \cong \operatorname{St} \left( \mathcal{I} \right)/M$ we get that $\operatorname{St} \left( \mathcal{I} \right)/N \cong \operatorname{St} \left( \mathcal{I} \right)/M$, which concludes Step 1.

Step 2: 
Since $\operatorname{St} \left( \mathcal{I} \right)/M \cong \operatorname{St} \left( K[t]/(f) \right) \cong \chev(K[t]/(f))$ is finite,  we can conclude from Step 1 that $N=M$.

To show this, note that we have $N,M \lhd \operatorname{St} \left( \mathcal{I} \right), N \subseteq M$ thus $M/N \lhd \operatorname{St} \left( \mathcal{I} \right)/N$. Hence by the third isomorphism theorem $(\operatorname{St} \left( \mathcal{I} \right)/N )/(M/N) \cong \operatorname{St} \left( \mathcal{I} \right)/M$. Hence $\lvert \operatorname{St} \left( \mathcal{I} \right)/N \rvert / \lvert M/N \rvert = \lvert \operatorname{St} \left( \mathcal{I} \right)/M \rvert$. But $\lvert \operatorname{St} \left( \mathcal{I} \right)/N \rvert = \lvert \operatorname{St} \left( \mathcal{I} \right)/M \rvert$, thus $\lvert M/N \rvert = 1$ which implies $N=M$.

\end{proof}

\begin{corollary} \label{cor: trivial cohom of kernel}
Assume the setting from above. Let $p=\operatorname{char}(K)$ and let $\Lambda$ be a group with no non-trivial elements of order $p$. Then $H^1(\ker(\varphi_f),\Lambda) = 0$, where the cohomology is considered with trivial action of $\ker(\varphi_f)$ on $\Lambda$.
\end{corollary}
\begin{proof}
By \Cref{thm: kernel gen by order p elements}, we have that $\chev(K[t],(f))\cap \mathcal{I}$ is normally generated by order $p$ elements, and hence the same holds for $B:=\ker(\varphi_f)$ (since the two groups are isomorphic because $\phi$ is injective). 
Since $\Lambda$ has no non-trivial elements of order $p$ the only possible morphism in $\hom_{\operatorname{Grp}}(B,\Lambda)$ is the trivial one. 
This implies that $H^1(B,\Lambda)$ is trivial, since the elements of the first non-Abelian cohomology group are precisely conjugacy classes of the elements of $\hom_{\operatorname{Grp}}(B,\Lambda)$.
\end{proof}

The vanishing of first cohomology of congruence KMS complexes is summarized in the following corollary. 
\begin{corollary} \label{cor: homology of KMS complex is trivial}
Let $\mathring \Phi$ by a spherical root system of rank $n\geq 4$ and $K$ a finite field with $\lvert K \rvert \geq 5$ and $\operatorname{char}(K) = p \neq 2$. Let $f \in K[t]$ be an irreducible polynomial of degree at least 2. Let $X(\mathring \Phi, K, f)$ denote the corresponding congruence KMS complex as in \Cref{def:congKMS}.  Furthermore, let $\Lambda$ be a group with no non-trivial elements of order $p$. Then 
$$H^1(X(\mathring \Phi, K, f),\Lambda)=0.$$ 
\end{corollary}

\begin{proof}
This follows from \Cref{thm: homology of quotients} together with \Cref{prop: simple connectedness}, \Cref{lemma: quotients have correct form} and \Cref{cor: trivial cohom of kernel}.
\end{proof}

\begin{remark}
A similar argument, using \Cref{thm: homology of quotients} and \Cref{prop: simple connectedness}, could be applied to other finite quotients of KMS groups used in the KMS complex construction as well, provided that one can show that the cohomology of the kernel of the quotient is trivial.
\end{remark}

\section{Non-Abelian cone functions}
\label{sec: n-a cone functions}
In the following section, we introduce non-Abelian cone functions, a tool to bound the coboundary expansion constants of strongly symmetric simplicial complexes. The definition is very similar to the work of Dinur and Dikstein in \cite{DDswap} and we claim no originality here. Our contribution in this chapter is the description of new ways to construct cone functions with bounded radius in \Cref{subsec: cone functions of joins} and \Cref{subsec: cone functions adding vertices}. These constructions are very similar to the ones for Abelian cone functions in \cite[Chapter 2.5, Chapter 3]{CSEofKMS}.

\subsection{Definition and relation to coboundary expansion}
\label{subsec: def cone and relation CBE}
Throughout, $X$ is an $n$-dimensional simplicial complex.

Given two vertices $u,v$ of $X$,  a \textit{path from $u$ to $v$} is a sequence of $k$-tuple of vertices $(u_0 ;... ; u_k)$ where $u = u_0$ and $v =u_k$ such that for every $0 \leq i \leq k-1$,  $(u_i,u_{i+1}) \in X_\ord (1)$ (we use the $;$ symbol to distinguish between the notation of a path and that of an ordered simplex).  The constant $k$ is called the \textit{length of the path}.  Given two paths $P_0  =(u_0; ... ;u_k),  P_1 = (u_k ; ... ;u_m)$,  we define $P_0 \circ P_1 = (u_0; ... ;u_m)$.  Also,  given a path $P = (u_0; ... ;u_k)$, we denote $P^{-1} = (u_k ;...; u_0)$.

A \textit{loop around $v$} is a path from $v$ to itself.  We note that $(v)$ is also considers as a loop and we will call it the \textit{trivial loop around $v$}.  Fix $v_0$ to be a vertex in $X$.  We define the following symmetric relations between two loops $P,  P '$ around $v_0$:
\begin{itemize}
\item We denote $P \mathop{\sim}\limits^{(BT)} P '$ if there exist vertices $u,v$ and paths $Q_0, Q_1$ from $v_0$ to $u$ such that $P = Q_0 \circ (u ; v; u) \circ Q_1^{-1}$ and $P ' = Q_0 \circ (u) \circ Q_1^{-1}$ ($(BT)$ stands for ``backtracking'').
\item We denote $P \mathop{\sim}\limits^{(TR)} P '$ if there is a triangle $(u,v,w)$ and paths $Q_0, Q_1$ from $v_0$ to $u$ and $w$ correspondingly such that $P = Q_0 \circ (u ; v; w) \circ Q_1^{-1}$ and $P '= Q_0 \circ (u ; w) \circ Q_1^{-1}$ ($(TR)$ stands for ``triangle'').
\end{itemize}

\begin{remark}
Note that the relations above are not defined to be transitive.  Indeed,  in the sequel,  it will be important to keep track of how many times the relation $\mathop{\sim}\limits^{(TR)} $ is used to pass from one loop to another.
\end{remark}

\begin{definition}
A \textit{$0$-cone} on $X$ is a pair  $\Cone_X = (v_0,  \lbrace P_u \rbrace_{u \in X(0)})$ such that
\begin{enumerate}
\item $v_0 \in X(0)$.
\item For every $u \in X(0), u \neq v_0$,  $P_u$ is a path from $v_0$ to $u$.  For $v_0$,  $P_{v_0}$ is the trivial loop around $v_0$.
\end{enumerate}

A \textit{$1$-cone} on $X$ is a triple $\Cone_X = (v_0,  \lbrace P_u \rbrace_{u \in X(0)},  \lbrace T_{(u,v)} \rbrace_{(u,v) \in X_\ord (1)} )$ such that
\begin{enumerate}
\item $v_0 \in X(0)$.
\item For every $u \in X(0), u \neq v_0$,  $P_u$ is a path from $v_0$ to $u$.  For $v_0$,  $P_{v_0}$ is the trivial loop around $v_0$.
\item For every $(u,v) \in X_{\ord} (1)$,  $T_{(u,v)}$ is a sequence of loops around $v_0$,  $T_{(u,v)} = (P_0,P_1,\dots,P_m)$ such that the following holds:
\begin{enumerate}
\item $P_0 = P_u \circ (u ; v) \circ P_v^{-1}$.
\item For every $0 \leq i \leq m-1$,  either $P_i \mathop{\sim}\limits^{(TR)}  P_{i+1}$ or $P_i \mathop{\sim}\limits^{(BT)} P_{i+1}$.
\item $P_m$ is the trivial loop around $v_0$.
\end{enumerate}
\end{enumerate}

We call $T_{(u,v)} = (P_0,P_1,...,P_m)$ as above a \textit{contraction} and denote 
$$\vert T_{(u,v)} \vert = \vert \lbrace 0 \leq i \leq m-1 : P_i \mathop{\sim}\limits^{(TR)}  P_{i+1} \rbrace \vert . $$  
We further define the \textit{radii of the cone} as follows:
\begin{enumerate}
\item $\Rad_0 (\Cone_X) = \sup_{u \in X(0)} \operatorname{length}(P_u)$. 
\item $\Rad_1(\Cone_X) = \sup_{(u,v) \in X_{\ord} (1)} \vert T_{(u,v)} \vert$.
\item $\Rad(\Cone_X) = \begin{cases}\max \{\Rad_0(\Cone_X), \Rad_1(\Cone_X) \} & \text{ if } \Cone_X \text{ is a 1-cone} \\ \Rad_0(\Cone_X) &  \text{ if } \Cone_X \text{ is a 0-cone}  \end{cases}$.
\end{enumerate}
\end{definition}

\begin{remark}
The definition of $\vert T_{(u,v)} \vert$ in \cite{DDswap} is slightly different than ours.  Namely,  they define another relation $\sim_1$ that is composed of applying $\mathop{\sim}\limits^{(TR)}$ once and $\mathop{\sim}\limits^{(BT)}$ as many times as needed and define  $\vert T_{(u,v)} \vert$ as how many such relations $\sim_1$ are needed to bring $P_0$ to $P_m$ such that $P_m$ can be made trivial by only applying only the relations $\mathop{\sim}\limits^{(BT)}$.  We leave it to the reader to verify that both definitions coincide.
\end{remark}

We denote by $\Aut (X)$ the group of simplicial automorphisms of $X$.  The complex $X$ is called \textit{strongly symmetric} if $\Aut (X)$ acts transitively on $X (n)$,  i.e.,  for every $\sigma_1,  \sigma_2 \in X(n)$ there is $g \in \Aut (X)$ such that $g. \sigma_1 = \sigma_2$.   In \cite{DDswap},  the following result was proven:
\begin{lemma}\cite[Lemma 1.6,  Lemma 4.3]{DDswap} \label{lemma: cone implies CBE}
Let $X$ be a pure $n$-dimensional finite simplicial complex, $n\geq 2$.  Assume that there is a 1-cone  on $X$ such that $\Rad_i (\Cone_X) < \infty, i=0,1$ and that $X$ is strongly symmetric,  then for every group $\Lambda$,
$$h^1_{\cobound} (X,  \Lambda) \geq \frac{1}{{n+1 \choose 3} \Rad_1 (\Cone_X)}.$$
\end{lemma}

\subsection{Cone functions of joins}
\label{subsec: cone functions of joins}
We denote by $\sqcup$ the disjoint union of sets.  Given two non-empty finite simplicial complexes $Y_1$ and $Y_2$ with vertex sets $V(Y_1)$ and $V(Y_2)$,  the \textit{join} of $Y_1$ and $Y_2$,  denoted $Y_1 * Y_2$,  is the simplicial complex with the vertex set $V (Y_1) \sqcup V(Y_2)$ and simplices $\sigma \sqcup  \tau$ for every $\sigma \in Y_1$ and $\tau \in Y_2$.  We note that  $Y_1 * Y_2$ includes all the simplices of the form $\tau = \emptyset \sqcup \tau$ and $\sigma = \sigma \sqcup \emptyset$ for every $\sigma \in Y_1$ and $\tau \in Y_2$.

\begin{proposition}
\label{cone of a join - single vertex prop}
Let $Y$ be a non-empty simplicial complex and $w$ a vertex not in $Y$. Then there exists a 1-cone function $\Cone_{\{w\}*Y}$ with $\Rad_i(\Cone_{\lbrace w \rbrace * Y}) \leq 1$ for $i = 0,1$.  
\end{proposition}

\begin{proof}
Set $X = \{w\} * Y$. First, assume that $\dim Y \geq 1$. We note that for every $(u,v) \in Y_{\ord} (1)$ it holds that $(w,u,v) \in X_{\ord} (2)$ and thus if $P$ is the trivial loop around $w$, and $P_u,P_v$ are two paths starting in $w$ and ending in $u,v$ respectively, then
$$P_u \circ (u ; v) \circ P_v^{-1} \mathop{\sim}\limits^{(TR)} P.$$
Thus,  we define a function $\Cone_{\lbrace w \rbrace * Y} = (v_0,  \lbrace P_u \rbrace_{u \in X(0)},  \lbrace T_{(u,v)} \rbrace_{(u,v) \in X_\ord (1)} )$ as follows: 
\begin{itemize}
\item $v_0 = w$.
\item For every $u \in X(0)$,  $P_u = (w ; u)$.
\item For every $(u,v) \in X_{\ord} (1)$,  $T_{u,v} = (P_0 = P_u \circ (u ; v) \circ P_v^{-1},  P_1 = P_{w})$.
\end{itemize}
It follows that $\Rad_i(\Cone_{\lbrace w \rbrace * Y}) =1, i=0,1$ as needed.

If $\dim Y=0$ then $\dim X=1$. But we can still define a 1-cone function of radius 1: The $v_0$ and $P_u$ are defined as before. The only edges in $X_{\ord}(1)$ are of the form $(w,u)$ (or $(u,w)$) with $u \in Y(0)$. Thus $T_{(w,u)}$ consists of $P_0 = (w) \circ (w,u) \circ (u,w) \sim_{(BT)} P_1 = (w)$ (and $T_{(u,w)}$ of $P_0 = (w,u) \circ (u,w) \circ (w) \sim_{(BT)} (w)$).  
\end{proof}

\begin{remark}
We note that in the above construction,  the cone radius is still $\leq 1$, even if $Y$ has no edges or if $Y$ is disconnected.
If $Y$ is empty, we can still say that $\lbrace w \rbrace * \emptyset$ has a 0-cone function, since it will be a connected graph consisting of one single vertex. 
\end{remark}

\begin{lemma} \label{lemma: join of two scs without assumption}
	Let $Y_1,Y_2$ be non-empty, finite simplicial complexes. Then $Y_1*Y_2$ admits a 0-cone function with $\Rad_0(\Cone_{Y_1*Y_2}) =2$.
\end{lemma}
\begin{proof}
	We have that the 1-skeleton of $Y_1*Y_2$ contains the complete bipartite graph on $Y_1(0) \sqcup Y_2(0)$ as a subgraph. In particular, every two vertices are connected by a path of length at most 2.
\end{proof}

\begin{proposition} \label{prop: join given one 0-cone}
	Let $Y_1$ be a $n$-dimension simplicial complex, $n \geq 1$, with 0-cone function $\Cone_{Y_1} = (v^1_0, (P^1_u)_{u \in Y_1(0)})$. Let $Y_2$ be a non-empty, finite simplicial complex. Then $Y_1*Y_2$ admits a 1-cone function $\Cone_{Y_1*Y_2}$ with $\Rad_0(\Cone_{Y_1*Y_2}) = \Rad_0(\Cone_{Y_1})$ and $\Rad_1(\Cone_{Y_1*Y_2}) \leq 2 \Rad_0(\Cone_{Y_1}) +1$. 
\end{proposition}

\begin{proof}
	We want to define a 1-cone $(v_0, (P_u)_{u \in Y_1(0) \sqcup Y_2(0)}, (T_{(u,v)})_{(u,v) \in (Y_1*Y_2)_\ord(1)} )$ for $Y_1*Y_2$. First, note that $(Y_1*Y_2)_\ord(1) = {Y_1}_\ord(1) \sqcup {Y_2}_\ord(1) \sqcup \{(u,v) \mid u \in Y_i, v \in Y_j, \{i,j\} = \{1,2\}\}$.
	\begin{itemize}
		\item We set $v_0 = v_0^1$.
		\item For $u \in Y_1(0)$, we set $P_u = P_u^1$.
		\item For $u \in Y_2(0)$, we set $P_u = (v_0,u)$.
		\item For $(u,v) \in {Y_2}_\ord(1)$, we set $T_{(u,v)} = (P_0 = (v_0;u;v;v_0), P_1 = (v_0;v;v_0), P_2 = (v_0))$ with $P_0 \sim_{(TR)} P_1$ and $P_1 \sim_{(BT)} P_2$.

	 \item Fix $u \in Y_1(0), v \in Y_2(0)$ and let $P_u = (v_0;u_1; \dots; u_n; u)$. Then we define $T_{(u,v)}$ as follows, see also \Cref{fig: prop join given one 0 cone case 2}. 
	\begin{align*}
		P_0 = (v_0;u_1; \dots; u_{n}) &\circ (u_n,u) \circ (u,v) \circ (v,v_0) \\
		&\wr_{(TR)} \text{with triangle } (u_n,u,v)\\
		P_1=(v_0;u_1; \dots; u_{n}) &\circ (u_n,v) \circ (v,v_0) \\
		&\wr_{(TR)} \\
		&\vdots \\
		&\wr_{(TR)} \\
		P_n=(v_0,u_1) \circ &(u_1,v) \circ (v,v_0) \\
		&\wr_{(TR)} \text{with triangle } (v_0,u_1,v)\\
		(v_0,v) \circ & (v,v_0)\\
		&\wr_{(BT)} \\
		&(v_0)
	\end{align*}
	Observe that $\lvert T_{(u,v)} \rvert = \operatorname{length}(P_u)$.
	\begin{figure}[h] 
		\centering
			\tikzset{
  offset round/.code=
    \tikz@addmode{%
      \pgfsetroundjoin         \pgfgetpath   \tikz@temp
      \pgfsetpath\pgfutil@empty\pgfoffsetpath\tikz@temp{#1}}}
		\begin{tikzpicture}
  \draw[Mycolor1, ->] plot[samples at={(1,0),(1,-1),(1,-2),(1,-3),(-1,-1.5), (1,0)}] (\x) [offset round=+2.5mm];
    \draw[Mycolor2, ->] plot[samples at={(1,0),(1,-1),(1,-2),(-1,-1.5), (1,0)}] (\x) [offset round=+2mm];
	\draw[Mycolor3, ->] plot[samples at={(1,0),(1,-1),(-1,-1.5), (1,0)}] (\x) [offset round=+1.5mm];
	  \node (a) at (-1,-1.5) {$v$};
  \node (b) at (1,0) {$v_0$};
  \node (c) at (1,-1) {$u_1$};
  \node (e) at (1,-2) {$u_{n}$};
    \node (f) at (1,-3) {$u$};
	\node[Mycolor1] (n) at (-1.55,0) {$P_0,$};
	\node[Mycolor2] (o) at (-1,0) {$P_1,$};
	\node[Mycolor3] (p) at (-0.45,0) {$P_n$};
  \graph { (a) -- {(b), (c),  (e), (f)}; (b) -- (c) --[dotted] (e) -- (f) };
\end{tikzpicture}
	\caption{Visualization of the paths in $T_{(u,v)}$ for $u \in Y_1(0), v \in Y_2(0)$.}
	\label{fig: prop join given one 0 cone case 2}
	\end{figure}
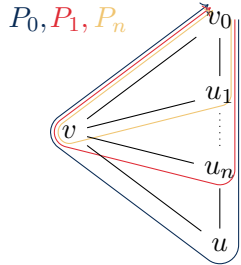

	\item Lastly, fix $(u,v) \in {Y_1}_\ord(1)$ and denote $P_u = (v_0;u_1; \dots; u_n; u), P_v = (v_0;v_1;\dots; v_n;v)$. Then we define $T_{(u,v)}$ as follows, similar to the above case, but now doing the contractions on ``both sides'' of the path. We also need to fix an arbitrary vertex $x \in Y_2(0)$. See also \Cref{fig: prop join given one 0 cone case 3}.
	\begin{align*}
		P_0 = (v_0;u_1; \dots; u_{n}) &\circ (u_n,u) \circ (u,v) \circ (v,v_n) \circ (v_n; \dots; v_1;v_0) \\
		&\wr_{(TR)} \text{ with triangle } (u,x,v) \\
		P_1 = (v_0;u_1; \dots; u_{n}) &\circ (u_n,u) \circ (u,x) \circ (x,v) \circ (v,v_n) \circ (v_n; \dots; v_1;v_0) \\
		&\wr_{(TR)} \text{with triangle } (u_n,u,x) \text{ and } (v,v_n,x)\\
		P_3 = (v_0;u_1; \dots; u_{n}) &\circ (u_n,x) \circ (x,v_n) (v_n; \dots;v_0) \\
		& \vdots \\
		(v_0,u_1) \circ &(u_1,x) \circ (x,v_1) \circ (v_1,v_0) \\
		&\wr_{(TR)} \text{with triangle } (v_0,u_1,x) \text{ and } (x,v_1,v_0)\\
		(v_0,x) \circ & (x,v_0)\\
		&\wr_{(BT)} \\
		&(v_0)
	\end{align*}
	\begin{figure}[h] 
\centering
\begin{tikzpicture}
		\tikzset{
  offset round/.code=
    \tikz@addmode{%
      \pgfsetroundjoin         \pgfgetpath   \tikz@temp
      \pgfsetpath\pgfutil@empty\pgfoffsetpath\tikz@temp{#1}}}  
  \draw[Mycolor1, ->] plot[samples at={(1,3),(-0.5,1.5),(0,0),(2,0),(2.5,1),(2.5,2),(1,3)}] (\x) [offset round=-2.5mm];
    \draw[Mycolor2, ->] plot[samples at={(1,3),(-0.5,1.5),(0,0),(1,1.5),(2,0),(2.5,1),(2.5,2),(1,3)}] (\x) [offset round=-2mm];
	\draw[Mycolor3, ->] plot[samples at={(1,3),(-0.5,1.5),(1,1.5),(2,0),(2.5,1),(2.5,2),(1,3)}] (\x) [offset round=-1.5mm];
	\draw[Mycolor4, ->] plot[samples at={(1,3),(-0.5,1.5),(1,1.5),(2.5,1),(2.5,2),(1,3)}] (\x) [offset round=-1mm];
	\node (a) at (2,0) {$v$};
  \node (b) at (0,0) {$u$};
  \node (c) at (1,3) {$v_0$};
  \node (d) at (-0.5,1.5) {$u_{1}$};
  \node (e) at (2.5,2) {$v_{1}$};
    \node (f) at (2.5,1) {$v_{2}$};
    \node (g) at (1,1.5) {$x$};
	\node[Mycolor1] (n) at (-1.5,3) {$P_0,$};
	\node[Mycolor2] (o) at (-1,3) {$P_1,$};
	\node[Mycolor3] (p) at (-0.5,3) {$P_2,$};
	\node[Mycolor4] (q) at (0,3) {$P_3$};
  \graph { (g) -- {(a),(b), (c),(d), (e), (f)}; (a)--(b) -- (d) --(c); (a) -- (f) -- (e) -- (c) };
\end{tikzpicture}
\caption{Visualization of the paths in $T_{(u,v)}$ for $(u,v) \in {Y_1}_{\ord}(1)$ (for an example with $\operatorname{length}(P_u)=2, \operatorname{length}(P_v)=3$).}
\label{fig: prop join given one 0 cone case 3}
	\end{figure}
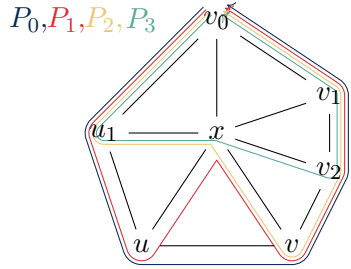
	Observe that the paths $P_u, P_v$ don't necessary have to be of the same length. In that case, the induction on one side halts earlier then on the other side. Counting the total number of $(TR)$ operations made, we have $\lvert T_{(u,v)} \rvert = \operatorname{length}(P_u) + \operatorname{length}(P_v) +1$.
\end{itemize}
\end{proof}

\begin{proposition}
\label{cone of a join - general prop}
Let $Y_1,  Y_2$ be finite, non-empty, simplicial complexes. Assume that there is a $1$-cone function $\Cone_{Y_1}$ with finite cone radius. Then there is a $1$-cone function $\Cone_{Y_1*Y_2}$ such that,  
\begin{align*}
\Rad_{0} (\Cone_{Y_1*Y_2})  & = \Rad_0(Y_1) \\
\Rad_1(\Cone_{Y_1*Y_2}) &= \max \{\Rad_0(Y_1), \Rad_1(Y_1),1\}.
\end{align*}
\end{proposition}

\begin{proof}
Recall that $(Y_1*Y_2)(0) = Y_1(0) \sqcup Y_2(0)$ and that any simplex of either $Y_i$ can be seen as a simplex in the join as well. In particular, $(Y_1*Y_2)_\ord(1) = {Y_1}_\ord(1) \sqcup {Y_2}_\ord(1) \sqcup \{(u,v) \mid u \in Y_i(0), v \in Y_j(0), i,j \in \{1,2\}, i \neq j \}$.
Let $\Cone_{Y_1} = (v_0^1, (P_u^1)_{u \in Y_1(0)}, (T_{(u,v)}^1)_{(u,v)\in {Y_1}_\ord(1))}$. We will define $\Cone_{Y_1 * Y_2} = (v_0,  \lbrace P_u \rbrace_{u \in Y_1*Y_2(0)},  \lbrace T_{(u,v)} \rbrace_{(u,v) \in (Y_1*Y_2)_\ord (1)} )$ as follows. 
\begin{itemize}
	\item We set $v_0 = v_0^1$.
	\item For $u \in Y_1(0)$, set $P_u = P_u^1$.
	\item For $v \in Y_2(0)$, set $P_v = (v_0,v)$. (Note that this is well defined since $v_0 \in Y_1$ is connected by an edge to any vertex of $Y_2$.)
	\item For $(u,v) \in {Y_1}_\ord(1)$, set $T_{(u,v)} = T_{(u,v)}^1$.
	\item For $(u,v) \in {Y_2}_\ord(1)$, set $T_{(u,v)} = (P_0 = (v_0,u)\circ (u,v) \circ (v,v_0), P_1 = (v_0,v) \circ (v,v_0), P_2 = (v_0))$. Note that $P_0 \sim_{(TR)} P_1$ and $P_1 \sim_{(BT)} P_2$.

\item
\begin{figure}[h] 
\centering
\begin{tikzpicture}
			\tikzset{
  offset round/.code=
    \tikz@addmode{%
      \pgfsetroundjoin         \pgfgetpath   \tikz@temp
      \pgfsetpath\pgfutil@empty\pgfoffsetpath\tikz@temp{#1}}} 
  \draw[Mycolor1, ->] plot[samples at={(1,0),(1,-1),(1,-2),(1,-3),(-1,-1.5), (1,0)}] (\x) [offset round=+2.5mm];
    \draw[Mycolor2, ->] plot[samples at={(1,0),(1,-1),(1,-2),(-1,-1.5), (1,0)}] (\x) [offset round=+2mm];
	\draw[Mycolor3, ->] plot[samples at={(1,0),(1,-1),(-1,-1.5), (1,0)}] (\x) [offset round=+1.5mm];
	  \node (a) at (-1,-1.5) {$v$};
  \node (b) at (1,0) {$v_0^1$};
  \node (c) at (1,-1) {$u_1$};
  \node (e) at (1,-2) {$u_{n}$};
    \node (f) at (1,-3) {$u$};
		\node[Mycolor1] (n) at (-1.5,0) {$P_0,$};
	\node[Mycolor2] (o) at (-1,0) {$P_1,$};
	\node[Mycolor3] (p) at (-0.5,0) {$P_n$};
  \graph { (a) -- {(b), (c),  (e), (f)}; (b) -- (c) --[dotted] (e) -- (f) };
\end{tikzpicture}
\caption{Visualization of the paths in $T_{(u,v)}$ for $u \in Y_1(0), v \in Y_2(0)$.}
\label{fig: prop cone of join general}
\end{figure}
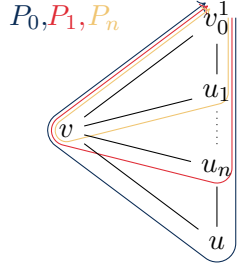

Lastly, we define the most involved case, where $u \in Y_1(0), v \in Y_2(0)$. We have $P_0 = P_u \circ (u,v) \circ (v,v_0)$.
We give names to the vertices in $P_u$: $P_u = (v_0;u_1;u_2; \dots; u_n; u)$. Note that each $u_i$ is also connected to $v$ by an edge. Next we define $P_1 = (v_0;u_1; \dots; u_n) \circ (u_n,v) \circ (v,v_0)$ and observe that $P_0 \sim_{(TR)} P_1$ via the triangle $(u_n,u,v)$.
We continue this process iteratively, see also \Cref{fig: prop cone of join general}: $P_i = (v_0;u_1;\dots; u_{n-i+1}) \circ (u_{n-i+1},v) \circ (v,v_0)$, hence $P_{i-1} \sim_{(TR)} P_i$ via the triangle $(u_{n-i+1}, u_{n-i +2},v)$ for $i = 1,\dots, n$. 
Then $P_n = (v_0;u_1;v;v_0)$, and we finish the contraction via $P_{n+1} = (v_0;v;v_0), P_{n+2} = (v_0)$ with $P_n \sim_{(TR)} P_{n+1}$ and $P_{n+1} \sim_{(BT)}P_{n+2}$.

In particular, we have used $(TR)$ once for every edge in $P_u$, hence $\lvert T_{(u,v)} \rvert = \operatorname{length}(P_u)$.

The definition of $T_{(v,u)}$ is analogues, just with reversed direction of the paths.

\end{itemize}
\end{proof}

To avoid having to make case distinctions between cone functions for $1$-dimensional complexes and complexes with dimension $\geq 2$, we introduce the following notation.

\begin{definition} \label{def: nac}
	Let $X$ be a finite simplicial complex. We say that $X$ has a non-Abelian cone function (NAC) if 
	\begin{itemize}
		\item $X \neq \emptyset$,
		\item if $\dim X = 1$ there exists a $0$-cone function for $X$,
		\item if $\dim X \geq 2$ there exists a 1-cone function for $X$.
	\end{itemize}
	Furthermore, we define the cone radius of $X$, denoted $\Rad(X)$ as 
	\begin{itemize}
		\item if $\dim X =0, X \neq \emptyset$ then $\Rad(X)=1$,
		\item if $\dim X \geq 1$ and $X$ has a NAC then $\Rad(X) = \min \{ \Rad(\Cone_X) \mid \Cone_X \text{ a 1-cone function of }X \}$,
		\item if $X$ does not admit a NAC, we set $\Rad (X) = \infty$.
	\end{itemize}

\end{definition}
 With this notation, we can summarize the results of this chapter. By summarizing the different cases, we loose some sharpness of the bounds on the cone radius. 
 \begin{corollary}
	Let $X,Y$ be two finite simplicial complexes. If both $X$ and $Y$ have a NAC then $X * Y$ has a NAC and $\Rad(X*Y)\leq 2 \max \{\Rad(X),\Rad(Y)\} +1$. 
 \end{corollary} \label{cor: summary of join results}
 \begin{proof}
	If $\dim X = \dim Y =0$ then $\dim X * Y=1$ and it has a 0-cone function by \Cref{lemma: join of two scs without assumption}. 
	If $\dim X = 1$ then the assumption implies that $X$ has a 0-cone function. Thus \Cref{prop: join given one 0-cone} gives the desired result. 
	If $\dim X \geq 2$, then by assumption there exists a 1-cone function for $X$. Hence \Cref{cone of a join - general prop} gives the desired result. 
 \end{proof}

\subsection{Constructing cone functions by adding vertices}
\label{subsec: cone functions adding vertices}

The idea is to construct a cone function for a given complex $X$, by starting with a given subcomplex for which we have a cone function and adding a set of vertices which is ``well-behaved''.  The exact formalism of this idea is given in Theorem \ref{adding vertices to cone thm} that can be thought of as a special version of the Mayer-Vietoris theorem for cone functions.  

\begin{definition}
For a complex $X$,  a subcomplex $X' \subseteq X$ is called full,  if for every $v_0,\dots,v_k \in X'$,  if $\lbrace v_0,\dots,v_k \rbrace \in X$, then $\lbrace v_0,\dots,v_k \rbrace \in X'$.
\end{definition}

\begin{theorem}
\label{adding vertices to cone thm}
Let $X$ be an $n$-dimensional simplicial complex, $n\geq 2$, and let $X '$ be a full subcomplex.  Assume there is a $1$-cone function $\Cone_{X'}$ and constants $R_j' \in \mathbb{N},   j=0,1 $ such that for every $j=0,1$,  $\Rad_j (\Cone_{X'})  \leq R_j '$. 
Also let $W \subseteq X(0)$ be a set of vertices such that the following holds:
\begin{enumerate}
\item $W \cap X' = \emptyset$. 
\item For every $w_1,  w_2 \in W$,  $\lbrace w_1, w_2 \rbrace \notin X (1)$.  
\item For every $w \in W$,  $\lk_X(w) \cap X'$ is a non-empty simplicial complex.
\item There are $0$-cone functions $\Cone_{\lk_X(w) \cap X'}$ for all $w \in W$ and a constant $R'' \in \mathbb{N}$ such that for every $w \in W$, $\Rad_0  (\Cone_{\lk_X(w) \cap X'})  \leq R ''.$
\end{enumerate} 
Let $X' \cup W$ be the full subcomplex of $X$ spanned by $X' (0) \cup W$.  Then there is a $1$-cone function $\Cone_{X' \cup W}$ such that $\Rad_0 (\Cone_{X' \cup W}) \leq R'_0 +1, \Rad_1 (\Cone_{X' \cup W}) \leq R'' (R_1 '+1)$.
\end{theorem}

\begin{proof}
	First, we fix some notation: for each $w \in W$ let $L_w := \lk_X(w) \cap X'$ and
	\begin{align*}
		\Cone_{X'} &= (v_0', (P_u')_{u \in X'(0)}, (T'_{(u,v)})_{(u,v)\in X'_\ord(1)})\\
		\Cone_{L_w} &= (v_0^w,(P_u^w)_{u \in L_w(0)}).
	\end{align*}
	Observe that $L_w(0)$ are those vertices of $X'$ which are neighbors of $w$, that $(X' \cup W) (0)=X'(0) \cup W$ and that edges in $X' \cup W$ are all edges of $X'$ and $\{\{u,w\} \mid w \in W, u \in L_w \}$.

	We define $\Cone_{X' \cup W} = (v_0, (P_u)_{u \in X'(0)\cup W}, (T_{(u,v)})_{(u,v)\in (X' \cup W)_\ord(1)})$ as follows.
	\begin{itemize}
		\item We set $v_0 = v_0' \in X'$.
		\item For $u \in X'(0)$, we set $P_u = P_u'$.
		\item For $w \in W$, we set $P_w = P'_{v_0^w} \circ (v_0^w,w)$.
		\item For $(u,v) \in X'_\ord (1)$, we set $T_{(u,v)} = T'_{(u,v)}$.
	
\item Lastly, we define the most involved case. Let $w \in W, u \in L_w$ (hence $(u,w)\in (X' \cup W)_\ord(1))$. Then we define $T_{(u,w)}$ as follows (and $T_{(w,u)}$ similarly by reversing all paths). We have $P_0 = P_u' \circ (u,w) \circ (w,v_0^w) \circ (P_{v_0^w}')^{-1}$. We give names to the vertices in $P_u^w=(v_0^w; u_1; \dots; u_n;u)$. We get the following series of contractions, see also \Cref{fig: thm adding vertices}.
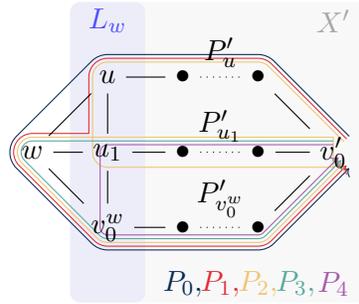
\begin{figure}[h!]
		\tikzset{
  offset round/.code=
    \tikz@addmode{%
      \pgfsetroundjoin         \pgfgetpath   \tikz@temp
      \pgfsetpath\pgfutil@empty\pgfoffsetpath\tikz@temp{#1}}} 
\centering	  
\begin{tikzpicture}
  \fill[gray!20,nearly transparent] (0.5,2) [rounded corners]-- (0.5,-2) [rounded corners]-- (4.5,-2) -- (4.5,2) -- cycle;
  \fill[blue!20,nearly transparent] (0.5,2) [rounded corners]-- (0.5,-2) [rounded corners]-- (1.5,-2) -- (1.5,2) -- cycle;

  \draw[Mycolor1, ->] plot[samples at={(4,0),(3,1),(2,1),(1,1),(0,0), (1,-1),(2,-1),(3,-1),(4,0)}] (\x) [offset round=-3mm];
  \draw[Mycolor2, ->] plot[samples at={(4,0),(3,1),(2,1),(1,1),(1,0),(0,0), (1,-1),(2,-1),(3,-1),(4,0)}] (\x) [offset round=-2.5mm];
  	\path plot[samples at={(4,0),(3,1),(2,1),(1,1),(1,0),(2,0),(3,0),(4,0)}] (\x) [offset round=- 2mm, spath/save=p];
	\path plot[samples at={(4,0),(3,0),(2,0),(1,0),(0,0),(1,-1),(2,-1),(3,-1),(4,0)}] (\x) [offset round=+-2mm, spath/save=q];
  \draw[->, Mycolor3, spath/use=p] -- (spath cs:q 0) [spath/use={q, weld}];
  \draw[Mycolor4, ->] plot[samples at={(4,0),(3,0),(2,0),(1,0),(0,0), (1,-1),(2,-1),(3,-1),(4,0)}] (\x) [offset round=-1.5mm];
   \draw[Mycolor5, ->] plot[samples at={(4,0),(3,0),(2,0),(1,0), (1,-1),(2,-1),(3,-1),(4,0)}] (\x) [offset round=-1mm];
    \node (a) at (0,0) {$w$};
  \node (b) at (1,1) {$u$};
  \node (c) at (1,0) {$u_1$};
  \node (d) at (1,-1) {$v_0^w$};
  \node (e) at (4,0) {$v_0'$};
    \node (f) at (2,1) {$\bullet$};
	\node (g) at (2,0) {$\bullet$};
	\node (h) at (2,-1) {$\bullet$};
	\node (i) at (3,1) {$\bullet$};
	\node (j) at (3,0) {$\bullet$};
	\node (k) at (3,-1) {$\bullet$};
	\node[blue!70] (l) at (1,1.75) {$L_w$};
	\node[gray!70] (m) at (4,1.75) {$X'$};
	\node[Mycolor1] (n) at (2,-1.75) {$P_0, $};
	\node[Mycolor2] (o) at (2.5,-1.75) {$P_1, $};
	\node[Mycolor3] (p) at (3,-1.75) {$P_2, $};
	\node[Mycolor4] (q) at (3.5,-1.75) {$P_3, $};
	\node[Mycolor5] (r) at (4,-1.75) {$P_4$};	
	
  \graph { (a) -- {(b), (c),  (d)}; (b) -- (f) --[dotted,"$P_u'$"] (i) -- (e); (c) -- (g) --[dotted,"$P_{u_1}'$"] (j) -- (e); (d) -- (h) --[dotted,"$P_{v_0^w}'$"] (k) -- (e); (b)--(c)--(d)   };
\end{tikzpicture}
\caption{Visualization of $T_{(u,w)}$ shown in the case $\operatorname{length}(P_u^w)=2$.}
\label{fig: thm adding vertices}
\end{figure}

\begin{align*}
	P_0 = P_u \circ (u,w) &\circ P_w^{-1} \\
	&\wr_{(TR)} \text{with triangle } (u,w,u_n)\\
	P_u \circ (u,u_n) &\circ (u_n,w) \circ P_w^{-1}  \\
	& \wr_{(BT)} \text{as often as necessary} \\
	P_u \circ (u,u_n) \circ P_{u_n}^{-1} & \circ P_{u_n} \circ (u_n,w) \circ P_w^{-1} \\
	&\wr \text{ apply sequence of relations from } T'_{(u,u_n)} \\
	(v_0) \circ P_{u_n} \circ &(u_n,w) \circ P_w^{-1} \\
	&\wr_{(TR)} \text{with triangle } (u_n,w,u_{n-1})\\
	P_{u_n} \circ (u_n,u_{n-1}) &\circ (u_{n-1},w) \circ P_w^{-1} \\
	& \wr_{(BT)} \text{as often as necessary} \\
	P_{u_n} \circ (u_{n},u_{n-1}) \circ P_{u_{n-1}}^{-1} & \circ P_{u_{n-1}} \circ (u_{n-1},w) \circ P_w^{-1} \\
	&\wr \text{ apply sequence of relations from } T'_{(u_n,u_{n-1})} \\
	(v_0) \circ P_{u_{n-1}} \circ &(u_{n-1},w) \circ P_w^{-1} \\
	&\vdots\\
	P_{u_1} \circ &(u_1,w) \circ P_w^{-1} \\
	&\wr_{(TR)} \text{with triangle } (u_1,w,v_0^w)\\
	P_{u_1} \circ (u_1,v_0^w) &\circ (v_0^w,w) \circ P_w^{-1}  \\
	& \wr_{(BT)} \text{as often as necessary} \\
	P_{u_1} \circ (u_1,v_0^w) \circ P_{v_0^w}^{-1} &P_{v_0^w} (v_0^w,w)P_w^{-1} \\
	&\wr \text{ apply sequence of relations from } T'_{(u_1,v_0^w)} \\
	(v_0) \circ P_{v_0^w} \circ &(v_0^w,w) \circ (w,v_0^w) P_{v_0^w}^{-1} \\
	& \wr_{(BT)} \text{as often as necessary} \\
	&(v_0)
\end{align*}
In total, for each edge in $P_u^w$ we have used at most $R_1' +1$ triangle relations. Hence in total $\lvert T_{(u,w)}\rvert \leq R''(R_1'+1)$.
\end{itemize}

\end{proof}

We need a slightly adapted version for 1-dimensional complexes. 
\begin{lemma}
\label{adding vertices to cone thm - dim 1}
Let $X$ be a one-dimensional simplicial complex, and let $X '$ be a full subcomplex.  Assume there is a $0$-cone function $\Cone_{X'}$ and constant $R' \in \mathbb{N}$ such that $\Rad_0 (\Cone_{X'})  \leq R'$. 
Also let $W \subseteq X(0)$ be a set of vertices such that the following holds:
\begin{enumerate}
\item $W \cap X' = \emptyset$. 
\item For every $w_1,  w_2 \in W$,  $\lbrace w_1, w_2 \rbrace \notin X (1)$.  
\item For every $w \in W$,  $\lk_X(w) \cap X'$ is non-empty.
\end{enumerate} 
Let $X' \cup W$ be the full subcomplex of $X$ spanned by $X' (0) \cup W$.  Then there is an $0$-cone function $\Cone_{X' \cup W}$ such that $\Rad_0 (\Cone_{X' \cup W}) \leq R' +1$
\end{lemma}
\begin{proof}
	Let $\Cone_{X'}=(v_0', (P_u')_{u \in X'(0)})$. We define $\Cone_{X' \cup W}=(v_0, (P_u)_{u \in X'(0)\cup W})$ as follows. 
	\begin{itemize}
		\item Set $v_0 = v_0'$.
		\item For $u \in X'(0)$, set $P_u = P_u'$.
		\item For $w \in W$, pick a vertex $v \in \lk_X(w) \cap X'$, and set $P_w = P_v \circ (v,w)$.
	\end{itemize}
\end{proof}

We can summarize both results in the following corollary. Note here that a 0-dimensional complex has a 0-cone function if and only if it contains exactly one vertex. 

\begin{corollary}\label{cor: adding vertices NAC version}
	Let $X$ be a finite simplicial complex of dimension at least 1 and let $X'$ be a full subcomplex of dimension at least 1. Assume that $X'$ has a NAC and $\Rad(X') \leq R'$ for some constant $R'\in \N$. Let $W \subseteq X(0)$ be a set of vertices and $R''\in \N$ such that
	\begin{enumerate}
\item $W \cap X' = \emptyset$. 
\item For every $w_1,  w_2 \in W$,  $\lbrace w_1, w_2 \rbrace \notin X (1)$.  
\item For every $w \in W$,  $\lk_X(w) \cap X'$ has a 0-cone function $\Cone_{\lk_X(w)\cap X'}$with $\Rad_0(\Cone_{\lk_X(w)\cap X'}) \leq R''$.
	\end{enumerate}
	Let $X' \cup W$ be the full subcomplex of $X$ spanned by $X' (0) \cup W$. Then $X'\cup W$ has a NAC and $\Rad(X'\cup W) \leq R''(R'+1)$.
\end{corollary}

\section{Cosystolic expansion with non-Abelian coefficients}
\label{sec: CSE with n-a coefficients}

In this section, we prove the following theorem which is the last key step in proving that classical congruence KMS complexes give rise to infinite families of bounded degree 1-coboundary expanders.

\begin{theorem} \label{thm: proper links are CBE}
	Let $n\geq 3$ be an integer. There exists an $\varepsilon >0$ such that for all prime powers $q > 2^{2n-1}$, every $n$-dimensional, $n$-classical KMS complex $X$ over $\mathbb{F}_q$ and every group $\Lambda$ all proper links of $X$ are $(\Lambda, \varepsilon)$ 1-coboundary expanders.
\end{theorem}

Together with the local-to-global machinery (\Cref{DD for cosys thm}) and the results on local-spectral expansion (\Cref{thm: KMS spectral thm}) this theorem gives rise to the following corollary.

\begin{corollary} \label{cor: cosystolic exp of classcial KMS}
	Let $n\geq 3$ be an integer. There exists $\varepsilon >0, q_0 \in \N$ such that for all prime powers $q \geq q_0$, every $n$-dimensional, $n$-classical KMS complex $X$ over $\mathbb{F}_q$ and every group $\Lambda$, we have $h_{\text{cs}}^1(X,\Lambda) \geq \varepsilon$.
\end{corollary}

Combining this with the result from \Cref{cor: homology of KMS complex is trivial}, we can then deduce the main result of this work.

\begin{corollary} \label{cor: main result CBE}
	Let $n \geq 3$ and let $\mathring \Phi$ be root system of type ${A}_{n},{B}_{n}, {C}_{n}$, or ${D}_{n}$ then there exists $q_0 \in \N,\varepsilon>0$ such that 
	\begin{itemize}
		\item  for every finite field $K$ with $|K|  \geq q_0$ and $\operatorname{char}(K)=:p \neq 2$, 
		\item for every group $\Lambda$ with no non-trivial elements of order $p$ (e.g.,  if $\Lambda$ is finite and $\vert \Lambda \vert < p$),
		\item for every irreducible polynomial $f \in K[t]$ of degree at least 2 
	\end{itemize}
the congruence KMS complex $X(\mathring \Phi, K,f)$ is a $(\Lambda,\varepsilon)$ 1-coboundary expander.

	In particular,  given a finite group $\Lambda$,  it holds for any sufficiently large finite field $K$ with $\operatorname{char}(K) > \vert \Lambda \vert$  and any infinite family of irreducible polynomials $f_s \in K[t]$ of degree at least 2 with $\lim_{s \to \infty}\deg(f_s) = \infty$ that the family $\lbrace X(\mathring \Phi, K, f_s) \rbrace_{s \in \N}$ is a family of $1$-coboundary expanders with $\Lambda$ coefficients. 
\end{corollary}

\begin{proof}
	Let $X_f=X(\mathring \Phi, K,f), (\varphi_f(U_{I\setminus \{i\}}))_{i \in I})$.
	First, \Cref{DD for cosys thm}, \Cref{thm: KMS spectral thm} and \Cref{thm: proper links are CBE} imply that $h^1_{\text{cs}}(X_f) > \varepsilon$ for some $\epsilon>0$ only depending on $n, A, q_0$. 
	Then \Cref{cor: homology of KMS complex is trivial} shows that $h^1_{\text{cs}}(X_f)=h^1_{\text{cb}}(X_f)$. 
	The bound on $h^0(X_f)$ follows from the fact that $X_f$ is a $\lambda$-local spectral expander for some $\lambda$ depending only on $A,n,q_0$ (\Cref{thm: KMS spectral thm}) which implies that the 1-skeleton is a $\lambda$-expander graph.
\end{proof}

\begin{remark}
	Note that the results of \Cref{thm: proper links are CBE} hold for all classical KMS complexes, while the result of \Cref{cor: homology of KMS complex is trivial} hold for all congruence KMS complexes. To get \Cref{cor: main result CBE}, both results need to hold, thus we consider the classical congruence KMS complexes.
\end{remark}

\Cref{thm: proper links are CBE} is an adaption of \cite[Theorem 1.1]{CSEofKMS} to the setting with non-Abelian coefficients. The proof structure is very similar,  and we just adapt it to the non-Abelian setting when necessary.

The idea of the proof is as follows:  First,  one observes that proper links of KMS complexes are isomorphic to joins of certain subcomplexes of spherical buildings, called opposite complexes. To further elaborate on this point: 
\begin{itemize}
\item Each proper link has an associate subdiagram of the Dynkin diagram of the underlying GCM.
\item Each connected component of the subdiagram corresponds to an irreducible spherical building. 
\item  Each spherical building admits an opposite relation between its simplices, and after fixing one ``fundamental chamber'' the set of simplices opposite this chamber is called the opposite complex. 
\item A proper link of the KMS complex is then isomorphic to the join of the opposite complexes inside the buildings corresponding to the connected components of the Dynkin diagram. This was shown in \cite[Chapter 3]{hdxfromkms} and we describe more details in \Cref{app: buildings and bn pairs}. 
\end{itemize}

Second,  given two complexes with a 1-cone function, we can construct a 1-cone function for their join as well. Hence it suffices to check that all opposite complexes of irreducible spherical buildings admit a NAC.

Third,  for  all ``classical'' buildings, i.e. those of type $A_n,C_n,D_n$,  we construct a NAC for their opposite complex.  This is done by using an inductive argument on a class of complexes which contains the desired opposite complexes. This argument was first used by Abramenko in \cite{Abr} to show $(n-1)$-sphericallity of the opposite complexes and was adapted to show the existence of Abelian cone functions in \cite{CSEofKMS}.  Below,  we prove the non-Abelian version of the induction argument. The proof of the fact that the relevant classes satisfy the requirements of the theorem can be found in the appendix for completeness, but it almost directly follows from the proof for the Abelian case. 

We note that the opposite complexes of buildings of type $D_n$ require a bit more technical work and therefore are dealt with in the second part of this section.

\subsection{Non-Abelian cone functions for opposite complexes of type $A_n$ and $C_n$}

The existence of non-Abelian cone functions of opposite complexes of buildings of type $A_n$ and $C_n$ will follow from the following theorem, which is a quantitative version of \cite[Lemma 22]{Abr}, see also \cite[Theorem 6.1]{CSEofKMS}. 

\begin{theorem} \label{thm: generalized induction argument}
Let $\mathbf{C}$ be a class of non-empty simplicial complexes. 
Assume that for any $n \in \mathbb{N}, n\geq 2$ there exists $\ell_{n} \in \mathbb{N}$ such that for all $\kappa \in \mathbf{C}, \dim \kappa = n$, there exists a filtration $\kappa_{0} \subset \kappa_{1} \subset \dots \subset \kappa_{\ell} = \kappa$ of $\kappa$ of length $\ell \leq \ell_n$ such that the following is satisfied (set $V_{i}:= \left\{ \text{vertices of } \kappa_{i} \right\} \setminus \left\{ \text{vertices of } \kappa_{i-1} \right\}$):
\begin{enumerate}
\item There is an $1$-cone function $\Cone_{\kappa_{0}}$ such that for $j=0,1$, 
 $$\Rad_{j}(\Cone_{\kappa_{0}}) \leq f(n),$$
	where $f: \mathbb{N} \to \mathbb{R}$ is a function depending only on $\mathbf{C}$. \label{gen L22: condition kappa0}
\item $\kappa_{i}$ is a full subcomplex of $\kappa$ and for all vertices $w_1,w_2 \in V_{i}$ we have $\{w_1,w_2\}\notin \kappa$, for all $0 \leq i \leq \ell$.
\label{gen L22: condition full and dif type}
\item  $\lk_{\kappa_{i}}(w) \cap \kappa_{i-1}$ is either in $\mathbf{C}$ or can be written as the join of two elements from $\mathbf{C}$ for any $w \in V_{i}, 1 \leq i \leq \ell$. \label{gen L22: condition links}
\end{enumerate}
Furthermore, we assume that there exists $\ell_1 \in \N$ such that for any $\kappa \in \mathbf{C}$ with $\dim(\kappa) = 1$, there exists a filtration $\kappa_{0} \subset \kappa_{1} \subset \dots \subset \kappa_{\ell} = \kappa$ of $\kappa$ of length $\ell \leq \ell_1$ such that the following is satisfied (set $V_{i}:= \left\{ \text{vertices of } \kappa_{i} \right\} \setminus \left\{ \text{vertices of } \kappa_{i-1} \right\}$):
\begin{enumerate}
\item There is an $0$-cone function $\Cone_{\kappa_{0}}$ such that 
 $$\Rad_{0}(\Cone_{\kappa_{0}}) \leq f(1),$$
\item $\kappa_{i}$ is a full subcomplex of $\kappa$ and for all vertices $w_1,w_2 \in V_{i}$ we have $\{w_1,w_2\}\notin \kappa$, for all $0 \leq i \leq \ell$,
\item $\lk_{\kappa_{i}}(w) \cap \kappa_{i-1}$ is non-empty for all $w \in V_i, 1 \leq i \leq \ell$.
\end{enumerate}  

Then for every $n \in \mathbb{N}, n\geq 2$ there exists a constant $\mathcal{R}(n)$ such that for every $\kappa \in \mathbf{C}$ with $\dim \kappa = n$ we have:
$\kappa$ has a $1$-cone function $\Cone_{\kappa}$ such that 
$$\Rad(\Cone_{\kappa}) \leq \mathcal{R}(n).$$
Furthermore, every $\kappa \in \mathbf{C}$ with $\dim(\kappa)=1$ admits a 0-cone function with $\Rad_0(\Cone_\kappa) \leq f(1)+\ell_1$.

We can describe $\mathcal{R}(n)$ recursively. We have $\mathcal{R}(0) = 1$. For $n\geq 1$, assume $\mathcal{R}(k)$ is known for $k <n$. Set 
$$S(n) = \max \lbrace 2, \mathcal{R}(n-1) \rbrace. $$
Then 
$$\mathcal{R}(n) = S(n)^{\ell_n} f(n) + \sum_{j=1}^{\ell_n} S(n)^j.$$
\end{theorem}

\begin{proof}
We will prove the result by induction on the dimension $n$.

For $n=1$, we can apply \Cref{adding vertices to cone thm - dim 1}  inductively to each step of the filtration to get a 0-come function $\Cone_\kappa$ with $\Rad_0(\Cone_\kappa) \leq f(1) +\ell$.

Fix $n\geq 2$ and $\kappa \in \mathbf{C}, \dim(\kappa)=n$ with filtration $\kappa_{0} \subset \dots \subset \kappa_{\ell}=\kappa$ as in the requirements.
We will prove by induction on $0 \leq i \leq \ell$, that there is an $1$-cone function $\text{Cone}_{\kappa_{i}}$ such that 
$$\text{Rad}(\text{Cone}_{\kappa_{i}}) \leq \mathcal{R}^{(i)}(n) $$
where $\mathcal{R}^{(i)}(n)$ is given inductively by
$$\mathcal{R}^{(0)}(n) = f(n), \qquad \mathcal{R}^{(i)}(n)= S(n)(\mathcal{R}^{(i-1)}(n)+1)=S(n)^if(n) + \sum_{j=1}^i S(n)^j.$$

For $i=0$ this follows from Assumption \ref{gen L22: condition kappa0}.
We proceed by induction on $i$. Fix $1 \leq i \leq \ell$ and assume there exists a constant $\mathcal{R}^{(i-1)}(n)$ and a $1$-cone function $\text{Cone}_{\kappa_{i-1}}$ such that $\text{Rad}(\text{Cone}_{\kappa_{i-1}}) \leq \mathcal{R}^{(i-1)}(n)$.
We want to apply \Cref{adding vertices to cone thm} to the following setting:

\begin{align*}
X=\kappa, X' = \kappa_{i-1},\  \operatorname{Cone}_{X'}= \operatorname{Cone}_{\kappa_{i-1}}, \ R' = \mathcal{R}^{(i-1)}(n), \ W = V_{i} = \kappa_{i}(0) \setminus \kappa_{i-1}(0).
\end{align*}
We check that the conditions of the theorem are satisfied.
\begin{enumerate}
\item $W\cap X' = \emptyset$ by definition.
\item For every $w_{1},w_{2} \in W: \left\{ w_{1},w_{2} \right\} \not\in \kappa_{i} \subset \kappa$ by assumption.
\item For every $w \in W, \lk_{X}(w) \cap X'$ is a non-empty simplicial complex, since by Assumption \ref{gen L22: condition links} we have $\lk_{X}(w) \cap X' \in \mathbf{C}$ or a join of two complexes in $\mathbf{C}$.  
\item By Assumption \ref{gen L22: condition links}, $\lk_{X}(w) \cap X'$ is either in $\textbf{C}$ or is a join of two elements of $\textbf{C}$. In the first case we use the induction hypothesis and get an $0$-cone function $\Cone_{\lk_{X}(w) \cap X'}$ such that 
$$\text{Rad}_{0}(\operatorname{Cone}_{\lk_{X}(w)\cap X'}) \leq \mathcal{R}(n-1).$$
(Note that any 1-cone function also gives rise to a 0-cone function trivially.)

In the second case, we can write $\lk_{X}(w) \cap X' = A *B$ for some $A,B \in \textbf{C}$.
Thus by \Cref{lemma: join of two scs without assumption}, there exists a 0-cone function $\Cone_{A*B}$ with $\Rad_0(\Cone_{A*B}) \leq 2$.

Hence we can set $R'' = \max\{2,\mathcal{R}(n-1)\} = S(n)$.
\end{enumerate}

Applying \Cref{adding vertices to cone thm} we get a $1$-cone function $\operatorname{Cone}_{X' \cup W}$ such that 
$$\text{Rad}(\operatorname{Cone}_{X' \cup W}) \leq R''(R'+1)$$
which in the original formulation means
$$\text{Rad}(\operatorname{Cone}_{\kappa_{i}}) \leq S(n)(\mathcal{R}^{(i-1)}(n)+1).$$
Thus we have $\mathcal{R}^{(i)}(n) = S(n)(\mathcal{R}^{(i-1)}(n)+1)$ and solving the recursion gives the second expression of $\mathcal{R}^{(i)}(n)$.

Now, since $\kappa_{\ell}=\kappa$, we get a $1$-cone function $\Cone_\kappa$ with $\Rad(\Cone_\kappa) \leq \mathcal{R}^{(\ell)}(n)$. We would like to set $\mathcal{R}(n)= \mathcal{R}^{(\ell)}(n)$. But the bound has to hold for all possible $\kappa \in \mathbf{C}$ of dimension $n$. Since $\mathcal{R}^{(i)}(n)$ is increasing in $i$, we take $\mathcal{R}(n) = \mathcal{R}^{(\ell_n)}(n) = S(n)^{\ell_{n}} f(n) + \sum_{j=1}^{\ell_{n}} S(n)^{j}$.

\end{proof}

\begin{remark} \label{rmk: idea going from non ab to ab in An and Cn}
	Assume that a class of simplicial complexes $\mathbf{C}$ satisfies the assumptions of the general induction argument in the Abelian case \cite[Theorem 6.1]{CSEofKMS}. To show that $\mathbf{C}$ satisfies the conditions of \Cref{thm: generalized induction argument}, it suffices to show that the $\kappa_0$ of the filtration of each complex from the Abelian version of the theorem also admits a non-Abelian 1-cone function (or 0-cone function in the case where it is one-dimensional).
The fact that this holds in the case of the classes containing the $A_n$ and $C_n$ opposite complexes follows, since the cone functions for the $\kappa_0$ are constructed using joins and the ``adding vertices'' construction, which both also work in the non-Abelian case. 
\end{remark}

We get the following corollaries, where we denote for a building $\Delta$ and $a \in \Delta$ by $\Delta^\circ (a)$ the complex opposite of $a$. More details on the opposite relation are explained in \Cref{app: buildings and bn pairs}. The idea of the proofs of the corollaries was explained in \Cref{rmk: idea going from non ab to ab in An and Cn}, for completeness, we provide the details in \Cref{app: An subcomplexes sec} and \Cref{app: Cn subcomplexes section}, .

\begin{corollary}[\Cref{app: cor An case}] 
For every $n \in \mathbb{N},  n \geq 2$ there is a constant $\mathcal{R} (n-1)$ such that for every building $\Delta$ of type $A_{n}, n\geq 2$ (which is a building of dimension $n-1$) over a field $K$ with $\lvert K \rvert \geq 2^{n-1}$ and $a\in \Delta$ a simplex, there exists a NAC $\Cone_{\Delta^0(a)}$ of $\Delta^0(a)$ such that  
$$\Rad(\Cone_{\Delta^0(a)}) \leq \mathcal{R} (n-1).$$
\end{corollary}

\begin{corollary}[\Cref{app: cor Cn case}] 
For every $n \in \mathbb{N},  n \geq 2$ there is a constant $\mathcal{R} (n-1)$ such that for every building $\Delta$ of type $C_{n}, n\geq 2$ (which is a building of dimension $n-1$) over a field $K$ with $\lvert K \rvert \geq 2^{2n-1}$ and $a\in \Delta$ a simplex, there exists a NAC $\Cone_{\Delta^0(a)}$ of $\Delta^0(a)$ such that  
$$\Rad(\Cone_{\Delta^0(a)}) \leq \mathcal{R} (n-1).$$
\end{corollary}

\subsection{Non-Abelian cone functions for opposite complexes of type $D_n$}

The difficulties in this case originate from the incidence geometric description of buildings of type $D_n$. In \cite{Abr}, Abramenko works with a weak $C_n$ type building and then describes the connection between it and the thick $D_n$ building (basically, one is a subdivision of the other). Thus, we first need to prove that a 1-cone function for the opposite complexes of a weak $C_n$ building gives rise to a 1-cone function for the corresponding thick $D_n$ case, analogues to \cite[Proposition 6.10]{CSEofKMS}. To do so, we briefly recall the geometric set up, more details can be found in the appendix.

\subsubsection{Buildings of type $D_n$ - geometric set up} \label{subsubsec: geom descr D_n}
For this section, let $(V,Q,f)$ be a pseudo-quadratic space but in the specific case where $K$ is a field, $\sigma=\operatorname{id}_K, \epsilon = 1, \Lambda = 0$, $V$ a $K$-vector space of dimension $m=2n \geq 4$, $Q$ an ordinary quadratic form and $f$ the corresponding symmetric bilinear form. We further assume that $V$ is non-degenerate and of Witt index $n$. 

We set $X = X(V)=\{0<U<V \mid U \text{ is totally isotropic}\}$. Then $\Delta = \flag X$ is a weak $C_n$ building, in particular it is not thick, since each totally isotropic subspace of dimension $n-1$ is contained in exactly two totally isotropic subspaces of dimension $n$.  Furthermore, the space $Y= \{U\in X \mid \dim U = n\}$ is partitioned into two sets $Y = Y_1\sqcup Y_2$ such that for $U_1,U_2 \in Y_i$ we have $n-\dim(U_1 \cap U_2)\in 2\Z$ for $i=1,2$ and for $U_1 \in Y_1, U_2 \in Y_2$ we have $n - \dim(U_1 \cap U_2) \in 2\Z +1$. In particular, if for $U_1,U_2 \in Y$ we have $\dim(U_1 \cap U_2) = n-1$ then $U_1$ and $U_2$ are in different sets of the partition. 

To obtain a thick building of type $D_n$ we set $\tilde{X}=\tilde{X}(V) = \{U \in X \mid \dim U \neq n-1\}$ and we say that $U,W \in \tilde{X	}$ are incident if and only if
$$ U \subseteq W \text{ or } W \subseteq U \text{ or } \dim (W \cap U) = n-1.$$
Then $\tilde{\Delta} = \flag \tilde{X}$ is the desired building. We will write $\orifl(\tilde{X})$ for $\flag \tilde{X}$ to stress that we consider $\tilde{X}$ not with the usual incidence relation given by inclusion, but with this new one. Note that $\orifl \tilde{X}$ is called the oriflamme complex of $\tilde{X}$.  

For $n\geq 4$ every building of type $D_n$ is of the form $\tilde{\Delta}$ for some field $K$ \cite[Proposition 8.4.3]{TitsBuildings}. We also consider the construction for $n=2,3$ in which case we get certain (but not all) buildings of type $D_2 = A_1 \times A_1$ and $D_3 = A_3$, which where already covered earlier.

To describe the opposite complexes, we need to introduce some further notation. 
\begin{definition} 
For arbitrary subspaces $U, W$ of $V$, we define
$$U\tilde{\pitchfork}_V W \iff U \pitchfork_V W \text{ or } (U,W \in \tilde{X}, U = U^\perp, W = W^\perp \text{ and } \dim(U \cap W) = 1).$$
For a set $\EE$ of subsets of $V$ we define 
\begin{align*}
X_{\EE}(V) &= \{U \in X \mid U \tilde{\pitchfork} \EE \}, \qquad T_{\EE}(V) = \flag X_{\EE}(V) \\
\tilde{X}_{\EE}(V) &= \{U \in \tilde{X} \mid U \tilde{\pitchfork} \EE \}, \qquad \tilde{T}_{\EE}(V) = \orifl \tilde{X}_{\EE}(V). 
\end{align*}
\end{definition}

\subsubsection{Cone functions for a subdivision}
With the above defined notation, we can now describe how to construct a 1-cone function for $\tilde{T}_{\EE}(V)$ given a 1-cone function for $T_{\EE}(V)$.

\begin{proposition} \label{prop: quantitative subdivision D_n}
Let $\tau = \{E_1<\dots < E_r\} \in \tilde{\Delta}$ and set $\EE(\tau) = \{E_i, E_i^\perp \mid 1 \leq i\leq r\}$. If $T_{\EE(\tau)}(V)$ has an $1$-cone function $\Cone_{T_{\EE(\tau)}(V)}$ with $\Rad_i(\Cone_{T_{\EE(\tau)}(V)}) \leq c, i=0,1$ for some $c \geq 0$ then there exists a $1$-cone function $\Cone_{\tilde{T}_{\EE(\tau)}(V)}$ with $\Rad_0(\Cone_{\tilde{T}_{\EE(\tau)}(V)}) \leq c$ and $\Rad_1(\Cone_{\tilde{T}_{\EE(\tau)}(V)}) \leq 2c$. 
\end{proposition}

To prove the proposition, we need the following definition and auxiliary lemmas. 

\begin{definition}\label{def: map varphi for subdivision}
	Let $\tau = \{E_1<\dots < E_r\} \in \tilde{\Delta}$ and set $\EE(\tau) = \{E_i, E_i^\perp \mid 1 \leq i\leq r\}$. Denote $T = T_{\EE(\tau)}(V)$ and $\tilde{T} =  \tilde{T}_{\EE(\tau)}(V)$. 
	Each $u \in X$ with $\dim u = n-1$ is contained in two $n$-dimensional spaces in $X$. We denote them by $W_u$ and $\bar{W}_u$.

	We define the following function: 
	\begin{align*}
		\varphi: T_{\ord}(1) &\to \text{ paths in } \tilde{T} \\
					(u,v) & \mapsto \begin{cases}
						(W_u;v) & \text{ if } \dim u = n-1, v \neq W_u \\
						(W_u) & \text{ if } \dim u = n-1, v = W_u \\
						(u;W_v) & \text{ if } \dim v = n-1, u \neq W_v \\
						(W_v) & \text{ if } \dim v = n-1, u = W_v \\
						(u;v) & \text{ else.}
					\end{cases}
	\end{align*}
\end{definition}
Note that this definition is well defined, since if $(u,v)$ is an edge in $T$ and $\dim v = n-1$ then $(u, W_v)$ is an edge in $\tilde{T}$ (either $u \subseteq v \subseteq W_v$ or $u \cap W_v = v$).

\begin{lemma} \label{lemma: aux lemma path}
	The map $\varphi$ preserves paths, i.e. given a path $P = (v_0; \dots;v_m)$ in $T$, we have that $\varphi(P) := \varphi (v_0,v_1) \circ \varphi (v_1,v_2) \circ \dots \circ \varphi (v_{m-1},v_m)$ is a path in $\tilde{T}$. 
\end{lemma}
\begin{proof}
Given two edges $(u,v), (v,w) \in T_{\ord}(1)$ such that the terminal vertex of one edge equals the initial vertex of the other. Then the terminal vertex of $\varphi((u,v))$ will be equal to the initial vertex of $\varphi((v,w))$ (if $\dim v \neq n-1$ it will be $v$, otherwise it will be $W_v$).
Hence, the concatenation of the images of the edges of a path is well defined and gives again a path.
\end{proof}

Note that it follows from the definition that $\varphi(P \circ Q) = \varphi(P) \circ \varphi(Q)$ and $\varphi(P^{-1})=\varphi(P)^{-1}$.
\begin{lemma} \label{lemma: aux lemma BT}
	Given two paths $P, P'$ in $T$, if $P \sim_{(BT)} P'$ then $\varphi(P) \sim_{(BT)} \varphi(P')$ or $\varphi(P) = \varphi(P')$.
\end{lemma}

\begin{proof}
	Since $P \sim_{(BT)} P'$ there exit paths $Q_0,Q_1$ and vertices $u,v$ in $T$ such that $P = Q_0 \circ (u;v;u) \circ Q_1^{-1}$ and $P' = Q_0 \circ (u) \circ Q_1^{-1}$.
	Thus, $\varphi(P)= \varphi (Q_0) \circ \varphi((u;v;u)) \circ \varphi(Q_1)^{-1}$ and similarly for $\varphi(P')$. 
	Now, if $\dim u, \dim v \neq n-1$ then $\varphi((u;v;u))=(u;v;u)$ and $\varphi((u)) = (u)$ which proves the statement in this case.
	
	If $\dim u =n-1$ then $\varphi((u))= (W_u)$ and $\varphi((u;v;u))= (W_u;v;W_u)$ or $\varphi((u;v;u)) = (W_u)$ and the statement follows.

	If $\dim v = n-1$ then $\varphi((u))=(u)$ and $\varphi((u;v;u))= (u;W_v;u)$ or $\varphi((u;v;u))=(W_v)=(u)$ (the second case occurs exactly when $u = W_v$).
\end{proof}

\begin{lemma} \label{lemma: aux lemma TR}
	Given two paths $P, P'$ in $T$, if $P \sim_{(TR)} P'$ then $\varphi(P) \sim_{(TR)} \varphi(P')$, $\varphi(P) \sim_{(BT)} \varphi(P')$ or $\varphi(P) = \varphi(P')$.
\end{lemma}

\begin{proof}
Since $P \sim_{(TR)} P'$ there exit paths $Q_0,Q_1$ and a triangle $(u,v,w)$ in $T$ such that $P = Q_0 \circ (u;v;w) \circ Q_1^{-1}$ and $P' = Q_0 \circ (u;w) \circ Q_1^{-1}$.
	Thus, $\varphi(P)= \varphi (Q_0) \circ \varphi((u;v;w)) \circ \varphi(Q_1)^{-1}$ and similarly for $\varphi(P')$.

	We will proceed by case distinction which we display in the following table. The last column shows the relationship between $\varphi(P)$ and $\varphi(P')$.

\begin{center}
\begin{tabular}{ |c|c|c||c|c||c| } 
 \hline
 $ u$ & $v$ & $ w$ & $\varphi((u;v;w))$& $\varphi((u;w))$ &  $\varphi(P)$ vs $\varphi(P')$ \\
 \hline  \hline
 $ \dim u \neq n-1$ & $\dim v \neq n-1$ & $\dim w \neq n-1$ & $(u;v;w)$ & $(u;w)$ & $\sim_{\text{TR}}$\\
 \hline \hline
 $u \neq W_w$ & $v \neq W_w$ & $\dim w =n-1$ & $(u;v;W_w) $ & $(u;W_w)$ & $\sim_{\text{TR}}$ \\ 
\hline
 $u \neq W_w$ & $v = W_w$ & $\dim w =n-1$ & $(u;W_w) $ & $(u;W_w)$ & $=$ \\ 
\hline
 $u = W_w$ & $v \neq W_w$ & $\dim w =n-1$ & $(u;v;W_w) $ & $(W_w)$ & $\sim_{\text{BT}}$ \\ 
 \hline \hline
$u \neq W_v$ & $\dim v =n-1 $ & $w \neq W_v$ & $(u;W_v;w) $ & $(u;w)$ & $\sim_{\text{TR}}$ \\ 
 \hline
 $u \neq W_v$ & $\dim v =n-1 $ & $w = W_v$ & $(u;W_v) $ & $(u;w=W_v)$ & $=$ \\ 
 \hline
 $u = W_v$ & $\dim v =n-1 $ & $w \neq W_v$ & $(W_v;w) $ & $(u=W_v;w)$ & $=$ \\ 
 \hline \hline
 $\dim u =n-1$ & $v \neq W_u $ & $w \neq W_u$ & $(W_u;v;w) $ & $(W_u;w)$ & $\sim_{\text{TR}}$ \\ 
\hline
$\dim u =n-1$ & $v = W_u $ & $w \neq W_u$ & $(W_u;w) $ & $(W_u;w)$ & $=$ \\ 
\hline
$\dim u =n-1$ & $v \neq W_u $ & $w = W_u$ & $(W_u;v;w) $ & $(W_u)$ & $\sim_{\text{BT}}$ \\ 
\hline

\end{tabular}
\end{center}	
\end{proof}

\textit{Proof of \Cref{prop: quantitative subdivision D_n}.}
We consider the setting of \Cref{def: map varphi for subdivision}. Assume that we have a 1-cone function $(v_0,(P_v)_{v \in T(0)}, (T_{(u,v)})_{(u,v) \in T_{\ord}}(1))$ of $T$. The goal is to construct a 1-cone function $(\tilde{v_0},(\tilde{P}_v)_{v \in \tilde{T}(0)}, (\tilde{T}_{(u,v)})_{(u,v) \in \tilde{T}_{\ord}(1)})$ for $\tilde{T}$.

For the apex, we set $\tilde{v_0} = v_0$ if $\dim v_0 \neq n-1$ and otherwise $\tilde{v_0}= W_{v_0}$.

Given a vertex $v \in \tilde{T}(0) = \tilde{Y}$, we set $\tilde{P}_v = \varphi(P_v)$. Since $v \in \tilde{Y}$, we have $\dim v \neq n-1$ and hence $\tilde{P}_v$ is a path from $\tilde{v_0}$ to $v$ of length at most the length of $P_v$. 

Let $(u,v) \in \tilde{T}_{\ord}(1)$. To define the contraction $\tilde{T}_{(u,v)}$, we distinguish two cases. First assume that $(u,v) \in T \cap \tilde{T}$, i.e. at least one of the two vertices has dimension less than $n$. Hence we can consider $T_{(u,v)} = (P_0, \dots, P_m)$ and define $\tilde{T}_{(u,v)} = (\varphi(P_0),\dots, \varphi(P_m))$. The auxiliary Lemmas \ref{lemma: aux lemma path}, \ref{lemma: aux lemma BT}, \ref{lemma: aux lemma TR} imply that this is a well defined contraction and $|\tilde{T}_{(u,v)}| \leq |T_{(u,v)}|$.

For the second case, let $(u,v) \in \tilde{T}_{\ord}(1)$ with $\dim u = \dim v =n$. Consider $T_{(u,u \cap v)} = (P_0^1, \dots, P_m^1)$ and $T_{(u \cap v,v)} = (P_0^2, \dots, P_\ell^2)$. Note that $\varphi((u;u\cap v; v))=(u,v)$. We set
\begin{align*}
	\tilde{P}_0 = \tilde{P}_u \circ \varphi((u,u\cap v) \circ (u \cap v,v)) \circ \tilde{P}_v^{-1} = \varphi(P_u \circ (u;u\cap v; v) \circ P_v^{-1}).
\end{align*}
By $|P_{u \cap v}|=:s$ many backtracking steps, we can transform $\tilde{P}_0$ to
\begin{align*}
	\tilde{P}_s = \varphi(P_u \circ (u,u\cap v) \circ P_{u \cap v}^{-1} \circ P_{u\cap v} \circ (u \cap v, v) \circ P_v^{-1}) = \varphi(P_0^1 \circ P_0^2).
\end{align*}
Applying the contraction steps from $T_{(u,u\cap v)}$ on the first half of the path and then the steps from $T_{(u \cap v,v)}$ on the second half we get the desired contraction, (again using \Cref{lemma: aux lemma BT} and \Cref{lemma: aux lemma TR}). In particular $|\tilde{T}_{(u,v)}| \leq |T_{(u, u \cap v)}| + |T_{(u \cap v,v)}|$. This also proves the bound on the 1-cone radius.

\qed

\begin{remark}
	The above reasoning also shows that if $T$ has a 0-cone function, then we can construct a 0-cone function for $\tilde{T}$. Further note that the dimensions of $T$ and $\tilde{T}$ coincide. Hence, we can also get a more general statement of \Cref{prop: quantitative subdivision D_n} saying that given a NAC for $T$, we can construct a NAC for $\tilde{T}$ and $\Rad(\tilde{T}) \leq 2\Rad(T)$.
\end{remark}

\subsubsection{Induction argument for the type $D_n$ case}

For the $D_n$ case, we need the following weaker and more specialized version of \Cref{thm: generalized induction argument}. In order to state it, we need the following definition. See \cite[Theorem 6.13]{CSEofKMS} for the Abelian version.
\begin{definition}
Let $K$ be a field, $V$ a $K$-vector space, $Y$ a collection of subspaces of $V$. Let $X = \flag Y$ and $U\in Y$. Then we define
\begin{align*}
 Y^{<U}&:=\{W\in Y \mid W < U \}, \ X^{<U} := \flag Y^{<U}\\
 Y^{>U}&:=\{W\in Y \mid W > U \}, \ X^{>U} := \flag Y^{>U}. 
\end{align*}
\end{definition}
\begin{remark}
Note that
$$\lk_X(U) = X^{<U} * X^{>U}.$$
\end{remark}

\begin{theorem} \label{thm: the big induction for the D_n case}
Let $\mathbf{C}$ be a class of tuples $(\kappa, (V,Q,f))$ where $(V,Q,f)$ is a thick pseudo-quadratic space and $\kappa = \flag Y$ where $Y\neq \emptyset$ is a set of vector subspaces of $V$. Assume that for any $n \in \N, n \geq 2$ there exist $\ell_n \in \N$ such that for all $\kappa \in \mathbf{C}, \dim \kappa = n$ there exists a filtration $\kappa_0 \subset \kappa_1 \subset \dots \subset \kappa_\ell = \kappa, \ell \leq \ell_n$ satisfying:
\begin{enumerate}
\item There exists an $1$-cone function $\Cone_{\kappa_0}$ such that for $j = 0,1$
$$\Rad_j(\Cone_{\kappa_0}) \leq f(n) ,$$
where $f: \N \to \R$ is a function depending only on $\mathbf{C}$.
\item Every $\kappa_i$ is a full subcomplex of $\kappa$, i.e. $\kappa_i = \flag Y_i$ for some $Y_i \subseteq Y$ and all vector subspaces in $Y_i \setminus Y_{i-1}$ have the same dimension. 
\item For every $U \in Y_i \setminus Y_{i-1}$ we have that $\lk_{\kappa_i}(U) \cap \kappa_{i-1} = \kappa_{i-1}^{<U} * \kappa_{i-1}^{>U}$ is such that $\kappa_{i-1}^{<U}, \kappa_{i-1}^{>U} \in \mathbf{C}$ or one (or both) have again a filtration $R_0 \subset \dots R_h, h \leq \ell_{i-1}$ which starts with an element of $\mathbf{C}$ (i.e. $R_0\in \mathbf{C}$), such that Condition 2. is satisfied and such that for every $W \in R_j(0) \setminus R_{j-1}(0)$ we have $R_{j-1}^{<W},R_{j-1}^{>W} \in \mathbf{C}$. 
\end{enumerate} 
Furthermore, we assume that there exists $\ell_1 \in \N$ such that for any $\kappa \in \mathbf{C}$ with $\dim(\kappa) = 1$, there exists a filtration $\kappa_{0} \subset \kappa_{1} \subset \dots \subset \kappa_{\ell} = \kappa$ of $\kappa$ of length $\ell \leq \ell_1$ such that the following is satisfied (set $V_{i}:= \left\{ \text{vertices of } \kappa_{i} \right\} \setminus \left\{ \text{vertices of } \kappa_{i-1} \right\}$):
\begin{enumerate}
\item There is an $0$-cone function $\Cone_{\kappa_{0}}$ such that 
 $$\Rad_{0}(\Cone_{\kappa_{0}}) \leq f(1),$$
\item $\kappa_{i}$ is a full subcomplex of $\kappa$ and for all vertices $w_1,w_2 \in V_{i}$ we have $\{w_1,w_2\}\notin \kappa$, for all $0 \leq i \leq \ell$.
\item $\lk_{\kappa_{i}}(w) \cap \kappa_{i-1}$ is non-empty for all $w \in V_i, 1 \leq i \leq \ell$.
\end{enumerate}  

Then for every $n \in \mathbb{N}, n\geq 2$ there exists a constant $\mathcal{R}(n)$ such that for every $\kappa \in \mathbf{C}$ with $\dim \kappa = n$ we have:
$\kappa$ has a $1$-cone function $\Cone_{\kappa}$ such that for $j=0,1$ we  have
$$\Rad_{j}(\Cone_{\kappa}) \leq \mathcal{R}(n).$$
Furthermore, every $\kappa \in \mathbf{C}$ with $\dim(\kappa)=1$ admits a 0-cone function with $\Rad_0(\Cone_\kappa) \leq f(1)+\ell_1$.
\end{theorem}
\begin{proof}
The proof is analogues to the one of \Cref{thm: generalized induction argument}.
The main difference is that to get a cone function for $\lk_{\kappa_i}(U) \cap \kappa_{i-1}$ it might happen that we cannot just use the induction hypothesis, but that we need to do another induction on the filtration of $\kappa_{i-1}^{<U}$ or $\kappa_{i-1}^{>U}$, using \Cref{adding vertices to cone thm} in the step of the induction. This gives an even weaker and more complicated bound on the cone radius of the final cone function, thus we decided to not keep track of it. But it will still only depend on $n$ and the class $\mathbf{C}$, not on the field $K$. 
\end{proof}

\begin{remark} \label{rmk: from Abelian to non-Abelian in Dn}
	Similar to \Cref{thm: generalized induction argument}, if we are give a class of simplicial complexes $\mathbf{C}$ which is know to satisfy the assumptions of the Abelian version of this theorem, \cite[Theorem 6.13]{CSEofKMS}, then it suffices to check that the $\kappa_0$ appearing in the filtrations of the Abelian version also admit a non-Abelian 1-cone function (resp. 0-cone function if $\kappa_0$ is one-dimensional).

\end{remark}

We get the following corollary, details about the terminology and the proof can be found in \Cref{app: sec dn case}.

\begin{corollary}[\Cref{app: dn case cor}]
For every $n \in \N, n\geq 2$ there exist a constant $\mathcal{R}(n-1)$ such for every building $\tilde{\Delta}$ of type $D_n$ over a field $K$ with $\lvert K \rvert \geq 2^{2n-1}$ and $a \in \tilde{\Delta}$ there exists a NAC $\Cone_{\tilde{\Delta}^0(a)}$ of ${\tilde{\Delta}^0(a)}$ such that
$$\Rad(\Cone_{\tilde{\Delta}^0(a)}) \leq 2 \mathcal{R}(n-1).$$
\end{corollary}

\subsection{Proof of \Cref{thm: proper links are CBE}}

The above bounds on the cone radii of NAC of opposite complexes yields the proof of \Cref{thm: proper links are CBE}:
\begin{proof}
Let $n\geq 3$ be an integer and $q > 2^{2n-1}$ be a prime power.  For every $n$-dimensional, $n$-classical KMS complex $X$ over $\mathbb{F}_q$ every link of a vertex of $X$ is either an opposite complex of a classical type,  i.e.,  of type $A_n,  B_n,  C_n$ or $D_n$ or a join of two such opposite complexes.  By Corollaries \ref{app: cor An case}, \ref{app: cor Cn case},  \ref{app: dn case cor} it follows that the link is a strongly symmetry simplicial complex that has a non-Abelian cone function whose radii are bounded independently of $q$ (in the case of the join, each factor has a non-Abelian cone function whose radii are bounded independently of $q$ and thus by \Cref{cone of a join - general prop},  the join also has such cone function).  Thus by \Cref{lemma: cone implies CBE} it follows that every link is a $1$-coboundary expander with respect to any group $\Lambda$ and that there is a bound on the expansion that is independent of $q$.
\end{proof}

\appendix
\section{Buildings, their opposite complexes and relation to KMS complexes} \label{app: sec buildings}
This section is almost identical with \cite[Chapter 5]{CSEofKMS}. We include it to improve the self-containedness of this paper and to recall the necessary notations.

\subsection{Buildings and BN-pairs} \label{app: buildings and bn pairs}
One way to define a building is as a simplicial complex with a family of subcomplexes called apartments. Each apartment is isomorphic to the same Coxeter complex. The type of the Coxeter complex, which can be described by a root system/ a GCM, is also what we call the type of the building. Each two simplices of the building have to be contained in a common apartment and each two apartments are isomorphic with an isomorphism that acts as the identity on the intersection.  

One way to construct buildings is via groups with a BN-pair (and for thick, irreducible buildings of dimension 2 and larger, all buildings can be constructed that way, as was proven by Jacques Tits in \cite{TitsBuildings}).

Given a group $G$, a BN-pair of $G$ is a pair of subgroups $B,N \leq G$ satisfying certain properties, see \cite[Definition 6.55]{AB2008}. In particular, $G = BN$ and setting $T=B \cap N$ we have that $W= N/T$ admits a set of generators $S$, such that $(W,S)$ is a Coxeter system. $W$ is also called the Weyl group of the BN-pair. 
To describe the building, we need to define the parabolic subgroups $P_J:= B \langle J \rangle B$ for $J \subset S$ (here $\langle J \rangle$ is the subgroup generated by the elements of $J$ inside $W$). The building can then we described as the following coset complex:
$$\Delta(G,B) := \CC(G; (P_{S\setminus \{j \}})_{j \in S}).$$

Since $P_J \cap P_L = P_{J \cap L}$ for any $J,L \subseteq S$ we can see $\Delta(G,B)$ also as the poset $\{gP_J \mid g \in G, J \subseteq S\}$ ordered by reversed inclusion. Then the simplex $\{g P_{S\setminus \{j\}} \mid j \in J\} \in \CC(G; (P_{S\setminus \{j \}})_{j \in S})$ corresponds to $gP_{S\setminus J}$ for $J \subseteq S, g \in G$.
In particular, the set of maximal simplices, called chambers, corresponds to $G/B$.

If the Weyl group is finite, we call the building spherical. This coincides with our definition of calling a GCM/root system spherical. A finite Weyl group has a unique longest element denoted by $w_0$, where longest means that is generated by the largest number of generators, counted with multiplicities. 

Two chambers $gB, hB, g,h\in G$ are called opposite if and only if $g^{-1}h \in Bw_0B$, denoted by $gB \text{ op } hB$.

Given a simplex $\sigma \in \Delta(G,B)$, the complex opposite $\sigma$ is defined as follows: 
$$\Delta^0(\sigma):= \{\tau \in \Delta(G,B) \mid \text{ there exist chambers } c,d \in \Delta(G,B) \text{ such that } \tau \subseteq c, \sigma \subseteq d, c \text{ op } d\}.$$

To see the connection between opposite complexes and the links of KMS complexes, note that every Chevalley group has a BN-pair in the following way. Let $A=(A_{ij})_{i,j \in I}$ be a spherical GCM of rank at least 2, $K$ a field and $G= \chev_A(K) = \langle x_\alpha(\lambda) \mid \alpha \in \Phi, \lambda \in K \rangle$, using the notation of \Cref{def: chevalley group}.
We set
\begin{alignat*}{3}
U_\alpha &:= \langle x_{\alpha}(\lambda) \mid \lambda \in K \rangle, \alpha \in \Phi,  \qquad &&U^+ &&:= \langle U_ \alpha \mid \alpha \in \Phi^+ \rangle, \qquad T := \langle h_\alpha(\lambda) \mid \alpha \in \Phi, \lambda \in K\rangle,\\
B &:= U^+ T, \qquad &&N &&:= \langle T, n_{-\alpha_i}(-\lambda^{-1}); i \in I, \lambda \in K^* \rangle.
\end{alignat*}
For example \cite[Theorem 7.115]{AB2008} shows that this is indeed a BN-pair. 

Let $C_0= 1B$ denote the fundamental chamber of $\Delta(G,B)$ where $G$ is a Chevalley group as described above. Recall that $w_0$ denotes the longest element on the Weyl group $W$. It satisfies $Bw_0B = U^+ w_0 B$. We can describe the opposite complex in this set-up as 
\begin{align*}
\Delta^{0}(C_{0}) &= \{ gP_{J} \mid g \in G, J \subseteq S: \exists h \in U^+w_{0}B: gP_{J}\supseteq hB  \}\\
&= \{ gP_{J} \mid \exists h \in U^+w_{0}B: h \in gP_{J} \} \\
&= \{ hP_{J} \mid h \in U^+w_{0}B \} \\
&= \{ uw_{0}P_{J} \mid u \in U^+ \}.
\end{align*}

This leads to the following proposition. 
\begin{proposition}{\cite[Chapter 3.2]{hdxfromkms}} \label{prop: opposite complex and U+}
Let $G$ be a Chevalley group with BN-pair as described above. Then 
$$\CC(U^+, (U_{I\setminus \{i\}})_{i \in I}) \cong \Delta^0(C_0).$$
\end{proposition}

If we assume that $A$ is 2-spherical and $\lvert K \rvert \geq 4$ we have $U^+ = \langle U_{\alpha_i} \mid i \in I \rangle$ (see e.g. the comment before \cite[Lemma 3.7]{hdxfromkms}). Combining this with the above proposition and Fact \ref{fact: links in CC} gives the following corollary. 
\begin{corollary}
Let $X$ be a KMS complex, and $\sigma \in X$ a face of dimension less than or equal to $ \dim X-2$ of type $\emptyset \neq T \subset I$. Then 
$$\lk_X(\sigma) \cong \Delta^0(C_0)$$
where $\Delta$ is the building obtained from the BN-pair structure of the Chevalley group of type $A_{I \setminus T}$.
\end{corollary}

\subsection{Geometric constructions of buildings}
In the following section, we describe how buildings of type $A_n$, $C_n$ and $D_n$ can be constructed as flag complexes of certain sets equipped with an incidence relation. We furthermore describe the opposite complexes in this set-up. The general framework is as follows.

\begin{definition}
Let $X$ be a set and let $I \subseteq X \times X$ be a reflexive and symmetric relation. For $a,b \in X$ we write $aIb$ if and only if $(a,b)\in I$ and call $a$ and $b$ incident. The structure gives rise to the following simplicial complex, called the flag complex of $X$: 
$$\flag(X) = \{\sigma \subseteq X \mid \forall a,b \in \sigma: aIb \}.$$
\end{definition} 

\subsubsection{Buildings of type $A_n$} \label{subsubsec: geom descr A_n}
Let $n\in \N, n\geq 1$. Any building of type $A_n$ can be described in the following way, see e.g. \cite{TitsBuildings}. 
Let $K$ be a field or skew field, $V$ an $(n+2)$-dimensional vector space over $K$. Set $X:= X(V) := \{0 <U<V \mid U \text{ is a } K \text{-subspace of } V \}$. Two subspaces $U,W \in X$ are called incident if $U \subseteq W$ or $W \subseteq U$. Then $\Delta = \flag X$ is a building of type $A_n$. 

To describe the opposite complex in this context, we need the following notation, see \cite[Definition 9]{Abr}.
\begin{definition} \label{def: A_n case X_E}
\begin{enumerate}
\item Two subspaces $U,W$ of $V$ are called transversal (in $V$) if $U\cap W = 0$ or $U+W = V$. If this is the case, we write $U\pitchfork_V W$ or simply $U \pitchfork W$.
\item Let $\EE$ be a set of subspaces of $V$. Then we write $U \pitchfork \EE$ if $U\pitchfork E$ for all $E \in \EE$. We set
$$X_{\EE}(V) := \{U \in X \mid U \pitchfork \EE \} \qquad T_{\EE}(V) = \flag X_{\EE}(V).$$
\end{enumerate}
\end{definition}  
We get the following result. 
\begin{proposition}{\cite[Corollary 12]{Abr}}
For any simplex $\sigma = \{E_1 < \dots < E_r\} \in \Delta$ set $\EE(\sigma) = \{E_i \mid 1 \leq i \leq r \}$. Then 
$$\Delta^0(\sigma) = T_{\EE(\sigma)}(V).$$
\end{proposition}

\begin{definition} \label{def: class C_A}
We denote by $\mathbf{C}_A$ the class of all simplicial complexes $T_{\EE}(V)$ where $\dim(T_{\EE}(V)) = n\geq 0$, $K$ is a (skew) field, $V$ an $(n+2)$-dimensional vector space over $K$, $\mathcal{E}$ is a finite set of subspaces of $V$ such that if we set $e_{j}=\lvert \left\{ E \in \mathcal{E} \mid \dim E=j \right\} \rvert$ we have $\sum_{j=1}^{n+1} {{n}\choose{j-1}} e_{j} \leq \lvert K \rvert$. Furthermore, $T_{\EE}(V)$ is defined as described in \Cref{def: A_n case X_E}. 
\end{definition}

\subsubsection{Hermitian and pseudo-quadratic forms}
In this part, we recall facts about hermitian and pseudo-quadratic forms that are needed to describe the buildings of type $C_n$ and $D_n$. We follow \cite[Chapter 5]{Abr}, who in turn mostly follows \cite{TitsBuildings}.

Let $K$ be a skew field, $\sigma: K \to K, a \mapsto a^\sigma$ be an involution, i.e. an anti-automorphism of $K$ such that $\sigma^2 = \operatorname{id}_K$. Let $\epsilon \in \{1,-1\} \subset K$. If $\sigma \neq \operatorname{id}_K$ then we require $\epsilon = -1$. Let $V$ be a right $K$-vector space of dimension $m \in \N \cup \{\infty\}$.
Let $K_{\sigma, \epsilon}:=\left\{\alpha-\alpha^\sigma \varepsilon \mid \alpha \in K\right\}$ and $ K^{\sigma, \epsilon}:=\left\{\alpha \in K \mid \alpha+\alpha^\sigma \varepsilon=0\right\}$. Let $\Lambda$ be a form parameter relative to $(\sigma, \varepsilon)$, i.e. $\Lambda$ is a subgroup of $(K,+)$ satisfying $K_{\sigma, \varepsilon} \subseteq \Lambda \subseteq K^{\sigma, \varepsilon}$ and $\alpha^\sigma \Lambda \alpha \subseteq \Lambda$ for all $\alpha \in K$.

Let $f: V \times V \to K$ be a $(\sigma, \epsilon)$-hermitian form, i.e. $f$ is biadditive and for all $x,y\in V, a,b \in K$ we have $f(xa,yb) = a^\sigma f(x,y) b$ and $f(y,x) = f(x,y)^\sigma \epsilon$.
Let $Q: V \to K/\Lambda$ be a $(\sigma,\epsilon)$-quadratic form with associated $(\sigma, \epsilon)$-hermitian form $f$, i.e. $Q(xa) = a^\sigma Q(x) a +\Lambda$ and $Q(x+y) -Q(x) -Q(y) = f(x,y) + \Lambda$ for all $x,y \in V, a \in K$. 
If $\Lambda = K$ we require that $f$ is alternating. If $\Lambda \neq K$ then $f$ is uniquely defined by $Q$.

For $M\subseteq V$ we write $M^\perp := \{x \in V \mid f(x,M) = 0\}$.
A subspace $U<V$ is called 
\begin{itemize}
\item non-degenerate if $U \cap U^\perp = 0$,
\item totally degenerate if $U \subseteq U^\perp$,
\item anisotropic if $0 \notin Q(U \setminus \{0\})$, 
\item isotropic if $0 \in Q(U \setminus \{0\})$,
\item totally isotropic if $U \subseteq U^\perp$ and $Q(U) = 0$.
\end{itemize}   
We require that there is at least one finite-dimensional maximal totally isotropic subspace. In that case, all maximal totally isotropic subspace have the same dimension. We denote this dimension by $n$ and call it the Witt index of $(V,Q,f)$. We require $0 < n< \infty$.
Furthermore, we require that $V^\perp = 0$. 
We say that the triple $(V,Q,f)$ is a pseudo-quadratic space if it is of the form described above. 

If it additionally satisfies $(m, \Lambda) \neq (2n, 0)$ then we say $(V,Q,f)$ is a thick pseudo-quadratic space.

Note that if $(V,Q,f)$ is a thick pseudo-quadratic space and $K$ is a finite field, then $\dim V \leq 2n+1$.
In general, if one only considers finite fields, especially of characteristic different from 2, the set-up can be substantially simplified, see \cite[Remark 9.3,9.4]{AB2008}. 

\subsubsection{Buildings of type $C_n$} \label{subsubsec: geom descr C_n}
Let $n \in \N, n\geq 1$ and $(V,Q,f)$ be a thick pseudo-quadratic space with Witt index $n$. Set $X := X(V) = \{0 < U < V \mid U \text{ is totally isotropic}\}$ with incidence relation given by containment as in the $A_n$ case and $\Delta = \flag X$. Then $\Delta$ is a thick building of type $C_n$ and every classical $C_n$ building can be obtained in this way (here classical means that the links of type $A_2$ correspond to Desarguesian planes, which is always the case for $n\geq 4$), see \cite[Theorem 8.22]{TitsBuildings}. 

\begin{definition}\label{def: C_n case N(E)}
Set $\mathcal{U} := \{0\leq U \leq V \mid \dim U < \infty \}, \mathcal{U}^\perp := \{U^\perp \mid u \in \mathcal{U}\}$ and $\mathcal{W}:= \mathcal{U} \cup \mathcal{U}^\perp$.  

Let $\EE$ be a finite subset of $\mathcal{W}$ such that $\EE^\perp = \EE$. Assume $K = \mathbb{F}_q$ is a finite field. 
Set $\mathcal{E}_j:=\{E \in \mathcal{E} \mid \operatorname{dim} E=j\}, e_j:= \lvert \mathcal{E}_j \rvert$ and $e_h^{(s)}:=\sum_{j=0}^{2 s}\binom{2 s}{j} e_{h+j}$ for $h \in \mathbb{N}, s \in \mathbb{N}_0$ and $h+2 s<m$. We define 
\begin{itemize}
\item $N(\mathcal{E}):=e_1^{(n-1)}$ if $f$ is alternating $(\Rightarrow m=2 n)$ and $Q=0$, 
\item $N(\mathcal{E}):=\left(e_1^{(n-1)}\right)^2$ if $m=2 n$ and $\sigma \neq \mathrm{id}$, 
\item $N(\mathcal{E}):=2 e_2^{(n-1)}$ if $m=2 n+1$,
\item $N(\mathcal{E}):=\max \left\{e_2^{(n-1)}+e_3^{(n-1)}+1,2 e_3^{(n-1)}\right\}$ if $m=2n+2$.
\end{itemize}
\end{definition}

\begin{definition}\label{def: X_E in C_n case}
For a set $\EE$ of subspaces of $V$ we define, similar to \Cref{def: A_n case X_E}, 
$$X_{\EE}(V) := \{U\in X \mid U \pitchfork \EE\}, \qquad T_{\EE}(V) = \flag X_{\EE}(V).$$ 
\end{definition}

\begin{proposition}{\cite[Corollary 15]{Abr}}
For any simplex $\tau = \{E_1 < \dots < E_r\} \in \Delta$ set $\EE(\tau) = \{E_i, E_i^\perp \mid 1 \leq i\leq r\}$. Then $\Delta^0(\tau) = T_{\EE(\tau)}(V)$.

\end{proposition}

\begin{definition}\label{def: class C_C}
With $\mathbf{C}_C$ we denote the class of simplicial complexes $T_{\EE}(V)$ where $n=\dim(T_{\EE}(V) +1 \geq 1$, $K$ a (skew) field, $(V,Q,f)$ is a thick pseudo-quadratic space with Witt index $n$, $\EE$ a finite subset of $\mathcal{W}$ (see \Cref{def: C_n case N(E)}) such that $\EE^\perp = \EE$. If $K$ is a finite field, we further require that $\lvert K \rvert \geq N(\EE)$. We set $T_{\EE}(V)$ as in \Cref{def: X_E in C_n case}.
\end{definition}

\subsubsection{Buildings of type $D_n$} \label{subsubsec: geom descr D_n}
For this section let $(V,Q,f)$ be a pseudo-quadratic space but in the specific case where $K$ is a field, $\sigma=\operatorname{id}_K, \epsilon = 1, \Lambda = 0$, $V$ a $K$-vector space of dimension $m=2n \geq 4$, $Q$ an ordinary quadratic form and $f$ the corresponding symmetric bilinear form. We further assume that $V$ is non-degenerate and of Witt index $n$. 

We set as before $X = X(V)=\{0<U<V \mid U \text{ is totally isotropic}\}$. Then $\Delta = \flag X$ is a weak $C_n$ building, in particular it is not thick, since each totally isotropic subspace of dimension $n-1$ is contained in exactly two totally isotropic subspaces of dimension $n$.  Furthermore, the space $Y= \{U\in X \mid \dim U = n\}$ is partitioned into to sets $Y = Y_1\sqcup Y_2$ such that for $U_1,U_2 \in Y_i$ we have $n-\dim(U_1 \cap U_2)\in 2\Z$ for $i=1,2$ and for $U_1 \in Y_1, U_2 \in Y_2$ we have $n - \dim(U_1 \cap U_2) \in 2\Z +1$. In particular, if for $U_1,U_2 \in Y$ we have $\dim(U_1 \cap U_2) = n-1$ then $U_1$ and $U_2$ are in different sets of the partition. 

To obtain a thick building of type $D_n$ we set $\tilde{X}=\tilde{X}(V) = \{U \in X \mid \dim U \neq n-1\}$ and we say that $U,W \in \tilde{X	}$ are incident (write $UIW$) if and only if
$$ U \subseteq W \text{ or } W \subseteq U \text{ or } \dim (W \cap U) = n-1.$$
Then $\tilde{\Delta} = \flag \tilde{X}$ is desired building. We will write $\orifl(\tilde{X})$ for $\flag \tilde{X}$ to stress that we consider $\tilde{X}$ not with the usual incidence relation given by inclusion, but with this new one. Note that $\orifl \tilde{X}$ is called the oriflamme complex of $\tilde{X}$.  

For $n\geq 4$ every building of type $D_n$ is of the form $\tilde{\Delta}$ for some field $K$ \cite[Proposition 8.4.3]{TitsBuildings}. We also consider the construction for $n=2,3$ in which case we get certain (but not all) buildings of type $D_2 = A_1 \times A_1$ and $D_3 = A_3$, which where already covered earlier.

To describe the opposite complexes, we need to introduce some further notation. 
\begin{definition} 
For arbitrary subspaces $U, W$ of $V$, we define
$$U\tilde{\pitchfork}_V W \iff U \pitchfork_V W \text{ or } (U,W \in \tilde{X}, U = U^\perp, W = W^\perp \text{ and } \dim(U \cap W) = 1).$$
For a set $\EE$ of subsets of $V$ we define (different from \Cref{def: A_n case X_E})
\begin{align*}
X_{\EE}(V) &= \{U \in X \mid U \tilde{\pitchfork} \EE \}, \qquad T_{\EE}(V) = \flag X_{\EE}(V) \\
\tilde{X}_{\EE}(V) &= \{U \in \tilde{X} \mid U \tilde{\pitchfork} \EE \}, \qquad \tilde{T}_{\EE}(V) = \orifl \tilde{X}_{\EE}(V) 
\end{align*}
\end{definition}

\begin{proposition}{\cite[Corollary 17]{Abr}}
Let $\tau = \{E_1<\dots < E_r\} \in \tilde{\Delta}$ and set $\EE(\tau) = \{E_i, E_i^\perp \mid 1 \leq i\leq r\}$ then $\tilde{\Delta}^0(\tau) = \tilde{T}_{\EE(\tau)}(V)$.
\end{proposition}
\section{Details for type $A_n$}
\label{app: An subcomplexes sec}

In this section, we fill in the details of the proof that the complex opposite the fundamental chamber in a spherical building of type $A_n$ has a non-Abelian cone function with bounded cone radius, by showing that the class $\mathbf{C}_A$ from \Cref{def: class C_A} satisfies the assumptions of \Cref{thm: generalized induction argument}.
The argument will follow \cite[Chapter 6.2]{CSEofKMS} closely.

We will need some facts that appear in the proof of \cite[Proposition 12]{Abr}.  
\begin{fact}
\label{Abr fact1}
If $T_{\EE}(V) \in \mathbf{C}_A$, then there exists $\ell \in X_{\mathcal{E}} (V)$ such that $\ell$ is a line,  i.e.,  $\dim (\ell) =1$. 
\end{fact}

Let $\ell \in X_{\mathcal{E}} (V)$ be a line.  Define 
$$Y_0 = \lbrace  U  \in T_{\mathcal{E}} (V) : \ell + U \in T_{\mathcal{E}} (V) \rbrace,$$
and further define $\kappa_0$ to be the subcomplex of $T_{\mathcal{E}} (V)$ spanned by $Y_0$.

Assume now that $\dim V = n+2$ and hence $\dim T_{\EE}(V) = n$ (and $\flag X$ is a building of type $A_{n+1}$). 
For $1 \leq i \leq n+1$,  define $Y_i \subseteq X_{\mathcal{E}} (V)$ as follows: 
$$Y_i = \lbrace U \in X_{\mathcal{E}} (V)\mid \dim (U) \geq n+2-i \text{ or } U \in Y_0 \rbrace.$$
For every such $i$,  let $\kappa_i$ be the subcomplex of $T_{\mathcal{E}} (V)$ spanned by $Y_i$.  
\begin{fact}
\label{Abr fact2}
We denote $Y = X_{\mathcal{E}} (V), \kappa = T_{\EE}(V)$ and for every $U \in T_{\mathcal{E}} (V)$ we denote the link of $U$ by $\lk_{\kappa}(U)$.  For every $1 \leq i \leq n+1$ and every $U \in Y_i \setminus Y_{i-1}$,  it holds that $\lk_\kappa(U) \cap \kappa_{i-1} = T_{\mathcal{E}' } (U) * T_{\overline{\mathcal{E}}} (V /U) $ such that 
\begin{itemize}
\item The set $\mathcal{E}'$ is a set of subsets of $U$ satisfying $\sum_{j=1}^{n+1} {{n}\choose{j-1}} e'_{j} \leq \lvert K \rvert$, hence $ T_{\mathcal{E}' } (U) \in \mathbf{C}_A$.
\item The set $\overline{\mathcal{E}}$ is a set of subsets of $V/U$ that satisfying $\sum_{j=1}^{n+1} {{n}\choose{j-1}} \overline{e}_{j} \leq \lvert K \rvert$ and hence $ T_{\mathcal{E}' } (V/U) \in \mathbf{C}_A$.
\end{itemize}
\end{fact}

\begin{lemma}
\label{X0 bound lemma}
There is a NAC $\Cone_{\kappa_0}$ such that 
$$\Rad (\Cone_{\kappa_0}) \leq n+2 .$$
\end{lemma}

\begin{proof}
Let $\Gamma_0$ be the subcomplex of $\kappa_0$ spanned by $Z_0=\lbrace  U  \in X_{\mathcal{E}} (V)  : \ell \leq U \rbrace$.  For $1 \leq i \leq n+1$,  we denote 
$$Z_i = \lbrace U   \in Y_0 \mid   \dim (U) \leq i \text{ or } U \in Z_0 \rbrace.$$
Further denote $\Gamma_i = \flag Z_i$. Note that $\Gamma_{n+1} = \kappa_0$.  We will show by induction that for every $0 \leq i \leq n+1$, there is a NAC $\Cone_{\Gamma_i}$ such that $\Rad ( \Cone_{\Gamma_i}) \leq i+1$.

For $i=0$,  we note that  $\Gamma_0$ is the ball of radius $1$ around $\ell$ in $\kappa$ and can be written as $\Gamma_0 = \lbrace \ell \rbrace * \lk_{\kappa}(\ell)$.  By Proposition \ref{cone of a join - single vertex prop},  there is NAC $\Cone_{\Gamma_0}$ such that   $\Rad (\Cone_{\Gamma_0}) \leq 1$.  

Let $1 \leq i \leq n+1$ and assume there is a NAC $\Cone_{\Gamma_{i-1}}$ such that $\Rad ( \Cone_{\Gamma_{i-1}}) \leq i$.  We will apply \Cref{cor: adding vertices NAC version} in order to bound the cone radius of $\Gamma_i$.  We note that any two $U_1,  U_2 \in Y_i \setminus Y_{i-1}, U_1 \neq U_2$ are of the same dimension and thus are not connected by an edge.  Let $U \in Z_i \setminus Z_{i-1}$.  We will show that $\lk_{\kappa_0}(U) \cap \Gamma_{i-1}$ is a join of the vertex $\lbrace U + \ell \rbrace$ with the complex spanned by all the other vertices in $\lk_{\kappa_0}(U) \cap \Gamma_{i-1}$.  

Note that $\lk_{\kappa_0}(U) \cap \Gamma_{i-1}$ is a clique complex,  thus it is enough to show that for every vertex $U'$ in $\lk_{\kappa_0}(U) \cap \Gamma_{i-1}$ that is not $\ell + U$,  there is an edge connecting $U'$ and $\ell + U$.  Let $U'$ be such vertex.  

First, we will deal with the case were $\dim (U') < i$.  In that case (using the fact that $\dim (U) =i$),  $U'$ being in $\lk_{\kappa_0}(U)$ implies that $U' \leq U$ and thus $U' \leq U + \ell$,  i.e., $U'$ is connected by an edge to $U + \ell$ as needed.  

Second, assume that $\dim (U') >i$.  From the definition of $Z_{i-1}$,  it follows that $U' \in Z_0$,  i.e., that $\ell \leq U'$.  Also,  $U'$ is in $\lk_{\kappa_0}(U)$ and thus $U \leq U'$.  It follows that $U+\ell \leq U'$,  i.e.,  $U'$ is connected by an edge to $U + \ell$ as needed.  

By Proposition \ref{cone of a join - single vertex prop},  for every $U \in Z_i \setminus Z_{i-1}$,  there is 0-cone function $\Cone_{\lk_{\kappa_0}(U) \cap \Gamma_{i-1}}$ such that $\Rad(\Cone_{\lk_{\kappa_0}(U) \cap \Gamma_{i-1}}) \leq 1$.  

By \Cref{cor: adding vertices NAC version},  it follows that there is a NAC  $\Cone_{\Gamma_i}$ such that $\Rad (\Cone_{\Gamma_i}) \leq i+1$ as needed. 
\end{proof}

\begin{corollary} \label{app: cor An case}
For every $n \in \mathbb{N} $ there is a constant $\mathcal{R} (n)$ such that for every $T_{\EE}(V) \in \mathbf{C}_A$ with $\dim T_{\EE}(V) = n$ there in a NAC $\Cone_{T_{\mathcal{E}} (V)}$ such that 
$$\Rad (\Cone_{T_{\mathcal{E}} (V)}) \leq \mathcal{R} (n).$$
In particular, if $\Delta$ is a building of type $A_n$ ($n \geq 2$) over a field $K$ with $\lvert K \rvert \geq 2^{n-1}$ and $a\in \Delta$ a simplex, then there exists a NAC $\Cone_{\Delta^0(a)}$ of $\Delta^0(a)$ such that  
$$\Rad(\Cone_{\Delta^0(a)}) \leq \mathcal{R} (n-1).$$
\end{corollary}

\begin{proof}
The result follows from \Cref{thm: generalized induction argument} together with \Cref{X0 bound lemma} and Fact \ref{Abr fact2}. In the notation of \Cref{thm: generalized induction argument}, we have $f(n) = n+2$ and $\ell_n = n+1$.
\end{proof}

\section{Details for type $C_n$} \label{app: Cn subcomplexes section}
In this section, we fill in the details of the proof that the complex opposite the fundamental chamber in a spherical building of type $C_n$ has a non-Abelian cone function with bounded cone radius. The class $\textbf{C}$ that will cover the case of opposite complexes in buildings of type $C_n$ will be the union of class $\mathbf{C}_A$ (see \Cref{def: class C_A}) and $\mathbf{C}_C$ (see \Cref{def: class C_C}).

Again, the argument will follow \cite[Chapter 6.3]{CSEofKMS} closely.

\begin{proposition}\label{C satisfies assumptions}
The class $\textbf{C}$ satisfies the requirements of \Cref{thm: generalized induction argument}.
\end{proposition}

\begin{corollary}\label{app: cor Cn case}
Let $\Delta=\flag X(V)$ be a classical $C_n$ building as described in \Cref{subsubsec: geom descr C_n} (note that $\dim \Delta = n-1$). Assume that $K$ is infinite or $K=\mathbb{F}_q$ and 
\begin{itemize}
\item $q \geq 2^{2n-2}$ if $f$ is alternating and $Q=0$,
\item $q \geq 2^{4n-4}$ if $m=2n$ and $\sigma \neq \operatorname{id}$,
\item $q \geq 2^{2n-1}$ if $m=2n+1$ or $2n+2$.
\end{itemize}
In particular, if $\Delta$ is the building coming from the BN-pair of a Chevalley group of type $B_n$ or $C_n$ over the field $\mathbb{F}_q$ we require $q \geq 2^{2n-1}$. 

Then for any $a \in \Delta$ there exists a NAC $\Cone_{\Delta^0(a)}$ of $\Delta^0(a)$ with 
$$\Rad(\Cone_{\Delta^0(a)}) \leq \mathcal{R}(n-1)$$
where $\mathcal{R}(n-1)$ is as in \Cref{thm: generalized induction argument} with $f(n)=2n+1$ and $\ell_n = 2n+2$.
\end{corollary}

The main part of the proof of \Cref{C satisfies assumptions} is covered by the following lemma.
\begin{lemma} \label{kappa_0 C_n case}
Let $(T_{\mathcal{E}}(V), \flag X) \in \textbf{C}_C$. Set $n =  \dim \flag X \geq 1$. By \cite{Abr}, there exists a 1-dimensional subspace $\ell \in X_{\mathcal{E}}(V)$. Define $X_{0}=\left\{ U \in X \mid U \text{ satisfies }1,2 \text{ or }3 \right\}$ where 
\begin{enumerate}
\item $\ell \leq U$ and $U \pitchfork \mathcal{E}$;
\item $\ell \nleq U \leq \ell^\perp$ and $U \pitchfork \mathcal{E} \cup (\mathcal{E}+\ell)$;
\item $U\nleq \ell^\perp, \dim U >1$ and $U \pitchfork \mathcal{E}\cup(\mathcal{E} \cap \ell^\perp) \cup (\mathcal{E} + \ell) \cup ((\mathcal{E} \cap \ell^\perp)+\ell)$.
\end{enumerate}
Then $\kappa_{0} := \flag X_{0}$ has a NAC $\Cone_{\kappa_0}$ such that 
$$\Rad(\Cone_{\kappa_0}) \leq 2n+3.$$
\end{lemma}
\begin{proof}
Let $(T_{\mathcal{E}}(V), \flag X) \in \textbf{C}_C$ as above. Thus $\dim \flag X = n$ and $\flag X$ is of type $C_{n+1}$. In particular, any totally isotropic subspace of $V$ has dimension at most $n+1$. 

We will apply \Cref{cor: adding vertices NAC version} inductively to two different filtrations of $\kappa_0$. The first one is defined as follows.
\begin{align*}
Y_0 &= \left\{ U \in X \text{ satisfying } 1. \right\} \\
Y_i &= Y_{i-1} \cup \left\{U \in X \mid \dim U = i, U \text{ satisfies } (ii)  \right\}=\left\{U \in X \mid \dim U \leq i, U \text{ satisfies } (ii)  \right\} \cup Y_{0} \\ & \text{ for } 1 \leq i \leq n+1. 
\end{align*} 
We set $\Gamma_i = \flag Y_i$.

The second filtration is given by
\begin{align*}
Z_0 &= Y_{n+1} \\
Z_i &= Z_{i-1} \cup \left\{ U \in X \mid \dim U = n+1-(i-1), U \text{ satisfies } (iii) \right\} \\
&= Z_{0} \cup \left\{ U \in X \mid \dim U \geq n+1-(i-1), U \text{ satisfies } (iii) \right\}
\end{align*}
for $1 \leq i \leq n+1$. We set $\zeta_i = \flag Z_i$. 
We have that $\Gamma_{n+1} = \zeta_0$ and $\zeta_{n+1} = \kappa_0$.

Next, we want to show that for each $\Gamma_i$ there exists a NAC $\Cone_{\Gamma_i}$ such that $\Rad(\Cone_{\Gamma_i}) \leq i+1$. We do this by induction on $i$. 

For $i=0$ we have $\Gamma_0 = \st_{\kappa}(\ell) = \{\ell\} * \operatorname{lk}_\kappa(\ell)$. Hence by \Cref{cor: summary of join results} there exists a NAC $\Cone_{\Gamma_0}$ with $\Rad(\Cone_{\Gamma_0}) \leq 1$.

For the following step, fix $1 \leq i \leq n+1$. We want to apply \Cref{cor: adding vertices NAC version} to the following set-up. 
Here $\kappa_{0}$ corresponds to the complex called $X$ in \Cref{cor: adding vertices NAC version}, $\Gamma_{i-1}$ corresponds to $X'$ , and $ Y_{i} \setminus Y_{i-1} = W_{i}$. 
We check that the conditions of \Cref{cor: adding vertices NAC version} are satisfied.

For Condition (i) note that $W_{i} \cap \Gamma_{i-1}=\emptyset$ by definition.
Condition (ii) holds since all subspaces in $W_{i}$ have the same dimension, hence they are not connected in $\Delta$ and thus not in $\kappa$.

For Condition (iii) let $w\in W_{i}$. Then $\lk_{\kappa_{0}}(w) = \flag\left\{ u \in \kappa_{0}(0) \mid u < w \text{ or } w > u  \right\}$. Furthermore, for $u \in Y_{i-1}$ we have $u+\ell \in X_{0}$, since $u \pitchfork \mathcal{E} + \ell \iff u + \ell \pitchfork \mathcal{E}$, see \cite{Abr}. 
Since $\ell \leq w+\ell$, we have $w + \ell \in Y_{0} \subset Y_{i-1}$. Next, we want to show that $\lk_{\kappa_{0}}(w) \cap \Gamma_{i-1}$ is a join of $\{w + \ell\}$ with the flag complex of all other vertices from $\lk_{\kappa_{0}}(w) \cap \Gamma_{i-1}$. Since $\lk_{\kappa_{0}}(w) \cap \Gamma_{i-1}$ is a clique complex, it suffices to check that every vertex different from $w + \ell$ is adjacent to $w + \ell$. Let $u \in (\lk_{\kappa_{0}}(w) \cap \Gamma_{i-1})(0) \setminus \{w + \ell\}$. Then we differentiate two cases: 
\begin{itemize}
\item Case 1: $\dim u \leq i-1$: then $u \leq w \leq w + \ell$.
\item Case 2: $\dim u >i$ then $u \in Y_{0}$, thus $w \leq u$ and $\ell \leq u$ thus $w + \ell \leq u$.
\end{itemize}
Hence any other vertex is connected to $w+ \ell$ by an edge and thus $$\Gamma_0 = \{w + \ell \} *\flag  \left((\lk_{\kappa_{0}}(w) \cap \Gamma_{i-1})(0) \setminus \{w + \ell\}\right).$$ By \Cref{cone of a join - single vertex prop} we have a 0-cone function $\operatorname{Cone}_{\lk_{\kappa_{0}}(w) \cap \Gamma_{i-1}}$ with $\operatorname{Rad}_0(\operatorname{Cone}_{\lk_{\kappa_{0}}(w) \cap \Gamma_{i-1}}) \leq 1$. 

Thus we can apply \Cref{cor: adding vertices NAC version} and get a NAC $\Cone_{\Gamma_i}$ satisfying $\operatorname{Rad}(\operatorname{Cone}_{\Gamma_{i}}) \leq i+1$.

We treat the second filtration similarly. We want to show that for every $0 \leq i \leq n+1$ there is a NAC $\Cone_{\zeta_i}$ such that $\Rad(\Cone_{\zeta_i}) \leq n+2 +i$. We again proceed by induction on $i$. 

For $i=0$ we have that $\zeta_{0}=\Gamma_{n+1}$ has cone function with radius $\leq n+2$ by the above reasoning.

Now fix $1\leq i \leq n+1$. We want to use \Cref{cor: adding vertices NAC version} with 
$\kappa_{0}$ corresponding to $X$, $\zeta_{i-1}$ to $X'$, and $ Z_{i}\setminus Z_{i-1} = W$.
Then conditions (i) and (ii) of \Cref{cor: adding vertices NAC version} are satisfied by the same argument as above. 

For Condition (iii), fix $w \in W$ and consider $\lk_{\kappa_{0}}(w) \cap \zeta_{i-1} = \flag \left\{ u \in Z_{i-1} \mid u <w \text{ or } u>w \right\}$.
We show that every vertex (different from $w \cap \ell^\perp$) is connected to $w \cap \ell^\perp$ hence 
$\lk_{\kappa_{0}}(w) \cap \zeta_{i-1}= \left\{ w \cap \ell^\perp \right\}* \flag \left\{ u \in Z_{i-1} \setminus \left\{ w \cap \ell^\perp \right\} \mid u <w \text{ or } u>w \right\}$.
Thus by \Cref{cone of a join - single vertex prop}, $\lk_{\kappa_{0}}(w) \cap \zeta_{i-1}$ has a 0-cone function with cone radius bounded by 1. 
Applying \Cref{cor: adding vertices NAC version} we get a NAC $\Cone_{\zeta_i}$ satisfying
$$\operatorname{Rad}(\operatorname{Cone}_{\zeta_{i}}) \leq \operatorname{Rad}(\operatorname{Cone}_{\zeta_{i-1}}) +1 \leq n+2 +i.$$

To show that every vertex (different from $w \cap \ell^\perp$) is connected to $w \cap \ell^\perp$, let $u \in \lk_{\kappa_{0}}(w) \cap \zeta_{i-1}(0) = \left\{ u \in Z_{i-1} \mid u <w \text{ or } u>w \right\}$. If $w \leq u$ then $w\cap \ell^\perp \leq u$ and thus $\{u,w\cap \ell^\perp \} \in \lk_{\kappa_{0}}(w) \cap \zeta_{i-1}$. If $u \leq w$, then $\dim u < \dim w = n+1-i+1$ hence $u \in Z_0$. Therefore $u$ satisfies either 1 or 2. But if it satisfies 1, then $\ell \leq u \leq w$, a contradiction to $w \in Z_i \setminus Z_{i-1}$. Hence $u$ satisfies 2, and in particular $u\leq \ell ^\perp$. Thus $u \leq w \cap \ell^\perp$ and hence $\{u,w\cap \ell^\perp \} \in \lk_{\kappa_{0}}(w) \cap \zeta_{i-1}$.

Since $\kappa_0 = \zeta_{n+1}$ we get the desired result. 
\end{proof}

\textit{Proof of \Cref{C satisfies assumptions}. }
We proof that every $(\kappa,\Delta) \in \mathbf{C}$ has a filtration satisfying the assumptions of \Cref{thm: generalized induction argument} by induction on $n=\dim \Delta$.

Let $n=0$. Hence $\Delta$ will be a building of type $A_1$ or $C_1$. By \cite{Abr} there exists at least on 1-dimensional subspace $\ell \in X_{\EE}(V)$ and thus $\kappa$ will be non-empty.   

For the inductive step, fix $n \geq 1$ and set $f(n) = 2n+3$.
We distinguish between two cases: 

If $(\kappa,\Delta)\in \mathbf{C}_A$, \Cref{app: An subcomplexes sec} gives the desired filtration with $f_A(n) = n+2 \leq f(n)$ and length of the filtration being $n+1$.

If $(\kappa, \Delta)\in \mathbf{C}_C$, we can use the notation as in \Cref{def: class C_C}, e.g. $(\kappa,\Delta) = (T_{\EE}(V), \flag X)$. As in \Cref{kappa_0 C_n case} we fix a 1-dimensional subspace $\ell \in X_{\EE}(V)$ and define $X_{0}=\left\{ U \in X \mid U \text{ satisfies }1,2 \text{ or }3 \right\} $ (where 1,2,3 are as in the lemma) and $\kappa_0=\flag X_0$. Hence, by \Cref{kappa_0 C_n case}, there exists a NAC $\Cone_{\kappa_0}$ such that 
$$\Rad(\Cone_{\kappa_0}) \leq 2n+3 = f(n).$$
The remaining filtration is defined as in the Abelian case, see \cite[Proposition 6.7]{CSEofKMS}, and hence satisfies the assumptions of \Cref{thm: generalized induction argument}.
\qed

\section{Details for type $D_n$} \label{app: sec dn case}

In this section, we fill in the details of the proof that the complex opposite the fundamental chamber in a spherical building of type $D_n$ has a non-Abelian cone function with bounded cone radius.

The argument will follow \cite[Chapter 6.4]{CSEofKMS} closely which itself relies on the work in \cite[Section 7]{Abr} and we also refer to that book for more details. 
We will define three subclasses $\mathbf{C}_1,\mathbf{C}_2,\mathbf{C}_3$ and first show that $\mathbf{C}_1,\mathbf{C}_2$ satisfy the assumptions of \Cref{thm: the big induction for the D_n case} to then conclude that $\mathbf{C}_D = \mathbf{C}_1 \cup \mathbf{C}_2 \cup \mathbf{C}_3$ satisfies them as well.

We fix $n$ to be the Witt index of $(V,Q,f)$ hence $\dim V = 2n$ and we are dealing with a $D_n$ building, which is $n-1$-dimensional.
Recall $X= \{0<U<V \mid U \text{ totally isotropic}\}$.
\begin{definition}{\cite[Definition 12]{Abr}}
Let $U,E <V$, $\dim  U <n, \dim E = n$. Assume that $E$ is not totally isotropic. Then 
$$U@E :\iff U^\perp \cap E \text{ is not totally isotropic.}$$
Set 
\begin{align*}
\mathcal{U}_n &= \mathcal{U}_n(V) = \{A<V \mid \dim A = n, A\notin X\} \\
\mathcal{M}(A) &= \{M<A \mid \dim M = n-1, M \in X\} \text{ for any } A \in \mathcal{U}_n.
\end{align*}
\end{definition}

\begin{definition} \label{def: class 1 D_n case}
We define the class $\mathbf{C}_1$ to consist of complexes $\flag X_{\EE; \mathcal{F}}(U;V)$ where 
\begin{itemize}
\item $(V,Q,f)$ is thick pseudo-quadratic space of Witt index $n$, in particular $\dim V = 2n$;
\item $U\in X$ with $\dim U = k \geq 2$;
\item $\EE$ a finite set of subspaces of $U$, set $e_j = \lvert \EE_j \rvert, 1 \leq j \leq k-1$ like before;
\item Let $\mathcal{F} = \mathcal{F}_1 \supset \mathcal{F}_2 \supset \dots \supset \mathcal{F}_{k-1}$ be finite subsets of $\mathcal{U}_n(V)$ satisfying
\begin{enumerate}
\item $U \cap F \in \EE \cup \{0\}$ for all $F\in \mathcal{F}$;
\item $\mathcal{F}_j^\perp = \mathcal{F}_j$;
\item $\dim(U\cap F) \leq k-1-j$ for all $F \in \mathcal{F}_j, 1 \leq j \leq k-1$.
\end{enumerate}
\item Assume that $\lvert K \rvert \geq \sum_{j=1}^{k-1} \binom{k-2}{j-1} e_j + 2^{k-1}s$, where $s = \lvert \mathcal{F} \rvert$.
\end{itemize} 
Set 
$$X_{\EE; \mathcal{F}}(U;V):= \{0< W < U \mid W \pitchfork_U \EE \text{ and } W @_V \mathcal{F}_j \text{ for } \dim W = j \}.$$
\end{definition}

\begin{lemma}\label{lemma: Y_0 for class 1 D_n case}
Let $\kappa = \flag X_{\EE; \mathcal{F}}(U;V) \in \mathbf{C}_1$ with $ k = \dim U \geq 2$, in particular $\kappa$ has dimension $k-2$. Then there exists a 1-dimensional subspace $\ell \in X_{\EE; \mathcal{F}}(U;V)$. We set $Y= X_{\EE; \mathcal{F}}(U;V)$ and define
$$Y_0 := \{A \in Y \mid A + \ell \in Y\}.$$
Then $\kappa_0 = \flag Y_0$ has a NAC $\Cone_{\kappa_0}$ with radius at most $k$.
\end{lemma}

\begin{proof}
Set $B_0 :=\{A \in Y \mid \ell \leq A\}, \beta_0 = \flag B_0$. Then $\beta_0 =\{\ell\} * \lk_Y(\{\ell\})$ and hence it has a NAC with cone radius bounded by 1.
Next, set $B_i= \{A \in Y_0 \mid \ell \leq A \text{ or } \dim A \leq i \}, \beta_i = \flag B_i$ for $1 \leq i \leq k-1$. Then all subspaces in $B_i \setminus B_{i-1}$ have dimension $i$. Let $A \in B_i \setminus B_{i-1}$. 
Then $\lk_{\beta_i}(A)\cap \beta_{i-1}= \{A+\ell\} * \flag \{W \in B_{i-1}\mid B \neq A + \ell, A <W \text{ or } W <A \}$ by the same argument as in \Cref{X0 bound lemma}. 
Applying \Cref{cor: adding vertices NAC version} inductively, we get a NAC $\Cone_{\beta_{k-1}} = \Cone_{ \kappa_0}$ with 
$$\Rad(\Cone_{\kappa_0}) \leq k.$$
\end{proof}

\begin{lemma}\label{lemma: class 1 D_n case}
The class $\mathbf{C}_1$ satisfies the assumptions of \Cref{thm: the big induction for the D_n case}. In particular, any complex in $\mathbf{C}_1$ has a NAC with cone radius depending only on the dimension of the complex.  
\end{lemma}

\begin{proof}
The result follows from \Cref{rmk: from Abelian to non-Abelian in Dn}, \Cref{lemma: Y_0 for class 1 D_n case} and \cite[Lemma 6.17]{CSEofKMS}.
\end{proof}

\begin{definition}\label{def: class 2 D_n case}
We define the class $\mathbf{C}_2$ to consist of complexes $\flag Z_{\EE; \mathcal{F}}(U;V)$ where
\begin{itemize}
\item $(V,Q,f)$ is a thick pseudo-quadratic space of Witt index $n$, in particular $\dim V = 2n$;
\item $U\in X$ with $\dim U = k \geq 2$;
\item $\EE$ a finite set of subspaces of $U$, set $e_j = \lvert \EE_j \rvert, 1 \leq j \leq k-1$ like before;
\item $\mathcal{F} \subset \mathcal{U}_n(V)$ with $\lvert \mathcal{F}\rvert = s < \infty$ such that 
\begin{enumerate}
\item $U\cap \mathcal{F} \subseteq \EE \cup \{0\}$,
\item $U \cap F^\perp = 0$ and $\dim(U \cap F) \leq 1$ for all $F \in \mathcal{F}$.
\end{enumerate}
\item Assume $\lvert K \rvert \geq \sum_{j=1}^{k-1} \binom{k-2}{j-1} e_j + 2s$.
\item Set 
$$Z_{\EE; \mathcal{F}}(U;V) := \{0<W<U \mid W \pitchfork_U \EE \text{ and } W @_V \mathcal{F}\}.$$
\end{itemize}
\end{definition}

\begin{lemma} \label{lemma: Y_0 for class 2 D_n case}
Let $\kappa=\flag Z_{\EE; \mathcal{F}}(U;V) \in \mathbf{C}_2$, $\dim U = k \geq 2$, in particular $\dim \kappa = k-2$. Set $Y = Z_{\EE; \mathcal{F}}(U;V)$. Then there exists $H <U$ with $\dim H = k-1$ and $H \in Y$. We define
$$Y_0:= \{B\in Y \mid B \cap H \in Y\}.$$
Then $\kappa_0 = \flag Y_0$ is $(k-2)$-dimensional and has a NAC $\Cone_{\kappa_0}$ with 
$$\Rad(\Cone_{X_0}) \leq k.$$
\end{lemma}

\begin{proof}
The proof is similar to the proofs of \Cref{X0 bound lemma} and \Cref{lemma: Y_0 for class 1 D_n case}, the difference being that we now look at the dual version, where we fix a hyperplane instead of a line.

Set $B_0 = \{A \in Y \mid A \leq H \}, \beta_0 = \flag B_0$. Then $\beta_0 = \{H \} * \lk_{Y_0}(\{H\})$ hence it has a NAC with cone radius bounded by $1$.

Next, set $B_i = \{A \in Y_0 \mid A \leq H \text{ or } \dim A \geq k-i\}, \beta_i = \flag B_i$ for $1 \leq i \leq k-1$. Then the subspaces in $B_i\setminus B_{i-1}$ are all of dimension $k-i$. Let $A \in B_i \setminus B_{i-1}$, then $\lk_{\beta_i}(A) \cap \beta_{i-1} = \{A \cap H \} * \flag \{W \in B_{i-1} \mid W \neq A \cap H, W < A \text{ or } A < W \} $. Hence it has a 0-cone function with cone radius bounded by 1. Thus we can apply \Cref{cor: adding vertices NAC version} inductively and obtain a NAC $\Cone_{B_{k-1}} = \Cone_{\kappa_0}$ with 
$$\Rad(\Cone_{\kappa_0}) \leq k.$$

\end{proof}

\begin{lemma}\label{lemma: class 2 D_n case}
The class $\mathbf{C}_2$ satisfies the assumptions of \Cref{thm: the big induction for the D_n case}. In particular, any complex in $\mathbf{C}_2$ has a NAC with cone radius bounded depending only on its dimension.  
\end{lemma}

\begin{proof}
The result follows from \Cref{rmk: from Abelian to non-Abelian in Dn}, \Cref{lemma: Y_0 for class 2 D_n case} and \cite[Lemma 6.20]{CSEofKMS}.
\end{proof}

Next, we will define a third class, which will be the one containing the relevant opposite complexes, see also \Cref{rmk: connection to op compl D_n}.

\begin{definition}\label{def: class 3 D_n case}
We define the class $\mathbf{C}_3$ to consist of complexes of the form $\flag Y_{\EE}(V)$ where
\begin{itemize}
\item $(V,Q,f)$ is a thick pseudo-quadratic space of Witt index $n$, in particular $\dim V = 2n$;
\item $\EE = \EE^\perp$ is a finite set of subspaces of $V$, set $\EE_j= \{E \in \EE : \dim E=j\}, e_j = \lvert \EE_j \rvert$;
\item $\lvert K \rvert \geq 2 \sum_{j=1}^{2n-1}\binom{2n-2}{j-1}e_j$;
\item $\widehat{\EE}_n:= \EE_n \cap \mathcal{U}_n(V)$;
\item $X = \{0<U<V \mid U \text{ is totally isotropic} \}$.
\end{itemize}
Set $$Y_{\EE}(V) :=\{U \in X \mid U \tilde{\pitchfork}\EE \text{ and if } \dim U < n \text{ then } U@ \hat{\EE}_n \}.$$
\end{definition}

\begin{lemma}\label{lemma: Y_0 for class 3 D_n case}
Let $\kappa= \flag Y_{\EE}(V) \in \mathbf{C}_3$ with $\dim V = 2n$. Then there exists a 1-dimensional subspace $\ell \in Y=Y_{\EE}(V)$. 
 Set
$$
\begin{aligned}
& \mathcal{F}:=\mathcal{E} \cup\left(\mathcal{E} \cap \ell^{\perp}\right) \cup(\mathcal{E}+\ell) \cup\left(\left(\mathcal{E} \cap \ell^{\perp}\right)+\ell\right) \text { and } \\
& \widehat{\mathcal{F}}_n:=\mathcal{F}_n \cap \mathcal{U}_n(V). 
\end{aligned}
$$
Consider the following conditions
\begin{enumerate}
\item $\ell \leq U, U \tilde{\pitchfork} \mathcal{E}$ and if $\dim U <n$ then $U @ \widehat{\mathcal{E}}_n$; 
\item $\ell \not \leq U \leq \ell^{\perp} \text{ and } \begin{cases} U \pitchfork \mathcal{E}, U @ \widehat{\mathcal{E}}_n & \text{ if } \dim U = n-1, \\
U \pitchfork \mathcal{E} \cup(\mathcal{E}+\ell), U @ \widehat{\mathcal{F}}_n & \text{ if } \dim U \leq n-2; \end{cases}$
\item $U \not \leq \ell^{\perp}, \operatorname{dim} U>1, U \widetilde{\pitchfork} \mathcal{F}$ and if $\dim U <n$ then $U @ \widehat{\mathcal{F}}_n$.
\end{enumerate}
We define $Y_0:=\{U \in X \mid U$ satisfies 1.,2., or 3. $\}$. Then $\kappa_0 = \flag Y_0$ has a NAC function $\Cone_{\kappa_0}$ with 
$$\Rad_j(\Cone_{\kappa_0}) \leq 2n+1 , \ 0 \leq j \leq 1.$$
\end{lemma}
\begin{proof}
Let $\kappa= \flag Y_{\EE}(V) \in \mathbf{C}_3$ with $\dim V = 2n$. By Step 1. \cite[Proposition 14]{Abr} there exists a 1-dimensional subspace $\ell \in Y$. Let $Y_0$ be defined as above. Recall that $X = \{0<U<V \mid U \text{is totally isotropic subspace}\}$. We define a filtration for $Y_0$ and we want to apply \Cref{cor: adding vertices NAC version} inductively.
We will use the fact stated in Step 3 of \cite[Proposition 14]{Abr}, that for every $U \in Y_0$, we have that $U\cap \ell^\perp \in Y_0$ and $(U\cap \ell^\perp)+\ell \in Y_0$.

Set
\begin{align*}
A_0 &:=\{U \in X \mid U \text{ satisfies 1.} \} \\
A_i &:= \{U \in X \mid U \text{ satisfies 1. or } (\dim U \leq i \text{ and } U \text{ satisfies 2.}) \}, \ 1 \leq i \leq n; \\
\alpha_i &:= \flag A_i, 0 \leq i \leq n.
\end{align*}
Then $A_0 = \{\ell\}* \flag \{U \in A_0 \mid U \neq \ell\}$ and hence has a NAC with cone radius bounded by 1. 
Clearly elements in $A_i \setminus A_{i-1}$ have the same dimension an are hence not adjacent. Let $U \in A_i \setminus A_{i-1}$. Then  
$$\lk_{\alpha_i}(U)\cap \alpha_{i-1} = \flag \{W \in X \mid (W \in A_0 \text{ and } U <W) \text{ or } (W \in A_{i-1} \text{ and } W <U)\}. $$ 
We want to show that each subspace in $\{W \in X \mid (W \in A_0 \text{ and } U <W) \text{ or } (W \in A_{i-1} \text{ and } W <U)\}$ is either contained in or contains $U + \ell$ (note that we have $U\cap \ell^\perp + \ell = U +\ell \in A_0 \subseteq A_{i-1}$).
If $W \in A_0$ and $U<W$ then $\ell \leq W$ and hence $U + \ell \leq W$. If $W \in A_{i-1}$ and $W<U$ then $W< U+ \ell$ as needed. 
The second part of the filtration is defined as follows. We set
\begin{align*}
B_i&:= \{U \in X \mid  U \in A_n \text{ or } (U \text{ satisfies } 3. \text{ and } \dim U \geq n-i+1 \},\\
\beta_i &:= \flag B_i, 0 \leq i \leq n.
\end{align*}
Then $A_n =B_0$. Let $U \in B_i \setminus B_{i-1}$. We want to show that $\lk_{\beta_i}(U) \cap \beta_{i-1}$ can be written as the join of $\{U \cap \ell ^\perp \}$ and $\flag \{W \in B_{i-1} \mid W\neq U \cap \ell^\perp, W < U \text{ or } U <W\}$.
First, note that $U \cap \ell^\perp \in Y_0$ and hence $U \cap \ell^\perp \in A_n \subseteq B_{i-1}$. Let $W \in B_{i-1}\setminus \{U \cap \ell^\perp\}$. If $W<U$ then $\dim W < n-i$, hence $W \in A_n\setminus A_0$, in particular $W \leq \ell^\perp$. Hence $W = W \cap \ell^\perp \leq U \cap \ell^\perp$.
If $U<W$ then $U \cap \ell^\perp \leq U \leq W$.   
Thus we get a filtration of $\kappa_0$ of length $2n$. By applying \Cref{cor: adding vertices NAC version} at each step, we get a NAC $\Cone_{\kappa_0}$ with 
$$\Rad_j(\Cone_{\kappa_0}) \leq 2n+1, 0 \leq j \leq 1.$$
  
\end{proof}

We now have all the necessary ingredients to show that the relevant class of simplicial complexes satisfies the conditions of \Cref{thm: the big induction for the D_n case}. The prove will be similar to the previous proofs of similar results. 

\begin{theorem} \label{thm: class for D_n case}
The class $\mathbf{C}_D= \mathbf{C}_1 \cup \mathbf{C}_2 \cup \mathbf{C}_3$ satisfies the conditions of \Cref{thm: the big induction for the D_n case}.
\end{theorem}
\begin{proof}
Let $\kappa = \flag Y \in \mathbf{C}_D$. If $\kappa \in \mathbf{C}_1 \cup \mathbf{C}_2$ then the result follows from \Cref{lemma: class 1 D_n case} and \Cref{lemma: class 2 D_n case}.

Now assume $\kappa=\flag Y_{\EE}(V) \in \mathbf{C}_3$. In this case the proof follows closely the proof of \cite[Proposition 14]{Abr}.
First of all, notice that we again have a 1-dimensional subspace $\ell \in Y_{\EE}(V) =:Y$ (see Step 1 \cite[Proposition 14]{Abr}).
Let $n=\frac{1}{2}\dim V = \dim \kappa +1$. We will again prove the existence of a filtration by induction on $n$, but this time starting with $n=2$. 

Let $n=2$, hence $\dim \kappa = 1$. Step 2. in \cite[Proposition 14]{Abr} shows that each vertex $U \in Y$ can be connected to $\ell$ by a path of length at most 5. Thus these paths define a $0$-cone function $\Cone_\kappa$ of $\kappa$ with radius at most 5.

Now assume $n\geq 3$. We define $\kappa_0 = \flag Y_0$ as in \Cref{lemma: Y_0 for class 3 D_n case}. In particular, we know that there exists a NAC $\Cone_{\kappa_0}$ with $\Rad_j(\Cone_{\kappa_0}) \leq 2n+1 = f(n)$ for all $0 \leq j\leq 1$. 

The fact that the rest of the filtration satisfies the conditions of \Cref{thm: the big induction for the D_n case} follows from the corresponfing proof in the Abelian case, see \cite[Theorem 6.23]{CSEofKMS}.
\end{proof}
The following remark from \cite{Abr} gives the connection between $X_\EE(V)$ and $Y_\EE(V)$. 
\begin{remark} \label{rmk: connection to op compl D_n}
Let $\tau = \{E_1, \dots E_r \}$ be a simplex of $\tilde{\Delta}= \orifl \tilde{X}$ and set $\EE = \EE(a)=\{E_i,E_i^\perp \mid 1 \leq i \leq r\}$. Then $\EE \cap \mathcal{U}_n = \emptyset$, $Y_{\EE}(V) = X_{\EE}(V)$, $e_n \leq 2, e_{n-1} = e_{n+1}=0$ and $e_i \leq 1$ for all other $1 \leq i \leq 2n-1$.
\end{remark}

\begin{corollary} \label{app: dn case cor}
Let $\tilde{\Delta}$ be a building of type $D_n$ over the field $K$, and $a \in \tilde{\Delta}$ be a simplex. 
If $\lvert K \rvert \geq 2^{2n-1}$, then $\tilde{\Delta}^0(a)$ has a NAC with
$\Rad_j(\Cone_{\tilde{\Delta}^0(a)}) \leq 2 \mathcal{R}(n-1), 0 \leq j \leq 1$, where $\mathcal{R}(n)$ does not depend on $K$.
\end{corollary}
\begin{proof}
By \Cref{thm: class for D_n case}, we have that $Y_{\EE(a)}(V)$ has NAC with radius $\leq \mathcal{R}(n-1)$. By \Cref{rmk: connection to op compl D_n}, we have $X_{\EE(a)}(V) = Y_{\EE(a)}(V)$. Thus \Cref{prop: quantitative subdivision D_n} yields a NAC for $\tilde{\Delta}^0(a) = \tilde{T}_{\EE(a)}(V)$ with cone radius bounded by $2\mathcal{R}(n-1)$.
\end{proof}

\bibliographystyle{alpha}
\bibliography{bibl}

\begin{thebibliography}{GdPVB25}

\bibitem[AB08]{AB2008}
P.~Abramenko and K.~S. Brown.
\newblock {\em Buildings}, volume 248 of {\em Graduate Texts in Mathematics}.
\newblock Springer, New York, 2008.
\newblock Theory and applications.

\bibitem[Abr96]{Abr}
Peter Abramenko.
\newblock {\em Twin buildings and applications to {S}-arithmetic groups},
  volume 1641 of {\em Lecture Notes in Mathematics}.
\newblock Springer-Verlag, Berlin, 1996.

\bibitem[ACR12]{ApteCR}
Himanee Apte, Pratyusha Chattopadhyay, and Ravi~A. Rao.
\newblock A local global theorem for extended ideals.
\newblock {\em J. Ramanujan Math. Soc.}, 27(1):1--20, 2012.

\bibitem[AH93]{AbelsHolz}
H.~Abels and S.~Holz.
\newblock Higher generation by subgroups.
\newblock {\em J. Algebra}, 160(2):310--341, 1993.

\bibitem[BLM24]{BLM}
Mitali Bafna, Noam Lifshitz, and Dor Minzer.
\newblock Constant degree direct product testers with small soundness.
\newblock In {\em 2024 {IEEE} 65th {A}nnual {S}ymposium on {F}oundations of
  {C}omputer {S}cience---{FOCS} 2024}, pages 862--869. IEEE Computer Soc., Los
  Alamitos, CA, [2024] \copyright 2024.

\bibitem[BM24]{BM}
Mitali Bafna and Dor Minzer.
\newblock Characterizing direct product testing via coboundary expansion.
\newblock In {\em S{TOC}'24---{P}roceedings of the 56th {A}nnual {ACM}
  {S}ymposium on {T}heory of {C}omputing}, pages 1978--1989. ACM, New York,
  [2024] \copyright 2024.

\bibitem[BMVY25]{BMV}
Mitali Bafna, Dor Minzer, Nikhil Vyas, and Zhiwei Yun.
\newblock Quasi-linear size {PCP}s with small soundness from {HDX}.
\newblock In {\em S{TOC}'25---{P}roceedings of the 57th {A}nnual {ACM}
  {S}ymposium on {T}heory of {C}omputing}, pages 45--53. ACM, New York, [2025]
  \copyright 2025.

\bibitem[Bou08]{bourbaki2008lie}
N.~Bourbaki.
\newblock {\em Lie Groups and Lie Algebras: Chapters 4-6}.
\newblock Number dln. 4-6 in Elements de mathematique [series]. Springer Berlin
  Heidelberg, 2008.

\bibitem[Car89]{carter1989simple}
R.~W. Carter.
\newblock {\em Simple groups of Lie type}, volume~22.
\newblock John Wiley \& Sons, 1989.

\bibitem[DD24a]{DD-cosys}
Yotam Dikstein and Irit Dinur.
\newblock {Coboundary and Cosystolic Expansion Without Dependence on Dimension
  or Degree}.
\newblock In Amit Kumar and Noga Ron-Zewi, editors, {\em Approximation,
  Randomization, and Combinatorial Optimization. Algorithms and Techniques
  (APPROX/RANDOM 2024)}, volume 317 of {\em Leibniz International Proceedings
  in Informatics (LIPIcs)}, pages 62:1--62:24, Dagstuhl, Germany, 2024. Schloss
  Dagstuhl -- Leibniz-Zentrum f{\"u}r Informatik.

\bibitem[DD24b]{DDcovers}
Irit Dinur and Yotam Dikstein.
\newblock Agreement theorems for high dimensional expanders in the small
  soundness regime: the role of covers.
\newblock \url{https://arxiv.org/abs/2308.09582}, 2024.

\bibitem[DD24c]{DDswap}
Irit Dinur and Yotam Dikstein.
\newblock Swap cosystolic expansion.
\newblock \url{https://arxiv.org/abs/2312.15325}, 2024.

\bibitem[DDL24]{DDL}
Yotam Dikstein, Irit Dinur, and Alexander Lubotzky.
\newblock Low acceptance agreement tests via bounded-degree symplectic {HDX}s.
\newblock In {\em 2024 {IEEE} 65th {A}nnual {S}ymposium on {F}oundations of
  {C}omputer {S}cience---{FOCS} 2024}, pages 826--861. IEEE Computer Soc., Los
  Alamitos, CA, [2024] \copyright 2024.

\bibitem[DM22]{DM}
Irit Dinur and Roy Meshulam.
\newblock Near coverings and cosystolic expansion.
\newblock {\em Arch. Math. (Basel)}, 118(5):549--561, 2022.

\bibitem[EK24]{EK}
Shai Evra and Tali Kaufman.
\newblock Bounded degree cosystolic expanders of every dimension.
\newblock {\em J. Amer. Math. Soc.}, 37(1):39--68, 2024.

\bibitem[GdPVB25]{hdxfromkms}
Laura Grave~de Peralta and Inga Valentiner-Branth.
\newblock High-dimensional expanders from {K}ac-{M}oody-{S}teinberg groups.
\newblock {\em European J. Combin.}, 126:Paper No. 104131, 26, 2025.

\bibitem[Hum90]{HumphreyReflGrp}
James~E. Humphreys.
\newblock {\em Reflection groups and {C}oxeter groups}, volume~29 of {\em
  Cambridge Studies in Advanced Mathematics}.
\newblock Cambridge University Press, Cambridge, 1990.

\bibitem[Kac83]{KacInfDimLieAlg}
V.~G. Kac.
\newblock {\em Infinite-dimensional {L}ie algebras}, volume~44 of {\em Progress
  in Mathematics}.
\newblock Birkh\"{a}user Boston, Inc., Boston, MA, 1983.
\newblock An introduction.

\bibitem[KKL14]{KKL}
Tali Kaufman, David Kazhdan, and Alexander Lubotzky.
\newblock Ramanujan complexes and bounded degree topological expanders.
\newblock In {\em 55th {A}nnual {IEEE} {S}ymposium on {F}oundations of
  {C}omputer {S}cience---{FOCS} 2014}, pages 484--493. IEEE Computer Soc., Los
  Alamitos, CA, 2014.

\bibitem[KM21]{KM}
Tali Kaufman and David Mass.
\newblock {Unique-Neighbor-Like Expansion and Group-Independent Cosystolic
  Expansion}.
\newblock In Hee-Kap Ahn and Kunihiko Sadakane, editors, {\em 32nd
  International Symposium on Algorithms and Computation (ISAAC 2021)}, volume
  212 of {\em Leibniz International Proceedings in Informatics (LIPIcs)}, pages
  56:1--56:17, Dagstuhl, Germany, 2021. Schloss Dagstuhl -- Leibniz-Zentrum
  f{\"u}r Informatik.

\bibitem[KO23]{KO2023high}
T.~Kaufman and I.~Oppenheim.
\newblock High dimensional expanders and coset geometries.
\newblock {\em European J. Combin.}, 111:Paper No. 103696, 31, 2023.
\newblock With a preface by Alexander Lubotzky.

\bibitem[KOW24]{KO24cbeofcoco}
Tali Kaufman, Izhar Oppenheim, and Shmuel Weinberger.
\newblock Coboundary expansion of coset complexes, 2024.

\bibitem[Mar18]{marquis2018introduction}
T.~Marquis.
\newblock {\em An introduction to {K}ac-{M}oody groups over fields}.
\newblock EMS Textbooks in Mathematics. European Mathematical Society (EMS),
  Z\"{u}rich, 2018.

\bibitem[ML63]{MacLane63}
Saunders Mac~Lane.
\newblock {\em Homology}, volume Band 114 of {\em Die Grundlehren der
  mathematischen Wissenschaften}.
\newblock Springer-Verlag, Berlin-G\"ottingen-Heidelberg; Academic Press, Inc.,
  Publishers, New York, 1963.

\bibitem[Opp18]{OppLocI}
Izhar Oppenheim.
\newblock Local spectral expansion approach to high dimensional expanders
  {P}art {I}: {D}escent of spectral gaps.
\newblock {\em Discrete Comput. Geom.}, 59(2):293--330, 2018.

\bibitem[OVB25]{CSEofKMS}
Izhar Oppenheim and Inga Valentiner-Branth.
\newblock New cosystolic high-dimensional expanders from {KMS} groups.
\newblock \url{https://arxiv.org/abs/2504.05823}, 2025.

\bibitem[Reh75]{rehmann1975prasentationen}
U.~Rehmann.
\newblock Pr{\"a}sentationen von {C}hevalleygruppen {\"u}ber k[t].
\newblock {\em preprint}, 1975.

\bibitem[Tit74]{TitsBuildings}
Jacques Tits.
\newblock {\em Buildings of spherical type and finite {BN}-pairs}, volume Vol.
  386 of {\em Lecture Notes in Mathematics}.
\newblock Springer-Verlag, Berlin-New York, 1974.

\bibitem[Wei13]{WeibelKBook}
Charles~A. Weibel.
\newblock {\em The {$K$}-book}, volume 145 of {\em Graduate Studies in
  Mathematics}.
\newblock American Mathematical Society, Providence, RI, 2013.
\newblock An introduction to algebraic $K$-theory.

\end{thebibliography}
\end{document}